\documentclass[a4paper]{article}

\usepackage{amsfonts}
\usepackage{amsmath}
\usepackage{amssymb}
\usepackage{amscd}
\usepackage{cancel}
\usepackage{mathtools}
\usepackage{here}
\usepackage{lipsum}
\usepackage{lineno}
\usepackage{tikz}
\usepackage[all]{xy}
\usepackage{xcolor}
\usepackage{bbm}
\usepackage{colortbl}
\usepackage{selinput}
\usepackage{hyperref}
\hypersetup{
    colorlinks=true,
    linkcolor=blue,
    filecolor=magenta,
    urlcolor=cyan,
    pdftitle={Overleaf Example},
    pdfpagemode=FullScreen,
    }
\newcommand{\lo}{\Omega}
\newcommand{\tg}{\approx}
\newcommand{\cz}{\,\exists\;}

\newcommand{\jia}{\oplus}
\newcommand{\z}{\mathbb{Z}}
\newcommand{\h}{\mathbb{H}}
\newcommand{\s}{\Sigma}
\newcommand{\hc}{\circ}
\newcommand{\n}{\{ }
\newcommand{\nn}{\} }
\newcommand{\ca}{\eta}
\newcommand{\cg}{\sigma}
\newcommand{\lt}{\varepsilon}
\newcommand{\ch}{ \bar{\nu}  }

\newcommand{\ck}{\zeta}
\newcommand{\cn}{\kappa}
\newcommand{\co}{\rho}

\newcommand{\cs}{\bar{\zeta}}
\newcommand{\ccs}{\bar{\sigma}}
\newcommand{\ct}{\bar{\kappa}}
\newcommand{\kj}{\phi}
\newcommand{\cx}{\delta}
\newcommand{\w}{\omega}
\newcommand{\wq}{\infty}
\newcommand{\af}{\alpha}
\newcommand{\sk}{\mathrm{sk}}
\newcommand{\Ker}{\mathrm{Ker}}
\newcommand{\cok}{\mathrm{Coker}}
\newcommand{\ord}{\mathrm{ord}}
\newcommand{\ind}{\mathrm{Ind}}
\newcommand{\sho}{\;\mathrm{Short}}
\newcommand{\m}{\;\mathrm{mod}\,}
\newcommand{\pa}{\partial}
\newcommand{\lon}{\;\mathrm{Long}}
\newcommand{\x}{\langle}
\newcommand{\xx}{\rangle}
\newcommand{\bz}{\subseteq}

\newcommand{\dyd}{\supseteq}
\newcommand{\xyd}{\subseteq}
\newcommand{\ty}{ \equiv}
\newcommand{\bzd}{\hfill\boxed{}}
\newcommand{\q}{ \Delta}
\newcommand{\y}{\mathbb{H}P^{2}}

\newtheorem{thm0}{Theorem}
\newtheorem{thm}{Theorem}[subsection]
\newtheorem{definition}[subsection]{Definition}
\newtheorem{cor}[thm]{Corollary}
\newtheorem{pro}[thm]{Proposition}
\newtheorem{lem}{Lemma}[subsection]

\newtheorem{rem}[thm]{Remark}

\newenvironment{proof}{\noindent{\bf Proof.}}
{\noindent \ \hfill$\Box$\par}
\begin{document}
\title{\bf{On the Homotopy Groups of the Suspended Quaternionic Projective Plane and  Applications}
\\
}
\author{Juxin Yang\thanks{School of Mathematical Sciences, Hebei Normal University,\,
Shijiazhuang, 050024, P.R. China;
and Yanqi Lake Beijing Institute of Mathematical Sciences and Applications (BIMSA),\;
Beijing, 101408, P.R. China. EMAIL:
yangjuxin@bimsa.com},\;\;Juno Mukai\thanks{Shinshu University,
3-1-1\,Asahi, Matsumoto, Nagano 390-862, Japan.\\
EMAIL:\;jmukai@shinshu-u.ac.jp} \;and Jie Wu\,\thanks{BIMSA,
Beijing, 101408, P.R. China.\;EMAIL:
wujie@bimsa.com}\;$^{\,,}$\thanks{Corresponding author}}
\date{}
\maketitle
{\bf Abstract}
In this paper, we determine  the \textit{2,3}-primary components of the homotopy groups $\pi _{r+k}(\Sigma ^{k}\mathbb{H}P^{2})$ for all $ 7\leq r\leq15$ and all $k\geq0$,   essentially we give the  determinations of the integral  homotopy groups $\pi _{r+k}(\Sigma ^{k}\mathbb{H}P^{2})$ for all $ 0\leq r\leq15$  and all $\;k\geq0$, which  in particular include the unstable ones.\;And we give the applications, including  constructing a suspended generalized Hopf fibration by a Toda bracket localized at 2,    two classification theorems  of  simply-connected  \textit{CW} complexes of  type of   suspended $\h P^{3}$ and some homotopy decompositions of  suspended self smash products localized at 3.
\section{Introduction}

Homotopy groups are situated  in the center of homotopy theory. In addition to the homotopy groups of spheres, the homotopy groups of finite \textit{CW} complexes, have been extensively studied. Experts give many fundamental results,\,such as J.\,H.\,C. Whitehead \cite{Jh1}, A.\,L. Blakers\,-W.\,S. Massey \cite{BlMa1}, I.\,M. James \cite{Ja1}, H. Toda \cite{To1}, M.\,E. Mahowald \cite{Ma1} and so on. Nowadays, the Blakers-Massey theorem, the James theorem  (\cite[Theorem 2.1]{Ja1}), the relative \textit{EHP} sequence (\cite{GW}) and Toda bracket methods are still critical    methods to determine the unstable homotopy groups of  finite \textit{CW} complexes. \\
\indent The projective spaces are core objects in algebraic topology.\,The homotopy theory of  the quaternionic projective spaces has been paid attention by  J.\,F. Adams \cite{Ad1},
D. Sullivan \cite{Su1} and I.\,M. James \cite{James76}. The homotopy groups of the suspended quaternionic projective spaces are concentrated by
A. Liulevicius \cite[1962]{Liulevicius} and the second author \cite{Mukai}, it is worthy mentioning that in the periods of \cite{Liulevicius}, the homotopy groups $\pi_{i}(S^{n})$ determined by H.\,Toda in \cite{Toda} were not well-known, A.\,Liulevicius did not use Toda's results but  used the  Adams spectral sequence,  determined many stable homotopy groups of the projective spaces.\\
\indent Let $\mathbb{H}$$P^{n}$ denote the quaternionic projective space of dimension $n$, it is the orbit space $(\h^{n+1}-\n0\nn)/(\h-\n0\nn)=S^{4n+3}/S^{3}$, and it is a \textit{CW} complex of type  $S^{4}\cup e^{8}\cup \cdots \cup e^{4n}.\,$ In this paper, we mainly determine the homotopy groups
$\pi _{r+k}((\Sigma ^{k}\mathbb{H}P^{2})_{(p)})$\;($p\in\n2,3\nn$) for all $ 7\leq r\leq15$ and all $k\geq0$,\;equivalently speaking, determine  the 2,\,3-primary components of these groups. It is well-known that $\pi_{8}(S^{5})\tg\z/24$ which tells us for any prime $p \notin \n2,3\nn$, $\s\h P^{2}\simeq S^{5}\vee S^{9}$  localized at $p$. Further we know if a  prime $p\notin \n2,3\nn$, the homotopy groups $\pi _{r+k}(\Sigma ^{k}\mathbb{H}P^{2})$ \;($ r\leq15$ and $\;k\geq0$)  localized at $p$ are  essentially known  in the  history due to Hilton-Milnor's formula and Toda. And for $ r\leq6,$   the groups are just homotopy groups of spheres which are also known. While for the groups  $\pi _{r+k}(\Sigma ^{k}\mathbb{H}P^{2})$\;($ 7\leq r\leq15$ ) localized at 2 or 3, especially for the unstable homotopy groups, generally saying, the situations  are more and more mysterious as  $r$ grows. In \cite{Gray},\;Brayton\,Gray gives a method to decide the homotopy type of  the homotopy fibre of the pinch map by his relative James construction, in some sense,  Gray's method (our Proposition\,\ref{grdl}) gives another view to understand the relative homotopy group method, namely, the James theorem. Gray's method is one of our fundamental methods to determine the homotopy groups of the two-cell complexes.\\\indent Since there're well-known isomorphisms $\pi_{n}(\y)\tg\pi_{n}(S^{11})\jia\pi_{n-1}(S^{3})$, so we will omit all the proofs for the determinations of
$\pi _{r+k}(\Sigma ^{k}\mathbb{H}P^{2})$  localized at 2 or 3 in the case $k=0$.\\
 \indent  Here are our main theorems.
\begin{thm0}
The \textit{2,3-primary components} of $\pi_{r+k}(\s^{k}\y)$\, ($ 7\leq r\leq15,\,$ $\;k\geq0$)  are summarized in the following table.\end{thm0}\;
\textcolor{white}{66666666666666} The \textit{2,3-primary components} of $\pi_{r+k}(\s^{k}\y)$\\\resizebox{12.04cm}{2.39cm}{
 \begin{tabular}{|c|c|c|c|c|c|c|c|c|c|c|c|c| }
\hline 
  $r\backslash\;k$   &0 &1&2&3 &4&5&6&7&8 &9&10&$\cdots$\\
  \hline
  7& 4+3& 0 & $\mathbf{
  -}$& & &  & & & &  & & \\
     \hline
   8&     2& $\wq$ & $\wq$& $\mathbf{
  -}$& & & &&& &&\\
    \hline
       9&  2 &0&$\wq$&2&$\mathbf{
  -}$&&&& &&&\\
    \hline
     10&  3 &0&2&2&2&$\mathbf{
  -}$&&& &&&\\
    \hline
     &     &&8+2&8+4&$\wq$+8&16+4&&& &&&\\        11    &$\wq$+3 &8+9&+9&+9&+4+9&+9& $\mathbf{
  -}$&&& && \\\hline
     &   &&&&&&&& &&&\\
    12 &$2^{2}$&$\wq+2$&$2^{2}$&$2^{3}$&$2^{4}$&$2^{3}$&$2^{2}$&$\mathbf{
  -}$& &&&

     \\
    \hline
      &   &&&&&&& &&&&\\
        13&   $2^{3}$ &$2^{3}$&$\wq+2^{2}$&$2^{3}$&$2^{4}$&$2^{3}$&$\wq+2^{2}$&$2^{2}$&$\mathbf{
  -}$&&&
      \\
    \hline
     &8+4&&&$\wq+2$&&&&& &&&\\

     14 &+2+$3^{2}$ &$2^{2}$+3&$2^{2}$+3&+3&$2^{2}$+3&$2^{2}$+3&$2^{2}$+3&$2^{2}$+3&2+3& $\mathbf{
  -}$&&

     \\
    \hline
     &4+$2^{2}$ &16+8&16+$2^{2}$&16+$2^{2}$&$\wq$+16&&&&$\wq$+128& &&\\
      15 &+3 &+9+3&+9&+27&+2+27&32+27&64+27&128+27&+27&128+27&$\mathbf{
  -}$&

     \\
    \hline
\end{tabular}}
\textcolor{white}{6}\\\\\\\textit{
In this table, $n$ indicates $\z/n$,  $n^{m}$ indicates $(\z/n)^{m}$,\,(for  positive $n$), $\wq$ indicates $\z$, $0$ indicates the trivial group, and $+$ indicates the direct sum. For $k\geq r-6$, $\pi_{r+k}(\s^{k}\y)$ are in the stable range, the 2,3-components of the groups $\pi_{r+k}(\s^{k}\y)$\;($k= (r-6)+1$) are just denoted by  bar ``$\mathbf{
  -}$", and the others in the stable range isomorphic to them respectively are omitted.}
\label{tab}

\textcolor{white}{.}

\begin{thm0}
Let $h: S^{11}\rightarrow \y$ be the homotopy class of the Hopf fibration whose homotopy cofibre is $\h P^{3}$. 
 After localization at 2, for some   odd integer $t$,  $\s h$ is contained in the Toda bracket
 $$ \n   j_{1},\nu_{5},2t\nu_{8}\nn
,$$
  where $S^{5}\stackrel{j_{1}}\longrightarrow \s\y$ is the inclusion, $\nu_{n}\in\pi_{3+n}(S^{n})$ ($n\geq4$) are the Hopf classes; what's more,  $\s h$ generates                        $\pi_{12}(\s\y)\tg\z/8$, but $\s^{\wq}h$ cannot generate any direct summand of $ \pi_{11}^{S}(\y)\tg\z/16\jia\z/4.$

\end{thm0}
\begin{thm0}

\label{0fldl} After localization at 3, suppose $k\geq1$ but $k\neq4$, then, up to homotopy, the $k$-fold  suspension of the simply connected homology $3$-dimensional  quaternionic projective spaces, that is,  \textit{CW} complexes of type  $$S^{4+k}\cup e^{8+k}\cup e^{12+k},$$ can be classified as the following, \\
\indent$ \Sigma ^{k}\mathbb{H}P^{3}, \;\quad \Sigma ^{k}\mathbb{H}P^{2}\cup_{3\Sigma^{k}h}e^{12+k},\;\quad\Sigma ^{k}\mathbb{H}P^{2}\vee S^{12+k},$\\
\indent$\Sigma^{k-1}A\vee S^{8+k},\; S^{4+k}\vee\Sigma^{4+k}\mathbb{H}P^{2},$ \;$(S^{4+k}\vee S^{8+k})\cup _{c_{k}} e^{12+k}\;$  \\and\;\;$S^{4+k}\vee S^{8+k} \vee S^{12+k}.$
\end{thm0}

\begin{thm0}
After localization at 3, suppose $Z$ is a simply-connected \textit{CW} complex, $k\geq1$ but $k\neq4$, if  $$ \widetilde{H}_{*}(Z;\z/3)\tg   \widetilde{H}_{*}( \s^{k}\h P^{3};\z/3)$$  as  Steenrod modules, then
$$\;\quad Z\simeq \s^{k}\h P^{3}.$$
\end{thm0}

\begin{thm0} After localization at 3, there exist homotopy equivalences,
$$\Sigma\mathbb{H}P^{2}\wedge\mathbb{H}P^{2}\,\simeq S^{13}\vee\Sigma^{5}\mathbb{H}P^{3},$$$$\Sigma\mathbb{H}P^{3}\wedge\mathbb{H}P^{3}\,\simeq\, \Sigma^{9}\mathbb{H}P^{3}\vee Y,$$

\noindent where $Y$ is a 6-cell  \textit{CW} complex and $\sk_{13}(Y)=\Sigma^{5}\mathbb{H}P^{2}$.
\end{thm0}
\noindent {\bf Acknowledgement}.
This work  is supported in part by Natural Science Foundation of China (NSFC\,grant\,no.\,11971144), High-level Scientific Research Foundation of Hebei Province and the start-up research fund from Yanqi Lake Beijing Institute of Mathematical Sciences and Applications.
\tableofcontents
\section{Preliminaries}

\subsection{Notations and some fundamental facts}
\indent\quad In this paper, all spaces, maps, homotopy classes are pointed,  basepoints and constant maps are denoted by $*$, and  homotopy classes of constant maps are denoted by $0$.
 If we  take the $p$-localization, we always use the original symbols  of  the spaces, maps and homotopy classes   to denote them after localization at $p$,\;and $\widetilde{H}_{\ast}(-)$  denotes the mod $p$ reduced homology.\;The homotopy fibre is called the fibre for short. For a map or a homotopy class $f$, we use $C_{f}$ to denote the homotopy cofibre of $f$, and the  homotopy cofibre is called the cofibre for short.\\
\indent   For a non-negative integer  $m$, let $\z/m$ denote $\z/m\z$. For a prime $p$, let $\z_{(p)}=\n \frac{a}{b} \;|$ $\,a, b\in\z $ are coprime and $p \nmid b$$\nn$, that is, the group or the ring of the $ p$-local integers;
  for a $\z_{(p)}$-module $A$  of form $\z_{(p)}$ or $\z/p^{k}$, we use $G=A\n\mathbbm{x}\nn$ to  denote a $\z_{(p)}$-module $G$ which is ismorphic to $A$ and generated by $\mathbbm{x}$, for example, $G=\z/4\n \mathbbm{x}\nn$ stands for $G\tg\z/4$ and $G$ is generated by
  $\mathbbm{x}$; and $(\z/m)^{k}$ denotes the direct sum of $k$-copies of $\z/m$,  we use $\jia$ to denote both the internal direct sum and  the external direct sum. And $\mathrm{ord}(\mathbbm{x})$ denotes the order of the element $\mathbbm{x}$ of a group. \\
\indent And after localization at a prime $p$, in a $\z_{(p)}$-module, $\x a_{1},a_{2}, \cdots,a_{n}\xx$ denotes the  $\z_{(p)}$-submodule generated by $a_{1},a_{2}, \cdots, a_{n}$. If we don't take the localization, in an abelian group,  $\x a_{1}, \;a_{2}, \cdots,\;a_{n}\xx$ denotes the $\z$-submodule generated by $a_{1},a_{2}, \cdots,a_{n}$, we know $\x\alpha,\beta, \gamma  \xx$ also denotes the stable Toda bracket, while these two meanings of $\x-,\;-,\;-\xx$ are always easy to distinguish. We use $\z_{+}$ to denote the set of  positive  integers. Suppose $p$ is a prime,   for a finitely generated  abelian group $G$, 
the \textit{p}-primary component of $G$ is $$G/\x g\in G\;|\; \ord(g)=q^{n} \;\text{for some prime\;} q\neq p  \;\text{and some}\;  n\in \z_{+}   \xx,$$ the so-called \textit{2,3}-primary component of $G$ is $$G/\x g\in G\;|\; \ord(g)=q^{n} \;\text{for some prime\;} q\notin\n 2, 3 \nn \;\text{and some}\; n\in \z_{+}   \xx.$$
Our result gives the homoptopy groups after localization at $p\in\n2,3\nn$, they are naturally corresponding to the \textit{2,3}-primary components these groups. \\
\indent
For  a commutative unitary ring  $R$,  let$$R\n a_{1}, \;a_{2}, \cdots,\;a_{n},\cdots,           \;        \nn,\; \;( |a_{k}|=j_{k})$$  denote the graded free $R$-module with bases  $ a_{1}, \;a_{2}, \cdots,\;a_{n},\cdots,           \;  $and    $a_{k}$ is of degree $j_{k}$. And we often use the obvious symbol  $X \stackrel{\xyd}\longrightarrow Y$ to denote the inclusion map  from $X$ to $Y$, where $X$,\,$Y$ are spaces or modules. For a space $X$, let $X^{\wedge n}$ denote the $ n$-fold self smash product of $X$, and  $X^{\wedge 0}=S^{0}$.\\
 \indent  The  symbols of the  generators of $\pi_{n+k}(S^{n})$  we use are all from    \cite{Toda},\;\cite{20STEM},\;\cite{2122STEM},\;\cite{2324STEM},\, \cite{Oda}\,and \cite{32STEM}, and mainly from  \cite{Toda}, the only two differences are: we denote Toda's  $E$ by $\s$,   the  suspension functor, and we denote Toda's  $\q$ by $P$, where $\q$ is the boundary homomorphism of   $EHP$ sequence.\;And by abuse of notation, sometimes we use the  same symbol to denote a map and its homotopy class. And for a \textit{CW} complex $ Z=Y\cup _{g}e^{m+1}$ where $Y=\sk_{m}(Z)$, we  denote the mapping\, cylinder $M_{g}$ by $M_{Y}$ to indicate $M_{Y} \simeq Y$, although  $M_{Y}$ and $M_{Y}/Y$  do not only depend on $Y$ up to homeomorphism, readers will see such  notations of the mapping cylinders will  make many benefits in this paper.\;\\\indent  Let $\alpha\in\pi_{n}(X),\beta\in\pi_{m}(S^{n})$\, where $n\geq2$ and let $k\in\z$,\;commonly and reasonably, $\alpha\beta$ is the the abbreviation of $\alpha\circ\beta$ and $ k\alpha\circ\beta$ is the the abbreviation of $(k\alpha)\circ\beta$;  here we claim that the symbol $k\alpha\beta$ is to denote $k(\alpha\beta)$.\;So it's necessary to  point out that, 
 $$k\alpha\beta\neq k\alpha\circ\beta\; \text{in general,}$$  
of course, $k\alpha\beta= k\alpha\circ\beta$ always holds  if $\beta$ is a suspension, or the codomain of $\alpha$  is $S^{7}$ or   a  group-like \textit{H}-space (\cite[p.\,118]{GW}). However, to write the equation (\cite[Proposition 1.4,\,p.\,11]{Toda})
$\n \af,\beta,\gamma\nn\hc\s\cx=-(\af\hc\n\beta,\gamma,\cx\nn)$ 
briefly, we denote $-(\af\hc\n\beta,\gamma,\cx\nn)$ by $-\af\hc\n\beta,\gamma,\cx\nn$, and  denote
$\n (-\af)\hc \mathbbm{x}\mid \mathbbm{x}\in \n\beta,\gamma,\cx\nn\nn$
 only by $(-\af)\hc\n\beta,\gamma,\cx\nn$, in this sense, we can write the above equation conveniently as $$\n \af,\beta,\gamma\nn\hc\s\cx=-\af\hc\n\beta,\gamma,\cx\nn$$ and  no confusion will arise. \\\indent
By the \textit{$p$-local Whitehead theorem}\,(\cite[Lemma 1.3]{Wilkerson}), we have the following two corollaries.
\begin{cor}\label{ygd}Let $p$ be a prime,  suppose $\alpha:\s X\rightarrow Y$ is a homotopy class, where $\s X$ and $Y$  are simply connected \textit{CW} complexes, then, after localization at $p$,  for any invertible element $t$ in the ring $ \z_{(p)}$,  the cofibres of  $t \alpha$ and   $ \alpha$ have the same homotoy type,\;that is,\;

\centerline{$C_{t\alpha}\simeq C_{\alpha}$}

\end{cor}
\begin{proof}
Denote $\mathrm{id}_{S^{1}}$ by $\iota_{1}$.\,We notice that,\;after localization at $p$,  for any invertible element $t$ in the ring $ \z_{(p)}$,   $t(\mathrm{id}_{\s X})=\mathrm{id}_{X}\wedge t\iota_{1}\in[\s X,\s X]$  has inverse $t^{-1}(\mathrm{id}_{\s X})=\mathrm{id}_{X}\wedge\, t^{-1}\iota_{1}\in[\s X,\s X]$. Then,  this corollary is immediately got by the  $p$-local Whitehead theorem and the naturality of cofibre sequences, we take the homology with $\z_{(p)}$ coefficients.
\end{proof}

\begin{cor}\label{lgd} 
 After localization at a prime $p$, suppose $ m,\,n_{1}, n_{2}\geq2$ are integers,  $j_{k}: S^{n_{k}}\rightarrow S^{n_{1}}\vee S^{n_{2}}$ $\;(k=1,2)$ are the inclusions,\; $\gamma_{k}\in \s\pi_{m-1}(S^{n_{k}-1}),$\,
 $(k=1,2)$. Then,  for any two invertible numbers $t_{1},\;t_{2}$ in the ring $ \z_{(p)}$,  
    the cofibres of  $t_{1} j_{1}\gamma_{1}+ t_{2}j_{2}\gamma_{2}$ and $j_{1}\gamma_{1}+ j_{2}\gamma_{2}$ have the same homotopy type, that is,
  $$C_{t_{1} j_{1}\gamma_{1}+ t_{2}j_{2}\gamma_{2}}\simeq C_{j_{1}\gamma_{1}+ j_{2}\gamma_{2}}.$$
\end{cor}

\begin{proof}

 Following Toda,\; let $ \iota_{k}$ be the homotopy class of the identity map of    $\emph{S}^{\,k}$.
 After localization at $p$, we consider the following  diagram,\\\\
\indent\qquad\qquad\qquad\qquad\xymatrix@C=3.4cm{
  S^{m} \ar[d]_{\iota_{m}} \ar[r]^{t_{1} j_{1}\gamma_{1}+ t_{2}j_{2}\gamma_{2}   } & S^{n_{1}}\vee S^{n_{2}} \ar[d]^{\frac{1}{\;\, t_{1}}\iota_{n_{1}} \vee \frac{1}{\;\, t_{2}}\iota_{n_{2}}  } \\
S^{m}\ar[r]^{ j_{1}\gamma_{1}+ j_{2}\gamma_{2}   } & S^{n_{1}}\vee S^{n_{2}}  }

 \noindent \\ by the given, $\gamma_{1},\gamma_{2}$ are suspensions, then,
\begin{eqnarray}
& 
&\notag(\frac{1}{\;\, t_{1}}\iota_{n_{1}} \vee \frac{1}{\;\, t_{2}}\iota_{n_{2}})(t_{1} j_{1}\gamma_{1}+ t_{2}j_{2}\gamma_{2})  \\\notag& 
=& (\frac{1}{\;\, t_{1}}\iota_{n_{1}} \vee \frac{1}{\;\, t_{2}}\iota_{n_{2}})\hc t_{1} j_{1}\gamma_{1}+(\frac{1}{\;\, t_{1}}\iota_{n_{1}} \vee \frac{1}{\;\, t_{2}}\iota_{n_{2}})\hc t_{2} j_{2}\gamma_{2}   \\\notag&
=& 
\frac{1}{\;\, t_{1}}\iota_{n_{1}} \hc t_{1} \gamma_{1}+\frac{1}{\;\, t_{2}}\iota_{n_{2}}\hc t_{2} \gamma_{2}  \\\notag&
=& 
j_{1}\gamma_{1}+ j_{2}\gamma_{2},
\end{eqnarray}
\noindent so, the diagram is commutative.   Then, this corollary is immediately got by the  $p$-local Whitehead theorem and the naturality of  the  cofibre sequences, checking \textit{mod p} homology.
 
\end{proof}

Some fundamental facts on homotopy groups are in the following.
\begin{rem}
   Let $X$ be an $(n-1)$-connected \textit{CW} complex  and  $\pi_{n}(X)\neq0$, then $\pi_{d+k}(\Sigma^{k} X)$ is  in the stable range\,$\Leftrightarrow d+k<2(n+k)-1$$\Leftrightarrow k\geq d-2n+2$.\; Hence $\pi_{d+k}(\Sigma^{k}\mathbb{H}P^{m})$  $(1\leq m\leq\infty) $  is  in  the stable range $\Leftrightarrow k\geq d-6$. What's more, this relation still holds for  $X=S^{0}$, in this case, we take  $n=0$, we have\;
    $\pi_{d+k}(S^{k})$ is  in the stable range\,$\Leftrightarrow  k\geq d-2n+2=d+2$.\;
\end{rem}\indent\quad
Roughly speaking, the following  (\cite[Corollary\, 5.1]{sggd}) proposition   tells us that   $\pi_{1}(-)$ can commute with colimit.
\begin{pro}\label{clhp1}
Let $I$ be a category with initial object and let $ X : I \rightarrow S_{*}$ be a diagram of
pointed and connected spaces. Then the fundamental group of the homotopy colimit given as a colimit of groups:$$\pi_{1}(\mathrm{hocolim}_{I}\,X)\tg \mathrm{colim}_{I}\pi_{1}(X).$$

\end{pro}

\subsection{On the cell structures of the fibres and cofibres}\label{sce6}

\begin{lem}\label{CW1}(\cite[p.\,237]{Cell})
Let $f: X\rightarrow Y$ be a fibration where $X$ and $Y$ both have homotopy types of \textit{CW} complexes and $Y$ is
path-connected,\; then the fibre of $f$ also has a homotopy type of \textit{CW} complex.\\
\end{lem}
\indent \quad The following lemma is immediately got by \cite[Theorem 2.3.1, p.\,62]{Cell}.
\begin{lem}\label{CW2}

 Let $f: X\rightarrow Y$ be a cellular  map between \textit{CW} complexs,\; then the mapping cylinder $M_{Y}$ and the mapping cone $ C_{f}$ of $f$ are both  \textit{CW} complexes.
 What's more,  $ X,\;Y$ are subcomplexes of $M_{Y}$, and  $Y$ is a subcomplex of $ C_{f}.$\quad $\boxed{}$

\end{lem}
\subsection{Extension problems of $\z_{(p)}$-modules}

\indent \quad We give the following lemma to prevent losing some solutions when we meet the  extension problems. We use $\text{gcd}(a,b)$ to denote the  greatest common factor of  integers $a,\,b.$
\begin{lem}\label{ext}Suppose $p$ is a prime,\begin{itemize}
\item [\rm (1)]
let $ 0\rightarrow A\stackrel{i}\rightarrow B\stackrel{q}\rightarrow  \z/p^{r}  \rightarrow0$  be an exact sequence of $\z_{(p)}$-modules, \;then  such  $\z_{(p)} $-modules  $B$ are given by $$ B\approx (A\jia \z_{(p)})/\langle     ( \zeta(x)   ,- p^{r}x )\mid x\in            \z_{(p)}                       \rangle,$$ where      $\zeta\in \mathrm{Hom} _{\z_{(p)}}( \z_{(p)} ,A  )$; if   the cohomology class        $ [\zeta]\in \mathrm{Ext}_{\z_{(p)}}^{1}(\z/p^{r}, A)$  runs over $\mathrm{Ext}_{\z_{(p)}}^{1}(\z/p^{r}, A)$, then, up to isomorphism,  the formula  gives all such   $\z_{(p)} $-modules  $B$  satisfying the exact sequence above.
 \item [\rm (2)]
let $ m,n\in \z_{+}$ and $t=\mathrm{min}\n m,n\nn$, for   the exact sequence  $\z_{(p)} $-modules, \\\\\indent \qquad\qquad  \qquad   \xymatrix@C=0.3cm{
  0 \ar[r] & \z_{(p)}\jia\z/p^{m} \ar[rr] && B\ar[rr]&& \z/p^{n} \ar[r] & 0, }\\\\
 up to isomorphism, all of such  $\z_{(p)}$-modules  $B$ are given by \\\\ \indent \qquad\qquad\qquad $B\approx  \z/p^{m}\jia\z/p^{i}\jia\z_{(p)}, \;\;(0\leq i \leq t-1),$       \\ \indent \qquad\qquad   or\;\, $ B\tg\z/p^{m+i-j}\jia\z/p^{j}\jia\z_{(p)}, \;\;(0\leq j\leq i\leq n$ and $j\leq t)$.
\end{itemize}
\end{lem}
\begin{proof}
    \begin{itemize}
        \item [\rm (1)] We only need to notice that \\\\ \indent \qquad\qquad\qquad\qquad$\cdots\rightarrow0 \rightarrow   \z_{(p)}\stackrel{\times p^{r}} \longrightarrow  \z_{(p)}\stackrel{\text{proj}.\;\;}\longrightarrow \z/p^{r}\rightarrow0$\\\\ is a projective resolution of  $\z/p^{r}$ of $\z_{(p)} $-modules, then (1) of our lemma  is immediately  got by \cite[Theorem 7.30, p.\,425]{Rotman}.

  \item [\rm (2)]We will use some basic techniques in homological algebra, see \cite[p.\,370]{Rotman}. For the  projective resolution of  $\z/p^{n},$\\\\ \indent \qquad\qquad\quad\xymatrix@C=0.3cm{
 \cdots \ar[rr]&&0 \ar[rr] ^{d_{2}}&&   \z_{(p)} \ar[rr]_{\times \;p^{n}}^{d_{1}} &&\z_{(p)}\ar[rr] && \z/p^{n} \ar[r] & 0, }\\\\
the 
deleted projective resolution is 
\\\\ \indent \qquad\qquad\quad\xymatrix@C=0.3cm{
 \cdots \ar[rr]&&0 \ar[rr] ^{d_{2}}&&   \z_{(p)} \ar[rr]_{\times \;p^{n}}^{d_{1}} &&\z_{(p)} \ar[r] & 0, }\\
 
 \noindent applying $\mathrm{Hom}(-, \z_{(p)}\jia\z/p^{m} )$, we have the cochain complex, 
$$0\rightarrow\mathrm{Hom}(  \z_{(p)}, \z_{(p)}\jia\z/p^{m} )\stackrel{d_{1}^{\ast}}\longrightarrow\mathrm{Hom} ( \z_{(p)}, \z_{(p)}\jia\z/p^{m} ) \stackrel{d_{2}^{\ast}}\longrightarrow0 \rightarrow\cdots,  $$
then,  
 \begin{eqnarray}
& 
&\notag\mathrm{Ext}_{\z_{(p)}}^{1}(\z/p^{n},  \z_{(p)}\jia\z/p^{m} ) \\\notag& 
=& \Ker(  d_{2}^{*} )/\mathrm{Im}(  d_{1}^{*} )  \\\notag&
=& \mathrm{Hom} (  \z_{(p)}, \z_{(p)}\jia\z/p^{m} ) /\mathrm{Im} (  d_{1}^{*} )  \\\notag&
\tg& \z/p^{n}\jia\z/p^{t},\; (t=\mathrm{min}\n m,n\nn).
\end{eqnarray}
For $ \varepsilon'_{1}, \varepsilon'_{2}\in \mathrm{Hom}(  \z_{(p)}, \z_{(p)}\jia\z/p^{m} ) $ where  $\varepsilon'_{1}(1)=(1,0),\; \varepsilon'_{2}(1)=(0,1),$ we put $\varepsilon_{1}=\varepsilon'_{1}+\mathrm{Im}(  d_{1}^{*} ),\;\varepsilon_{2}=\varepsilon'_{2}+\mathrm{Im}(  d_{1}^{*} )$. Then,$$\mathrm{Ext} _{\z_{(p)}}^{1}(\z/p^{n},  \z_{(p)}\jia\z/p^{m} )=\z/p^{n}\n \varepsilon_{1}\nn\jia \z/p^{t}\n\varepsilon_{2}  \nn.$$
We notice that if $\mathrm{gcd}(\lambda,p)=\mathrm{gcd}(\mu, p)=1$, 
 ($\lambda,\mu\in\z$), then \\\\\indent \qquad\;\;$((\z_{(p)}\jia\z/p^{m})\jia\z_{(p)})\,/ \,  \langle   ((\lambda p^{i} \varepsilon'_{1}(x), \mu p^{j} \varepsilon'_{2}(x)  )  ,- p^{n}x )  \mid x\in            \z_{(p)}       \rangle\\\indent \qquad=((\z_{(p)}\jia\z/p^{m})\jia\z_{(p)})\,/ \,  \langle  ( ( p^{i} \varepsilon'_{1}(x), p^{j} \varepsilon'_{2}(x)  )  ,- p^{n}x )  \mid x\in \z_{(p)}  \rangle$.\\\\Then by (1) of this lemma,\;  all of such   $\z_{(p)} $-modules  $B$ are given by$$B\tg((\z_{(p)}\jia\z/p^{m})\jia\z_{(p)})/   \langle   ( ( p^{i} \varepsilon'_{1}(x), p^{j} \varepsilon'_{2}(x)  )  ,- p^{n}x )  \;\mid \;x\in            \z_{(p)}        \rangle,$$where $0\leq i\leq n$ and $0\leq j\leq t=\mathrm{min}\n m,n\nn$.\;\\
\indent\qquad For such integers $i,j$:\\
 if $i<j$, equivalently, $ 0\leq i<j\leq t=\mathrm{min}\n m,n\nn$, also equivalently,  $ 0\leq i\leq t-1$ for some $j$. Then,\\\\
$((\z_{(p)}\jia\z/p^{m})\jia\z_{(p)})/  \langle  ( (p^{i} \varepsilon'_{1}(x), p^{j} \varepsilon'_{2}(x))  ,- p^{n}x )  \mid x\in            \z_{(p)}       \rangle$ \\
$\tg$\resizebox{4.1cm}{0.49cm}{
$  \frac{\z_{(p)}\n a,b,c  \nn}{\langle p^{m}b,\;p^{i}a+ p^{j}b-p^{n}c  \rangle   }$ }\\\\
=\resizebox{4.1cm}{0.55cm}{$\frac{\z_{(p)}\n a,\;b,\;c  \nn}{\langle p^{m}b,\;p^{i}(a+  p^{j-i}b-p^{n-i}c)     \rangle }$}\\\\
=\resizebox{5.1cm}{0.6cm}{$\frac{\z_{(p)}\n a+  p^{j-i}b-p^{n-i}c,\;\;b,\;\;c  \nn}{\langle p^{m}b,\;p^{i}(a+  p^{j-i}b-p^{n-i}c)     \rangle }$}\\\\
$\tg \z/p^{m}\jia\z/p^{i}\jia\z_{(p)}$;\\\\ \indent if $i\geq j$,\;equivalently, $ 0\leq j\leq i\leq n$ and $j\leq t$;  let $$b'=p^{i-j}a+  b-p^{n-j}c,$$then, \\\\
$((\z_{(p)}\jia\z/p^{m})\jia\z_{(p)})/  \langle   (( p^{i} \varepsilon'_{1}(x), p^{j} \varepsilon'_{2}(x)  )  ,- p^{n}x )  \mid x\in            \z_{(p)}       \rangle$ \\\\ $\tg$ \resizebox{4.1cm}{0.52cm}{
$ \frac{\z_{(p)}\n a,b,c  \nn}{\langle p^{m}b,\;p^{i}a+  p^{j}b-p^{n}c  \rangle   }$}\\\\=\resizebox{5.1cm}{0.57cm}{
$\frac{\z_{(p)}\n a,b,c  \nn}{\langle p^{m}b,\;p^{j}(p^{i-j}a+  b-p^{n-j}c)     \rangle } $}\\\\
=\resizebox{5.1cm}{0.57cm}{$\frac{\z_{(p)}\n a,\;\;p^{i-j}a+  b-p^{n-j}c,\;\;c  \nn}{\langle p^{m}b,\;\;p^{j}(p^{i-j}a+  b-p^{n-j}c)     \rangle }$}\\
=\resizebox{5.1cm}{0.57cm}{$\frac{\z_{(p)}\n a,\;\;b',\;\;c  \nn}{\langle p^{m+n-j}c-p^{m+i-j}a ,\;\;\;p^{j}b'    \rangle } $}\\\\
=\resizebox{5.1cm}{0.57cm}{$\frac{\z_{(p)}\n a,\;\;b',\;\;c  \nn}{\langle p^{m+i-j}(p^{n-i}c-a) ,\;\;\;p^{j}b'    \rangle} $}\\\\
=\resizebox{5.1cm}{0.57cm}{$\frac{\z_{(p)}\n  p^{n-i}c-a   ,\;\;b',\;\;c  \nn}{\langle p^{m+i-j}(p^{n-i}c-a) ,\;\;\;p^{j}b'    \rangle} $}\\\\
$\tg \z/ p^{m+i-j} \jia\z/p^{j} \jia \z_{(p)}$.\\

Hence the result holds.

\end{itemize}
\end{proof}

\subsection{Toda bracket}\label{sce6}

\textcolor{white}{5\;\,5}Toda bracket is an art of constructing  homotopy liftings   and  homotopy extensions of maps,  it plays a fundamental role on dealing with  composition relations of homotopy classes, it is deeply studied in \cite{Toda}.\;In this paper, we will freely use the well-known properties of Toda brackets shown in \cite[p.\,9-12]{Toda},
especially\;\cite[Proposition 1.2]{Toda}, \cite[Proposition 1.4]{Toda}\,and\, \cite[Proposition 1.6]{Toda}. 
\\\indent To give a detailed proof of our Corollary \ref{y1yl}, we need to   state some definitions about Toda brackets.

 For the sequence of spaces and maps \\\indent \xymatrix@C=0.55cm{
 &&&& L \ar[r]^{f} & M\ar[r]^{g\qquad\qquad} & N,\;\text{where}\; g\circ f\simeq*,  } 
 
\noindent if  a map $\overline{g}: C_{f}\rightarrow N$ makes the following commutative,  \\
\indent\qquad\qquad\qquad\qquad\qquad\xymatrix@C=0.5cm{
  M \ar[d]_{\subseteq} \ar[r]^{g } &  N      \\
  C_{f}\ar[ur]_{ \overline{g} }                     }

\noindent we say  $\overline{g}$  an extension  of $g$  with respect to $f$; a coextension of $ f$  with respect to $g$, that is,  $\widetilde{f}:\s L=C^{+}L\cup C^{-}L\rightarrow C_{g}$, is    defined as the following: we use
the cone functor $C(-)=\,-\,\wedge I$ where $I$ has basepoint 1, define $\widetilde{f}^{+}$  to be the composition  of the following, where $q$ is the quotient map and the map $C^{f}=f\wedge \mathrm{id}_{I}$ is extended over $f$, \\
\xymatrix@C=0.62cm{
  &&&C^{+}L\ar[r]^{ \mathrm{id}} &CL  \ar[r]^{C^{f}} & CM \ar[r]^{\subseteq\quad}_{j_{_{_{CM}}}\;\;} & N\vee CM \ar[r]^{\quad\;q} & C_{g}  ,} \textcolor{white}{.}\\
\noindent define $\widetilde{f}^{-}$  is a composition  of the following form,\\
  \xymatrix@C=0.69cm{
&&&    C^{-}L\ar[r]^{ \mathrm{id}} &CL  \ar[r]^{\overline{gf}} & N \ar[r]^{\subseteq}_{j_{_{_{N}}}} & C_{g}  ,}

\noindent here $\overline{gf}$ is an extension of $gf$ with respect to $\mathrm{id}_{L}$;\\
then $$\widetilde{f}^{+}| _{L}= \widetilde{f}^{-}| _{L},$$ define\,
    $\widetilde{f}:\s L=C^{+}L\cup C^{-}L\longrightarrow  C_{g}$ to be $$\widetilde{f}=\widetilde{f}^{+}\cup\widetilde{f}^{-}.$$ In general, up to homotopy, neither the extension nor the coextension of a map is unique.
To indicate another map, the above map  $\bar{g}$  is also denoted by $\mathrm{ext}_{f}(g)$, and the above map $\widetilde{f}$  is also denoted by $\mathrm{coext}_{g}(f)$. Suppose $\af\in [L,M]$, $\beta\in [M,N]$ and $\beta\af=0$,  an extension of $\af$ with respect to $\beta$  which is denoted by $\mathrm{ext}_{\af}(\beta)$ is defined to be $[\mathrm{ext}_{a}(b)]$  for some choice of $\mathrm{ext}_{a}(b)$, where $a$ can be any representative map of $\af$  and  $b$ can be any representative map of $\beta$; a coextension of $\beta$ with respect to $\af$ which is denoted by $\mathrm{coext}_{\beta}(\af)$ is defined to be $[\mathrm{coext}_{b'}(a')]$ for some choice of $\mathrm{coext}_{b'}(a')$, where $a'$ can be any representative map of $\af$  and  $b'$ can be any representative map of $\beta$. \\\indent
The following is a paraphrase of \cite[Proposition 1.7,\;p.\,13]{Toda}, which is essentially the definition of Toda brackets.
\begin{definition}(\cite{Toda})\label{ec}
For the following sequence of spaces and homotopy classes, where  $\af\hc \s^{n}\beta=\beta\gamma=0,$ 
\[\xymatrix@!C=0.83cm{ &W&  \s^{n}X\; & \s^{n}Y  &\s^{n}Z,& 
       \ar"1,4";"1,3" _{\;\s^{n}\beta_{}} \ar"1,5";"1,4" _{ 
 \s^{n}\gamma}
\ar"1,3";"1,2" _{\af}
     }\]
 The  Toda bracket indexed by $n$ which is denoted by  $\n\af,\s^{n}\beta,\s^{n}\gamma\nn_{n}$ is defined to be the collection of all compositions of form $$(-1)^{n}\mathrm{ext}_{_{\s^{n}\beta}}(\af)\hc \s^{n}\mathrm{coext}_{_{\beta}}(\gamma),$$
 that is, the collection of homotopy classes of form
 $$(-1)^{n}[\mathrm{ext}_{_{\s^{n}b}}(a)\hc \s^{n}\mathrm{coext}_{_{b}}(c)]\;\text{where}\;a\in\af,\;b\in\beta,\;c\in\gamma.$$ 
 \end{definition}

The Toda bracket $ \n     \alpha,\beta,\gamma              \nn  _{0} $ is  denoted by $ \n     \alpha,\beta,\gamma              \nn   $  for short. We shall  give a detailed proof of the following useful corollary, this corollary is well-known to experts majoring in determinations of homotopy groups, however,  so far, there is no clear proof, so in this paper we give a detailed proof of it. \\\indent In the following corollary, we use the same notation to denote a map and its homotopy class if there is no amibiguity.
\begin{cor}\mbox{}\hspace{-8pt}\label{y1yl}\;
For any cofibre sequence
\[\xymatrix@!C=0.43cm{ \cdots&\s X&  C_{f}\; & Y  &X,& 
\ar"1,2";"1,1"
       \ar"1,4";"1,3" _{i} \ar"1,5";"1,4" _{ 
 f}
\ar"1,3";"1,2" _{\;\,p}
     }\]
   \noindent the relation\;$\mathrm{id}_{\s X}\in\n  p,i,f   \nn $ always holds.
\end{cor}

\begin{proof} 
 We regard $p$ as the composition \\\\
 \xymatrix@C=0.9cm{&&
  C_{f} \ar[r]^{pinch\quad} & C_{f}/Y \ar[r]^{\mathrm{id}} & CX/X\ar[r]^{\psi\qquad}_{\cong\qquad} & C^{+}X\cup C^{-}X   }
  
\noindent where $\psi$ is the inverse of 
 $$\kj:C^{+}X\cup C^{-}X \stackrel{\cong} \longrightarrow CX/X,$$  $\kj(a)=*$ for any $a\in C^{+}X$, and $\kj(x\wedge t)= [x\wedge t]$ for any $x\wedge t\in C^{-}X.$ 
 We know $C_{i}$ can  retract to $\s X$ by  pinching $CY$ to $*$,  we take    $$\overline{p}:\;C_{i}=  (Y\cup_{f} CX)\cup_{i} CY\twoheadrightarrow CX/X \stackrel{\psi} \longrightarrow C^{+}X\cup C^{-}X$$  as this retraction. So, $\overline{p}|_{C_{f}}=p$, that is, $\overline{p}$ is  an extension of $p$ with respect to $i$. In the following, we explain that we can construct a map $\widetilde{f}: \s X\rightarrow C_{i}$ which is a coextension of $f$ with respect to $i$, and $\widetilde{f}$ is exactly a right  inverse of $\overline{p}$.\\
\indent  An extension of $i\circ f$ with respect to $\mathrm{id}_{X}$, that is, $\overline{\,i\circ f\,}$,\,  is taken as the composition $C^{-}X\hookrightarrow Y\vee CX\twoheadrightarrow Y\cup_{f} CX$, we take
 $\widetilde{\,f\,}=   \widetilde{\,f\,}^{+}\cup\widetilde{\,f\,}^{-} $, where     $\widetilde{\,f\,}^{-}=j_{_{_{C_{f}}}}\hc\,\overline{i\circ f}$,\;and $\widetilde{\,f\,}^{+}$ can be  anyone  satisfying  its definition, then $\overline{p}\hc\widetilde{\,f\,}^{+}$ is the constant map, successively,\; $\overline{p}\hc\widetilde{f}$ is the composition $$C^{+}X\cup C^{-}X\stackrel{\kj}\longrightarrow CX/X\stackrel{\psi} \longrightarrow  C^{+}X\cup C^{-}X,$$ this composition is $\mathrm{id}_{\s X}$.  Hence we get the result.
 
\end{proof}

We introduce some common terms on dealing with Toda brackets.
\begin{rem}\; 
\begin{itemize}
    \item [\rm(1)]Suppose $n\geq0$,  $\alpha\in [\s^{n}Y,Z      ]$,\;$  \beta\in [X,Y      ]$ and $\gamma\in [W,X]$   satisfy $\alpha\circ   \s^{n}\beta= \beta\circ\gamma=0$,   and $[  \s^{n+1}W,Z  ]    $ is  abelian, then
  $$\mathbbm{x}\in \n     \alpha,\,\s^{n}\beta,\,\s^{n}\gamma               \nn _{n} \;\;\mathrm{mod} \;A$$ 
  means $ \n     \alpha,\,\s^{n}\beta,\,\s^{n}\gamma             \nn _{n} $  contains $\mathbbm{x}$  and it is a coset of  the subgroup $A$, if $A$ is generated by $\n a_{\lambda}\nn _{\lambda\in\Lambda}$, then $\m A$ is also denoted by $\m\,a_{\lambda},\, (\lambda\in\Lambda)$.
  
 For $\mathbbm{a}\in G/H$ where $G$ is an abelian group and $H$ is its subgroup, we say $H$  the indeterminacy of the coset $\mathbbm{a}$, and it's denoted by $\ind(\mathbbm{a})=H$, so, $$\ind(f\hc\n     \alpha,\,\s^{n}\beta,\,\s^{n}\gamma              \nn _{n}\hc \s g)=f\hc\ind(\n     \alpha,\,\s^{n}\beta,\,\s^{n}\gamma              \nn _{n})\hc \s g.$$
 \item [\rm(2)] Suppose $G$ is an abelian group, $\n K_{i} \nn$ is a family of subgroups of $G$,  $\mathbbm{a}_{i}\in G/K_{i}$, and $$g\in\mathbbm{a}_{j_{1}}\xyd \mathbbm{a}_{j_{2}} \xyd \mathbbm{a}_{j_{3}}\xyd \cdots \xyd \mathbbm{a}_{j_{m}} \dyd \mathbbm{a}_{k_{1}} \dyd \mathbbm{a}_{k_{2}} \dyd \cdots \dyd \mathbbm{a}_{k_{n}}\ni g'.$$ We write   $$g\in\mathbbm{a}_{j_{1}}\xyd \mathbbm{a}_{j_{2}} \xyd \mathbbm{a}_{j_{3}}\xyd \cdots \xyd \mathbbm{a}_{j_{m}} \dyd \mathbbm{a}_{k_{1}} \dyd \mathbbm{a}_{k_{2}} \dyd \cdots \dyd \mathbbm{a}_{k_{n}}\ni g'\m A,$$ 
if \;\,$\ind(\mathbbm{a}_{j_{m}})=A$.
That is, $\m A$  corresponds to the largest coset.
  \item [\rm(3)]  Suppose $G$ is an abelian group and $B$ is a subgroup of \,$G$,  $g,g'\in G$,  then
  $$g\ty g'\m B
 \stackrel{\texttt{def.}}\iff 
g-g'\in B;$$ 
  suppose $U\xyd G$,  and $\mathbbm{x}\in G$, then $$\mathbbm{x}\ty U \m B\stackrel{\texttt{def.}}\iff \cz\mathbbm{y}\in U, \;\text{ such that}\,\;\mathbbm{x}\ty\mathbbm{y}\m B;$$
  $$0\ty U
 \stackrel{\texttt{def.}}\iff 
0\in U.$$
 \item [\rm(4)] A Toda bracket 
consisting of only one element is usually   identified with its element.
\end{itemize}

\end{rem}

\section{On the fibre of the pinch map}

\subsection{Gray's relative James construction}
\indent\quad Suppose  $A$ is a closed subspace of  $X$, to parallel to nowaday's familiar notion $J(X)$,  the ordinary James construction, we  use $J(X,A)$ to denote  the relative James construction,  which is denoted  by $(X,A)_{\wq}$ in \cite{Gray} by its founder  \textit{B.\,Gray}.
\begin{pro}\label{grdl}
 Let $X$ be a path-connected  \textit{CW}  complex and $A$ be its path-connected  subcomplex, let $i:A\hookrightarrow X$  be the inclusion, successively   $i:(A,A)\hookrightarrow (X,A)$,\;
then,
\begin{itemize}

\item[\rm(1)]\label{grdl1a}
there exists a fibre sequence\\
\xymatrix@C=1.2cm{
 & J(A) \ar[r]_{J(i)}^{\xyd} & J(X,A) \ar[r]^{} & X\cup_{i} CA\ar[r]^{p} &  \s A  }

\noindent here, $p$ is  the pinch map,    $ J(A)$ is the ordinary James construction, $ J(i)$ is the inclusion   extended over    $i:(A,A)\hookrightarrow (X,A)$.\\
\indent What's more, if  $K\stackrel{f}\longrightarrow L \stackrel{j_{_{L}}} \longrightarrow C_{f}\stackrel{q}\longrightarrow\s K$  \;is a cofibre sequence where  $K$,\;$L$ are path connected  \textit{CW}  complexes,\; then there exists a fibre sequence\\
\xymatrix@C=1.2cm{
  &J(K) \ar[r]_{J(i_{f})}^{\xyd} &  J(M_{L}, K ) \ar[r]^{} & C_{f} \ar[r]^{q} &  \s K,  }

\noindent here,\; $M_{L}$ is the mapping cylinder of $f$,
$ J(i_{f})$  is the inclusion   extended by  the inclusion  $i_{f} :(K,K)\hookrightarrow (M_{L}, K)$.
\item[\rm(2)]\label{grdl2a}   If  $X/A$ is path-connected, and  $\lo\s A$,  $X/A$ are \textit{CW} complexes of finite types, then\\
\centerline{$\widetilde{H}_{\ast}(J(X,A); \mathbbm{k})\approx \widetilde{H}_{\ast}(X;\mathbbm{k})\otimes H_{\ast}(\lo\s A;\mathbbm{k})$,   $ \;(\mathbbm{k}:\;field$).};
\item[\rm(3)]\label{grdl3a} filtration :\\ \indent$J_{0}(X,A)\subseteq J_{1}(X,A)\subseteq J_{2}(X,A)\subseteq\cdots,\quad J(X,A)=\bigcup_{n\geq 0}J_{n}(X,A)$,\\
 where \,$J_{0}(X,A)=*,\;\, J_{1}(X,A)=X$.
\item[\rm(4)]\label{grdl4a} if  $X=\s X' $,\; $A=\s A'$,\;then $J_{2}(X,A)=X \cup _{[1_{X},i]} C(X\wedge A')$, where $[1_{X},i]$  is the Whitehead product.
\item[\rm(5)]\label{grdl5a}  $J_{n}(X,A)/J_{n-1}(X,A)=X\wedge A^{\wedge (n-1)},$\;\\
\indent $ \s J(X,A)\simeq\s \bigvee_{n\geq 0}X\wedge A^{\wedge n}.$
\end{itemize}
\end{pro}
\begin{proof}
We noticed the following fundamental facts:\,Lemma\,\ref{CW1};\;Lemma\,\ref{CW2};  any point of a \textit{CW} complex is non-degenerate; $U\hookrightarrow V$ is a closed cofibration if and only if $(V,U)$ is an NDR pair.
Hence this proposition\, is essentially given by B.Gray.\\
 (1) From  \cite[Theorem 2.11]{Gray}, \cite[Lemma\, 4.1]{Gray} and \cite[Theorem 4.2]{Gray}, just regard $C_{f}$ as $M_{L}\cup CK$.\\(2) From \cite[Theorem 5.1]{Gray}.\\
 (3) From \cite[p.\,498]{Gray}.\\
 (4) From \cite[Corollary\, 5.8]{Gray}. \\
(5) From \cite[p.\,499]{Gray},   \cite[Theorem 5.3]{Gray}.

\end{proof}

\begin{rem}\label{Fk}
By Proposition\,\ref{grdl}(1),  the fibre of the pinch map $p_{k}:\s^{k}\h P^{2}\rightarrow S^{8+k}$  is  $J(M_{S^{4+k}},S^{7+k})$,    where $M_{S^{4+k}}$ is the mapping cylinder of  $\s^{k}f_{\mathrm{Hopf}}$,\; here $f_{\mathrm{Hopf}}:\;S^{7}\rightarrow S^{4} $ is the $\mathrm{Hopf}$ fibration.\end{rem}

\indent  From now on, we fix the symbol $F_{k}$ to denote the fibre of the pinch map\;\;\\$\s^{k}\h P^{2}\stackrel{p_{k}}\longrightarrow S^{8+k}$, that is,\; $F_{k}=J(M_{S^{4+k}},S^{7+k})$; and let   $p=\s^{\wq}p_{1}$.

\begin{cor}\label{bqjg}
Suppose  $2\leq n \leq m$, then for a cofibre sequence \;\\ \centerline{$S^{m}\stackrel{f}\longrightarrow S^{n} \longrightarrow C_{f}\stackrel{q}\longrightarrow S^{m+1}$,}\\ there exists a fibre sequence,
$$J(M_{S^{n}}, S^{m}) \longrightarrow C_{f}\stackrel{q}\longrightarrow S^{m+1},$$ and  $J(M_{S^{n}}, S^{m})$ has a \textit{CW} structure as  
$$J(M_{S^{n}}, S^{m})  \simeq S^{n}\cup e^{n+m}\cup e^{n+2m}\cup e^{n+3m}\cup\cdots,\; (\text{infinitely many cells}).$$
\end{cor}
\begin{proof}

By Proposition\,\ref{grdl}(1),
 for the cofibre sequence\;$S^{m}\stackrel{f}\longrightarrow S^{n} \longrightarrow C_{f}\stackrel{q}\longrightarrow S^{m+1}$, there exists a fibre sequence,\;
$J(M_{S^{n}}, S^{m}) \longrightarrow C_{f}\stackrel{q}\longrightarrow S^{m+1}$.\\
By Proposition\,\ref{grdl}(5),  \\\centerline{$\s J(M_{S^{n}}, S^{m})\simeq \s \bigvee _{r\geq0}   S^{n}\wedge S^{mr}\simeq \bigvee _{r\geq0}   S^{1+n+mr}$\;$_{.}$}
Thus,
\begin{eqnarray}
& 
&\notag\widetilde{H}_{*}(\s J(M_{S^{n}}, S^{m})) \\\notag& 
\tg&  \bigoplus _{r\geq0} \widetilde{H}_{*} ( S^{1+n+mr})   \\\notag&
=& \z\n     x_{1+n} ,\;     x_{1+n+m}      , \; x_{1+n+2m},\;  x_{1+n+3m},\;\cdots   \;            \nn,\;\;\;(| x_{1+n+mr}|=1+n+mr).
\end{eqnarray}
So,$$ \widetilde{H}_{*}( J(M_{S^{n}}, S^{m}))=\z\n     y_{n} ,\;     y_{n+m}      , \; y_{n+2m},\;  y_{n+3m},\;\cdots   \;            \nn,\;\;| y_{n+mr}|=n+mr.$$
Since  $2\leq n \leq m$, as the fibre of $ S^{n}\cup e^{n+m+1}\stackrel{q}\longrightarrow S^{m+1}$,\;$J(M_{S^{n}}, S^{m})$ is simply-connected. So,$$J(M_{S^{n}}, S^{m})  \simeq S^{n}\cup e^{n+m}\cup e^{n+2m}\cup e^{n+3m}\cup\cdots,\; \mathrm{(infinitely\, many\, cells)}.$$

\end{proof}

Similarly, we have,\;
\begin{cor}\label{bqjg1}
 For a cofibre sequence\;\;$X\stackrel{f}\longrightarrow Y \longrightarrow C_{f}\stackrel{q}\longrightarrow \s X$, where $X$ and $Y$ are path-connected \textit{CW} complexes,  then,
there exists an isomorphism of graded groups:\\ $$\widetilde{H}_{*}( J(M_{M_{Y}}, X);\z)  \tg \widetilde{H}_{*}(  \bigvee_{n\geq0} Y\wedge X^{\wedge n};\z).$$
$\bzd$
\end{cor}

We recall Proposition\,\ref{grdl}(1),\; for a map $K\stackrel{f}\longrightarrow L$,\;\;$J( i_{f}):J(K)  \hookrightarrow J(M_{L}, K)$ is the inclusion extended over the inclusion $i_{f}:(K,K)\hookrightarrow (M_{L}, K)$.
\begin{lem}\label{xh2}
  Suppose $2\leq n \leq m$, then,  for  a cofibre sequence $S^{m}\stackrel{f}\longrightarrow S^{n} \longrightarrow C_{f}\stackrel{q}\longrightarrow S^{m+1}\rightarrow\cdots$, there exist   commutative diagrams for each $k\geq0$,  where  $ j_{_{S^{n+1}}}: S^{n+1} \longrightarrow \s J(M_{S^{n}}    ,S^{m} )$ is the inclusion.  \\

\indent\qquad\qquad\qquad\xymatrix@C=2.4cm{
   \pi_{k}(J(S^{m}))  \ar[d]_{\tg} \ar[r]^{(J( i_{f}))_{*}   } &   \pi_{k}(    J(M_{S^{n}}    ,S^{m} ) )\ar[d]^{\s} \\
 \pi_{k+1}(S^{m+1})\ar[d]_{ (\s f)_{*}} \ar[r]^{       } & \pi_{k+1}(   \s J(M_{S^{n}}   ,S^{m} )) \ar[d]^{\mathrm{id}} \\
 \pi_{k+1}(S^{n+1}) \ar[r]^{( j_{_{S^{n+1}}})_{*}} &  \pi_{k+1}(   \s J(M_{S^{n}}   ,S^{m} ))   }

\end{lem}

\begin{proof}
There's commutative diagram,\\\indent\qquad\qquad\quad\qquad\resizebox{5.6cm}{0.4cm}{

\xymatrix@C=2.4cm{
 J(S^{m})  \ar[d]_{\mathrm{id}} \ar[r]^{J( i_{f})  }_{\xyd} &     J(M_{S^{n}}    ,S^{m}  )\ar[d]^{\xyd} \\
 J(S^{m})\ar[d]_{ \xyd} \ar[r]^{        \xyd  } &    J( J(M_{S^{n}}   ,S^{m} )) \ar[d]^{\mathrm{id}} \\
 J(M_{S^{n}}) \ar[r]^{ \xyd} &  J( J(M_{S^{n}}   ,S^{m} )).   }}\\

 \noindent

 Since in the homotopy category of pointed path connected \textit{CW} complexes, there exists isomorphism of functors $J(-)\cong \lo\s (-)$, and  $S^{n}\hookrightarrow M_{S^{n}}$ is a homotopy equivalence,
 and \;$  J(M_{S^{n}})\hookrightarrow J( J(M_{S^{n}}   ,S^{m} ))  $ is extended by  $j_{M_{S^{n}}}\;: M_{S^{n}}\hookrightarrow J(M_{S^{n}}   ,S^{m} )$, let $j_{_{_{S^{n}}}}:\;S^{n}\rightarrow J(M_{S^{n}}   ,S^{m} )$ \;be the inclusion,  \\
then we have the homotopy-commutative diagram, \\

\indent\qquad\qquad\quad\qquad\resizebox{5.6cm}{0.4cm}{\xymatrix@C=2.4cm{
 \lo S^{m+1}  \ar[d]_{\mathrm{id}} \ar[r]^{ \lo\s i_{f} }_{} &     J(M_{S^{n}}    ,S^{m} ) )\ar[d]^{\lo\s(\mathrm{id})} \\
 \lo S^{m+1}\ar[d]_{ \lo\s f} \ar[r]^{          } &    \lo\s J(M_{S^{n}}   ,S^{m} ) \ar[d]^{\mathrm{id}} \\
 \lo S^{n+1} \ar[r]^{ \lo\s  j_{_{_{S^{n}}}} } &  \lo\s J(M_{S^{n}}   ,S^{m} ).   }}\\\\

  \noindent  applying  the functor $\pi_{k}(-)$ to this diagram, we get the result.
\end{proof}

\begin{cor}\label{j2} Let $m\geq2$ be an integer, and let $Y$ be a \textit{CW} complex satisfying the three conditions:\\
\indent$\mathrm{(i)}$ $Y$ is simply-connected and not contractible;\\
\indent$\mathrm{(ii)}$  $H_{k}(Y;\z)$ is  finitely generated for each $k\geq0$, and the cell structure of $Y$ is minimal, that is, the cell structure of $Y$ is  consistent with its homology;\\
\indent$\mathrm{(iii)} \; \mathrm{dim}(Y)\leq m.$\\
Denote the dimension of
the bottom cell(s) of  $Y$ by  $r$, namely, $$r=\mathrm{min}\n q\in\z_{+} \;|\;H_{q}(Y;\z)\neq0\nn,$$ and set
$Z= Y\cup_{f} e ^{m+1}$ for some  $f$.  For the  fibre $J( M_{Y},S^{m})$ of $Z\stackrel{pinch}\longrightarrow S^{m+1},$ we have $$\mathrm{sk}_{j} (J( M_{Y},S^{m}))=J_{t}( M_{Y},S^{m}),\;\;j=m(t-1)+r-1.$$
\end{cor}

\begin{proof}

By the first one of  Proposition\,\ref{grdl}\,(5), after checking homology, we notice the steps of attaching higher cells by introducing on the dimensions of the skeletons give an isomorphism of functors which are the direct diagrams to define the colimits  $\mathop{\mathrm{colim}}\limits_{n}J_{n}( M_{Y},S^{m})$ and $\mathop{\mathrm{colim}}\limits_{n}\sk_{n}(J( M_{Y},S^{m}))$. Isomorphic functors have isomorphic colimits, hence the result holds.

\end{proof}

The following lemma is   well-known to the experts. For completeness, we give a detailed proof.

\begin{lem}\label{wn}

After localization at 2, \;$$\s\y\simeq S^{5}\cup_{\nu_{5}}e^{9}.$$
After localization at 3,\;$$\s\y\simeq S^{5}\cup_{\alpha_{1}( 5)}e^{9}.$$

\end{lem}

\begin{proof}
We know $\y$ is the cofibre of the map  $f_{\mathrm{Hopf}}:\;S^{7}\rightarrow S^{4}$,  and by checking cohomology ring we know $f_{\mathrm{Hopf}}$ has the Hopf invariant 1,  in the following the homotopy class of   $f_{\mathrm{Hopf}}$ are still denoted by  $f_{\mathrm{Hopf}}$ for simplicity.\\
\indent After localization at 2,  $f_{\mathrm{Hopf}}$ has the Hopf invariant 1,   by \cite[Proposition 8.1,\,p.\,82]{Toda},   $$ f_{\mathrm{Hopf}}\ty \nu_{4}\m 2\nu_{4},\s\nu',$$  thus $$ \s f_{\mathrm{Hopf}}\ty \nu_{5}\m 2\nu_{5},$$   by Corollary\,\ref{ygd} we have,  $ C_{\nu_{5}}\simeq C_{t\nu_{5}}$ for any odd $t$,\, we get the result of the first part.\\
\indent After localization at 3,  we know $\widetilde{H}_{*} (\s\y)=\z/3\n a,\,b  \nn,\;  |a|=5,\;|b|=9$,  with  Steenrod operation  $ P^{1}_{*}(b)=a$   which can be easily got by observing the cohomology ring  $\widetilde{H}^{*} (\h P^{2})$.
Since  $\pi_{8}(S^{5})=\z/3\n \alpha_{1}( 5) \nn$,\;(\cite{Toda}), and  $ P^{1}_{*}(b)=a$ suggests   $\s\y\;\cancel{\simeq} \;S^{5}\vee S^{9}$,  so
$\s f_{\mathrm{Hopf}}=\pm\alpha_{1}( 5)$, by Corollary\,\ref{ygd},  $\s\y\simeq S^{5}\cup_{\alpha_{1}( 5)}e^{9}$.

\end{proof}
\begin{rem} For the proof of Lemma\,\ref{wn}, we point out such a fundamental fact. That is, on  spheres of dimensions 1,3 or 7, for the generalized Hopf invariants defined by Toda,  the generalized Hopf invariants defined by G.W Whitehead, and the original Hopf invariants defined by Hopf,  the first two  correspond  to each other up to sign, and the first two are truly generalized the third. On this fact, see \cite[Theorem 8.17, p.\,540]{GW},  \cite[p.\,60]{Toda1952}  and \cite[Proposition 8.1]{Toda} for more details.

\end{rem}

By Proposition\,\ref{grdl}(4),   Corollary\,\ref{j2} and Lemma\,\ref{wn},\;  we have the following lemma.
\begin{lem}\label{sk0}

Suppose $k\in\z_{+}$,\, then,\;up to homotopy,\\ after localization at 2, $$\mathrm{sk}_{18+3k}( F_{k})=S^{4+k}\cup_{f_{k}} e^{11+2k},$$where  $ f_{k}=[\iota_{4+k},\nu_{4+k}]$ is the  Whitehead product;\\ after localization at 3, $$\mathrm{sk}_{18+3k}( F_{k})=S^{4+k}\cup_{g_{k}} e^{11+2k},$$ where $ g_{k}=[\iota_{4+k},\alpha_{1}( 4+k)]$ is  the  Whitehead product, and \\  $\alpha_{1}( 4+k)\in \pi_{7+k}(S^{4+k})\tg\z/3$ is a generator.$\hfill\boxed{}$

\end{lem}

\begin{lem}\mbox{}\hspace{-8pt}\label{wgs}\;\;\cite[Theorem 8.18, p.\,484]{GW}

 Suppose\;$\alpha\in \pi_{p+1}(X),\,\beta\in \pi_{q+1}(X),\,\gamma\in\pi_{m}(S^{p}),\,\delta\in\pi_{n}(S^{q})$ and $ \iota_{1}=[\mathrm{id}_{S^{1}}]$,\;then,\;$[ \alpha\circ\s\gamma,\beta\circ\s\delta ]=[\alpha,\beta]\circ(\gamma\wedge\delta\wedge\iota_{1})=[\alpha,\beta]\circ\s(\gamma\wedge\delta)$.$\hfill\boxed{}$
\end{lem}

\begin{cor}\label{fkgk}
Suppose $k\in\z_{+}$,\, then,\;\\ after localization at 2, $$f_{k}=\pm[\iota_{4+k},\iota_{4+k}] \hc  \nu_{2k+7}=\pm P(\iota_{2k+9})\hc \nu_{2k+7};$$
after localization at 3,
$$g_{k}=\pm[\iota_{4+k},\iota_{4+k}] \circ \alpha_{1}(2k+7).$$
\end{cor}

\begin{proof}\;By Lemma\,\;\ref{wgs},\\\\
\indent$f_{k}=[\iota_{4+k},\nu_{4+k}]\\\indent\quad=[ \iota_{4+k}\hc\s\iota_{3+k}, \iota_{4+k}\hc\s\nu_{3+k} ]\\\indent\quad=[\iota_{4+k},\iota_{4+k}] \hc\s (\iota_{3+k}\wedge \nu_{3+k})\\\indent\quad=\pm[\iota_{4+k},\iota_{4+k}] \hc  \nu_{2k+7}$, \;\\\\on $  [\iota_{4+k},\iota_{4+k}]=\pm P(\iota_{2k+9})$ and $P(\s^{2}\alpha_{2m-1})=P(\iota_{2m+1})\hc \alpha_{2m-1}$\,where\, $\alpha_{2m-1}\in \pi_{*}(S^{2m-1})$,  see \cite[Proposition\, 2.5]{Toda};\\\\
\indent $g_{k}=[\iota_{4+k},  \alpha_{1}(4+k)    ]\\\indent\quad=[\iota_{4+k}\hc \s\iota_{3+k},  \iota_{4+k}\hc\s\alpha_{1}(3+k)    ]\\\indent\quad=[\iota_{4+k},\iota_{4+k}] \circ \s(\iota_{3+k}\wedge\alpha_{1}(3+k))\\\indent\quad=\pm[\iota_{4+k},\iota_{4+k}] \circ \alpha_{1}(2k+7).$

\end{proof}

\begin{lem}\label{gjlx}
\indent After localization at 2, up to homotopy, \\  \indent\quad$\mathrm{sk}_{20}(F_{1})=S^{5}\vee S^{13}$,\;
$\mathrm{sk}_{23}(F_{2})=S^{6}\cup _{2\bar{\nu}_{6}} e^{15},$\\
\indent\quad $\mathrm{sk}_{26}(F_{3})=S^{7}\vee S^{17},$ $\;\mathrm{sk}_{29}(F_{4})=S^{8}\cup _{f_{4}} e^{19},$ \\
\indent\quad $\mathrm{sk}_{32}(F_{5})=S^{9}\cup _{\bar{\nu}_{9}\nu_{17}} e^{21},\;\mathrm{sk}_{35}(F_{6})=S^{10}\cup _{P( \nu_{21})} e^{23}$,\\
\indent\quad $\mathrm{sk}_{38}(F_{7})=S^{11}\cup _{  \cg_{11}\nu_{18}^{2}} e^{25},\;\;\mathrm{sk}_{41}(F_{8})=S^{12}\cup_{P( \nu_{23})}  e^{27}$,\\
where $f_{4}=\nu_{8}\sigma_{11}-2t'\sigma_{8}\nu_{15}$ for some odd $t'$.\\
\indent After localization at 3, up to homotopy,\\
\indent\quad $\mathrm{sk}_{20}(F_{1})=S^{5}\vee S^{13}$,\;\;$\mathrm{sk}_{23}(F_{2})=S^{6}\cup _{g_{2}}  e^{15},$\\
\indent\quad $\mathrm{sk}_{26}(F_{3})=S^{7}\vee S^{17},$ $\;\mathrm{sk}_{29}(F_{4})=S^{8}\cup_{g_{4}} e^{19},$ \\
\indent\quad  $\mathrm{sk}_{32}(F_{5})=S^{9}\vee S^{21},\;\mathrm{sk}_{35}(F_{6})=S^{10}\cup_{g_{6}} e^{23}$\\
  here, $g_{k}=\pm[\iota_{4+k},\iota_{4+k}] \circ \alpha_{1}(2k+7)\in\pi_{11+2k}(S^{4+k})$  \;($k=2,4\;or\;6)$ are all of order 3.
\end{lem}

\begin{proof}
In this proof  in the 2-local case, we will freely  use some well-known results on  $[\iota_{m},\iota_{m}] $ and some relations of the generators  of   $\pi_{*}(S^{n})$ due to \textit{H.\,Toda} in \cite{Toda}, however,\textit{ more conveniently}, one can also find   the results we used on   $[\iota_{m},\iota_{m}] $  in \cite[(1.3), p.\;3]{Oguchi}, and find   the relations of the generators  of   $\pi_{*}(S^{n})$ we used in  \cite[table, p.\;104]{Oguchi}.\\
After localization at 2: \\
\indent $f_{1}=[\iota_{5},\iota_{5}]\hc\nu_{9}=\nu_{5}\ca_{8}\nu_{9}=0$,\\
\indent $f_{2}=[\iota_{6},\nu_{6}]=-2\ch_{6}$,\\
\indent $f_{3}=[\iota_{7},\iota_{7}]\hc\nu_{13}=0\hc\nu_{13}=0$,\\
\indent $\pm f_{4}=[\iota_{8},\iota_{8}]\hc\nu_{15}= (\s \cg'-2\cg_{8})\hc\nu_{15}=   (\s \cg')\hc\nu_{15}   -2\cg_{8}\nu_{15}=x\nu_{8}\cg_{11}-2  \cg_{8}\nu_{15}   $,\;$(x:\;\mathrm{odd})$,\\ $f_{4}$ and $tf_{4}$\, ($t:\mathrm{odd})$ have the same cofibre, so we can replace the original $f_{4}$ by $f_{4}=\nu_{8}\sigma_{11}-2t'\sigma_{8}\nu_{15}$, ($t':\mathrm{odd}).$ \\
\indent $f_{5}=[\iota_{9},\iota_{9}]\hc\nu_{17}=(\cg_{9}\ca_{16}+\lt_{9}+\ch_{9})\nu_{17}=\ch_{9}\nu_{17}$,\\
\indent $f_{6}=\pm P(\nu_{21})$.\\
\indent $f_{7}=[\iota_{11},\iota_{11}]\hc\nu_{21}=(\cg_{11}\nu_{18})\nu_{21}=\cg_{11}\nu_{18}^{2}$.\\
After localization at 3: 
$$\pm g_{k}=[\iota_{4+k},  \alpha_{1}(4+k)    ]=[\iota_{4+k},\iota_{4+k}] \circ \alpha_{1}(2k+7)\in\pi_{11+2k}(S^{4+k}),\;k\geq1.$$
If $4+k\ty1\m2$,  by\,\cite[(13.1), Serre's isomorphism]{Toda}, we know\;\\$\s :\pi_{11+2k}(S^{4+k})\rightarrow\pi_{12+2k}(S^{5+k})$ are all monomorphisms.  Since the Whitehead product is in  $\mathrm{Ker}(\s),$   thus, $g_{1}=g_{3}=g_{5}=0$. If $4+k\ty0\m2$,  by\,\cite[(13.1), Serre's isomorphism]{Toda},  $g_{k}$   are the isomorphic images of  the 3-order elements $\alpha_{1}(2k+7)$, then $g_{k}\in\pi_{11+2k}(S^{4+k})$ are of order 3.
\end{proof}

\section{Main tools to compute the cokernels and kernels }
\subsection{Diagram Long$(m,k,t)$, sequence Long$(m,k)$ and so on}
We recall  Remark\,\ref{Fk},\;\;$p_{k}$,\,$F_{k}$ are the pinches and  fibres.
 \begin{definition}\label{dt}

  We fixed the symbols $i_{k} ,\partial_{k},\varphi_{k+t},w_{k+t}$ and  $\theta_{k+t}$ to denote the maps of the following homotopy commutative diagram with rows fibre sequences,  where $\pa_{k}$ is the composition $\lo S^{8+k}\stackrel{\simeq}\longrightarrow J(S^{7+k})\hookrightarrow J(M_{S^{4+k}}, S^{7+k})=F_{k}$,
  and define  $i=\s^{\wq}i_{1} $. \\\\
\indent\quad  \xymatrix@C=0.6cm{
\Omega S^{8+k}\ar[d]\ar[rr] ^{  \partial_{k}}&& F_{k} \ar[d]^{ \varphi_{k+t}   }\ar[rr]^{ i_{k}  } && \Sigma^{k}\mathbb{H}P^{2}\ar[d]_{w_{k+t}=\Omega^{t}\Sigma^{t}  }\ar[rr]^{ p_{k}  }&& S^{8+k} \ar[d]_{ \theta_{k+t} =\Omega^{t}\Sigma^{t} }\\
\Omega^{t+1} S^{8+k+t}\;\ar[rr] ^{ \Omega^{t}\partial_{k+t}  }&&  \Omega^{t}F_{k+t}\ar[rr] ^{  \Omega^{t}i_{k+t}  } &&\Omega^{t}\Sigma^{k+t} \mathbb{H}P^{2}\ar[rr] ^{  \Omega^{t} p_{k+t}}&& \Omega^{t}S^{8+k+t}}.
 \end{definition}

 \indent We recall Proposition\,\ref{grdl}(1),  for a map or a homotopy class $K\stackrel{f}\longrightarrow L $,  the inclusion $( K,K)\hookrightarrow (M_{L},K)$ is denoted by $i_{f}$. \begin{definition}
\label{4gm}\begin{itemize}
\item[\rm(1)]
There exists the following commutative diagrams with exact rows  induced by the above.\;If $k\geq1$,   $\pa_{k*}:=(J(\mathbbm{i}))_{*}\hc \mathbbm{t}$,   $\mathbbm{t}$ is the isomorphism $\pi_{m+1}( S^{8+k})\rightarrow   \pi_{m}( J(S^{8+k}))$,\;\;$\mathbbm{i}=i_{\nu_{4+k}}$ if localized at 2, and  $\mathbbm{i}=i_{\alpha_{1}(4+k)}$ if localized at 3.\;This diagram is called  $\mathrm{ Long}(m,k,t)$.\\   \resizebox{11.5cm}{0.45cm}{\xymatrix@C=0.15cm{
\pi_{m+1}( S^{8+k})\ar[d]^{\s^{t} }\ar[rr] && \pi_{m}(F_{k}) \ar[d]^{ \varphi_{k+t_{*}}   }\ar[rr]^{ i_{k_{*}}  } &&\pi_{m}( \Sigma^{k}\mathbb{H}P^{2})\ar[d]^{w_{k+t_{*}}=\s^{t} }\ar[rr]^{ p_{k_{*}}  }&&\pi_{m}( S^{8+k}) \ar[d]^{ \s^{t}}\ar[d]\ar[rr]^{  \partial_{k_{*_{m}}}}&& \pi_{m-1}(F_{k})\ar[d]^{ \varphi_{k+t_{*}}   }\\
    \pi_{m+1+t}( S^{8+k+t})\ar[rr] &&\pi_{m+t}(F_{k+t})\ar[rr] ^{  i_{k+t_{*}} \;\;\; } &&\pi_{m+t}(\Sigma^{k+t} \mathbb{H}P^{2})\ar[rr] ^{ \;\;\;\;p_{k+t_{*}}}&&\pi_{m+t}(S^{8+k+t})\ar[rr] &&\pi_{m-1+t}(F_{k+t}).}}{\textcolor{white} {.}}

\begin{equation}\label{longmkt}  \;\mathrm{ Long}(m,k,t)
\end{equation}
\item[\rm(2)] $\mathrm{Long}(m,k,t)$ induces the following commutative diagram with  short exact rows,\,the following diagram is called\; $\mathrm{Short}(m,k,t)$.\;Here,  $i_{k_{*}} $ induced \;$\cok(\pa_{k*_{m+1}})\tg i_{k_{*}}(\pi_{m}(F_{k}))$,   by abuse of notation, $ i_{k_{*}}(\pi_{m}(F_{k}))$ is usually  denoted by $\cok(\pa_{k*_{m+1}})$.\\\\ \xymatrix@C=0.26cm{
0\ar[rr] && i_{k_{*}}(\pi_{m}(F_{k})) \ar[d]^{\mathcal{A}(m,k,t) = \s^{t}}\ar[rr]^{ \subseteq} &&\pi_{m}( \Sigma^{k}\mathbb{H}P^{2})\ar[d]^{\mathcal{B}(m,k,t)= \s^{t} }\ar[rr]^{p_{k_{*}}}&&\Ker(\partial_{k_{*_{m}}})\ar[d]^{  \mathcal{C}(m,k,t)=   \s^{t}    }\ar[d]\ar[rr]&&0\\
0\ar[rr] &&i_{k+t_{*}}(\pi_{m+t}(F_{k+t})) \ar[rr] ^{ \subseteq} &&\pi_{m+t}(\Sigma^{k+t} \mathbb{H}P^{2})\ar[rr] ^{\;\;\,p_{k+t_{*}}}&&\Ker(\partial_{k_{*_{m+t}}})\ar[rr] &&0}

 \begin{equation}\label{shortmkt} \mathrm{Short}(m,k,t)
\end{equation}
\item[\rm(3)]The following short exact sequence is called\;$\mathrm{Short}(m,k),$\\\\
\indent \qquad \xymatrix@C=0.3cm{ 0 \ar[r] &  \cok(\partial_{k_{*_{m+1}}})\ar[rr]&&\pi_{m}( \Sigma^{k}\mathbb{H}P^{2}) \ar[rr] && \Ker(\partial_{k_{*_{m}}})  \ar[r]& 0}
{\textcolor{white} {  555}}
 \begin{equation}\label{shortmk}  \;\mathrm{Short}(m,k)
\end{equation}
The following short exact sequence is still called\;$\mathrm{Short}(m,k).$\\\\\indent  \xymatrix@C=0.3cm{
 &&&& 0 \ar[r] &i_{k_{*}}(\pi_{m}(F_{k}))\ar[rr]^{\bz}&&\pi_{m}( \Sigma^{k}\mathbb{H}P^{2}) \ar[rr] && \Ker(\partial_{k_{*_{m}}}) \ar[r] & 0, }
 {\textcolor{white} {.}}\\\\ by abuse of notation,  $ i_{k_{*}}(\pi_{m}(F_{k}))$ is usually still denoted by $\cok(\pa_{k*_{m+1}})$.
  \item[\rm(4)] The following exact sequence is called $\mathrm{Long}(m,k),$\;\\\\
\resizebox{10.9cm}{0.5cm}{
\indent\xymatrix@C=0.85cm{
  \pi_{m+1}( S^{8+k})\ar[r]^{\; \quad \partial_{k_{*_{m+1}}}}& \pi_{m}(F_{k}) \ar[r]^{ i_{k_{*}}\quad  }  & \pi_{m}( \Sigma^{k}\mathbb{H}P^{2})\ar[r]^{p_{k_{*}} } & \pi_{m}( S^{8+k})  \ar[r]^{  \partial_{k_{*_{m}}}} & \pi_{m-1}(F_{k}).  }}

   \begin{equation}\label{longmk} \;\mathrm{  \;Long}(m,k)
\end{equation}
\end{itemize}

 \end{definition}

\subsection{Some diagrams on skeletons of the fibres}

 \begin{lem}\label{xwjd}Let  $m\in\z_{+}$, suppose $K= Y\cup_{f} e ^{m+1}$ is a \textit{CW} complex,\;where  $Y=\sk _{m}(K)$ is a path-connected \textit{CW} complex,
  and suppose $\pa$ is the composition of$$\lo S^{m+1}\stackrel{\simeq}\longrightarrow J(S^{m})\stackrel{\xyd}\longrightarrow J( M_{Y},S^{m}).$$
 Then, for the fibre sequence,\\
 
 \indent \qquad \qquad  \qquad\xymatrix @C=0.5cm{
  \lo S^{m+1} \ar[r]^{\pa\;\;} &   J( M_{Y},S^{m}) \ar[r] & K \ar[r] & S^{m+1},   }\\\\
 there exists a homotopy-commutative diagram,\\
  \indent \qquad \qquad  \qquad\qquad\qquad \xymatrix{
    S^{m} \ar[d]_{\xyd} \ar[r]^{f} & Y \ar[d]^{\subseteq} \\
    \lo S^{m+1}  \ar[r]^{\pa\;} &  J( M_{Y},S^{m}). }
    \end{lem}

    \begin{proof}
By Proposition\,\ref{grdl}(1), the following is homotopy-commutative,\\
\xymatrix@!C=2.5cm{   & S^{m} &   &  &  \\
& S^{m} \ar[r] ^{\subseteq}  &J(S^{m})\ar[r] _{\simeq} &\lo S^{m+1}&\\
 Y \ar[r]  ^{\subseteq}_{\simeq} &M_{Y}\ar[r] ^{\subseteq}  &J( M_{Y},S^{m})   &  &\\
      \ar"2,2";"1,2" ^{\mathrm{id}} \ar"1,2";"2,4"^{E}  \ar"2,2";"3,2" ^{\subseteq}  \ar"2,3";"3,3" _{\subseteq}   \ar"2,4";"3,3"^{\partial}    \ar"2,2";"3,1"_{  f   }
     }

hence the result holds.

\end{proof}

\begin{lem}\label{com} Suppose $k\geq1$, then  there exists the homotopy-commutative diagram, \\\indent \qquad\qquad\qquad\qquad\xymatrix@C=1.6cm{
  S^{4+k} \ar[d]_{\lo^{t}\s^{t}} \ar[r]^{j_{k}} & F_{k} \ar[d]^{\phi_{k+t}} \\
  \lo^{t}S^{4+k+t} \ar[r]^{\lo^{t}j_{k+t}} & \lo^{t}F_{k+t}   }

 \noindent  it induces the commutative diagram,\\
\indent \qquad\qquad\qquad\xymatrix@C=1.4cm{
 \pi_{m}( S^{4+k}) \ar[d]_{\s^{t}} \ar[r]^{j_{k_{*}}} &  \pi_{m}(F_{k}) \ar[d]^{\phi_{k+t_{*}}} \\
 \pi_{m+t}(S^{4+k+t}) \ar[r]^{j_{k+t_{*}}} & \pi_{m+t}(F_{k+t})   }
\begin{equation}\label{com1}
    \mathrm{COM}(m,k,t)
\end{equation}
 \noindent We call the above commutative diagram of homotopy groups  $\mathrm{COM}(m,k,t)$.

\end{lem}

\begin{proof}
For any triad $(m,k,t)$, firstly, we consider$\lon(4+k,k,t)$ given by diagram (\ref{longmkt}).\;In this case, $  \phi_{k+t} $ induces $ \pi_{4+k}(    F_{k} ) \tg \pi_{4+t}(\lo^{t}F_{k+t})$.\;

By the naturality of the   Hurewicz isomorphisms and  the adjoin, we have the commutative diagram,\\
\indent\qquad\qquad\qquad\xymatrix@C=2cm{
  H_{4+k+t}(S^{4+k+t}) \ar[d]_{\tg} \ar[r]^{j_{k+t_{*}}} & H_{4+k+t}(F_{4+k+t}) \ar[d]^{\tg} \\
   \pi_{4+k+t}(S^{4+k+t}) \ar[d]_{\tg} \ar[r]^{j_{k+t_{*}}} & \pi_{4+k+t}(F_{4+k+t}) \ar[d]^{\tg} \\
 \pi_{4+k}(\lo^{t}S^{4+k+t})  \ar[r]^{(\lo^{t}j_{k+t})_{*}} &  \pi_{4+k+t}(\lo^{t}F_{4+k+t})   } {\textcolor{white} {.}}\\\\
\noindent By Corollary\,\ref{bqjg}, $F_{4+k+t}=S^{4+k+t}\cup e^{11+2k+2t}\cup\cdots$, \\
  thus $$H_{4+k+t}(S^{4+k+t})\stackrel{ j_{k+t_{*}}    }\longrightarrow  H_{4+k+t}(F_{4+k+t})$$ is an isomorphism, 
 successively, $$\pi_{4+k}(\lo^{t}S^{4+k+t})\stackrel{ j_{k+t_{*}}    }\longrightarrow   \pi_{4+k+t}(\lo^{t}F_{4+k+t}) $$ is an isomorphism.\\
Therefore, all of  the 4 maps of the following diagram induce isomorphisms between $\pi_{4+k}(-)$. \\
\indent \qquad\qquad\qquad\qquad\xymatrix@C=1.6cm{
  S^{4+k} \ar[d]_{\lo^{t}\s^{t}} \ar[r]^{j_{k}} & F_{k} \ar[d]^{\phi_{k+t}} \\
  \lo^{t}S^{4+k+t} \ar[r]^{\lo^{t}j_{k+t}} & \lo^{t}F_{k+t}   }

\noindent By the cellular approximation theorem,  any map  $f:S^{4+k}\rightarrow\lo^{t}F_{k+t}  $  can decompose as the following up to homotopy, \\

\indent\qquad\qquad\xymatrix@C=1.9cm{ S^{4+k} \ar[r]^{( x(f))\mathbf{1}_{ S^{4+k}}} & S^{4+k} \ar[r]^{  \mathrm{incl}\qquad\qquad}_{\subseteq\qquad\qquad} &\lo^{t}F_{k+t},\qquad x(f)\in\z,  }{\textcolor{white} {.}}\\\\
\noindent  (regarding $ x(-)$ as a $\z$-value function),\\
 thus\;$$\pi_{4+k}(\lo^{t}F_{k+t})=\z\n [\mathrm{incl}] \nn,$$
$f$ induces $f_{*}:\pi_{4+k}(S^{4+k})\rightarrow \pi_{4+k}(\lo^{t}F_{k+t}),\;\;\iota_{4+k}\mapsto x(f) \cdot [\mathrm{incl}]$.\\
Thus, $$x(  \phi_{k+t}   \hc j_{k} )=\pm1,  \;\;\;x(       \lo^{t}j_{k+t}\hc(\lo^{t}\s^{t}) )=\pm1.$$
 If $x(  \phi_{k+t}   \hc j_{k} )=x(       \lo^{t}j_{k+t}\hc(\lo^{t}\s^{t}) )$, our lemma is established,  if  $x(  \phi_{k+t}   \hc j_{k} )=-x(       \lo^{t}j_{k+t}\hc(\lo^{t}\s^{t}) )$,  replace  $j_{k}$ by  $-j_{k}$, that is,  denote the original  $-j_{k}$ by $j_{k}$.

\end{proof}

 \begin{rem} \label{xt}

We fix the notations $ j_{k}, j'_{k}$  to denote the inclusions of the following homotopy-commutative diagram,\\
\indent \xymatrix{
  &&&S^{4+k}\ar[d]_{j'_{k}} \ar[dr]^{j_{k}}        \\
 &&& F_{k}\ar[r]_{ i_{k}\quad }  & \s ^{k}\h P^{2}.\;\;\;           }\\\\

We fix the notations $ \widehat{ j}_{k},\;\widetilde{ j}_{k}$       to denote the inclusions of the following homotopy-commutative diagram,\\
\indent\qquad \qquad \qquad \qquad \xymatrix@C=1.7cm{
  S^{4+k}\ar[dr]_{j'_{k}} \ar[r]^{\widetilde{j}_{k}}
                & sk_{3k+17}(F_{k}) \ar[d]^{\widehat{j}_{k}}  \\
                & F_{k}\;.         }

 \noindent The  homotopy-commutative properties of them are easy to get, similarly to the proof of  Lemma\,\ref{com},  we  need only to notice that
  $\pi_{4+k}(\s^{k}\y) =\z\n [j_{k}] \nn$,   $\pi_{4+k}(F_{k}) =\z\n [j'_{k}] \nn$, $\pi_{4+k}(sk_{3k+17}(F_{k})) =\z\n [\;\widetilde{j}_{k}\;] \nn$  and  $\widehat{ j}_{k}$ induces an isomorphism of $H_{4+k}(-)$, (of course, if necessary, replace $i_{k}$ by $-i_{k}$, still denoted by $i_{k}$, or replace  $j'_{k}$ by $-j'_{k}$,\,still denoted by\;$j_{k}$).\\
 \indent However, by abuse of notation, without giving rise to ambiguity,  $$ j_{k}, j'_{k}, \widehat{ j}_{k},\;\widetilde{ j}_{k} \;\text{and}\; \s^{\wq} j_{1},  \s^{\wq}j'_{1}, \s^{\wq}\widehat{ j}_{1},\;\s^{\wq}\widetilde{ j}_{1}$$ are all denoted by $j$ for simplicity.

 \end{rem}

\subsection{Diagrams BUND$(m+1,k)$  and D$(m+1,k)$}
\indent  By Lemma\,\;\ref{xwjd} and  Lemma\,\;\ref{wn}, we have the following.
\begin{cor}\label{ab}
Suppose $\s\mathbbm{y} \in \pi_{8}(S^{5})$ and $k\in\z_{+}$,
   after localization at  $p=2\;or\;3$,  there exists homotopy-commutative diagram,\\
\indent\qquad\qquad\qquad\qquad\qquad\quad \xymatrix@C=1.5cm{
  S^{7+k} \ar[d]_{\lo\s} \ar[r]^{\s^{k}\mathbbm{y} } & S^{4+k} \ar[d]^{j_{k}} \\
  \lo S^{8+k} \ar[r]^{\pa_{k}} & F_{k}   }
\indent\\\\ \noindent take  $\s^{k}\mathbbm{y}=\nu_{4+k}$ if they are localized at 2,  take  $\s^{k}\mathbbm{y}=\alpha_{1}(4+k)$ if they are  localized at 3; it induces the commutative diagram   after localization at  $p=2\;or\;3$, \\
 \indent\qquad\qquad\qquad\qquad\xymatrix@C=1.5cm{
  \pi_{m}(S^{7+k}) \ar[d]_{\s} \ar[r]^{(\s^{k}\mathbbm{y})_{*}  } &  \pi_{m}( S^{4+k}) \ar[d]^{j_{k*}} \\
  \pi_{m+1}(S^{8+k}) \ar[r]^{\pa_{k_{*m+1}}} &   \pi_{m}(F_{k})   }
  \begin{equation}\label{bund1}\mathrm{BUND}( m+1,k)\end{equation}
 \noindent We call this diagram  \;$\mathrm{BUND}(m+1,k)$. $\bzd$
\end{cor}

\subsection{Sequcences PISK$(m,k)$ and diagrams D$(m+1,k)$}
Similarly  to Corollary\,\;\ref{ab},\;we have

\begin{cor}\label{sk111} Let $k\in\z_{+},$ then there exists the  fibre sequence,\\

\indent\qquad\xymatrix@C=0.3cm{
  \lo S^{11+2k} \ar[r]^{d\qquad} &J(M_{S^{4+k}}, S^{10+2k})  \ar[r]^{} & \sk_{3k+17}(F_{k}) \ar[r]^{\;\;\,q} & S^{11+2k},  }

\noindent they induce the following exact sequence,\;($d_{*}=d_{k*m+1},d_{k*m}$ respectively), and we call the sequence  $\mathrm{PISK}(m,k),$ \,\\\\ \resizebox{11.9cm}{0.36cm}{
\xymatrix@C=0.3cm{
                  \pi_{m+1} (S^{11+2k}) \ar[r]^{d_{*}\qquad} &\pi_{m}(J(M_{S^{4+k}}, S^{10+2k}))  \ar[r]^{} &\pi_{m}( \sk_{3k+17}(F_{k})) \ar[r]^{\;\;\;\;\,q_{*}} & \pi_{m}(S^{11+2k})  \ar[r]^{d_{*}\qquad\quad}& \pi_{m-1}(J(M_{S^{4+k}}, S^{10+2k})).  }}

\begin{equation}\label{pisk1}
    \mathrm{PISK}(m,k)
\end{equation}

There exists the commutative diagram and we call it  $\mathrm{D}(m+1,k).$\\
\indent\qquad\qquad\qquad\xymatrix@C=1.4cm{
\pi_{m}(S^{10+2k}) \ar[d]_{\s} \ar[r]^{\mathbbm{b}_{k}} &\pi_{m}(S^{4+k}) \ar[d]^{j_{*}} \\
  \pi_{m+1}(S^{11+2k}) \ar[r]^{d_{*}=d_{k*_{m+1}}\qquad} & \pi_{m}(        J(M_{S^{4+k}}, S^{10+2k})       )  }
\begin{equation}\label{dkm}
    \mathrm{D}(m+1,k)
\end{equation}
 \noindent we take $\mathbbm{b}_{k}=f_{k}$ if localized at 2,\\
 \noindent we take $\mathbbm{b}_{k}=g_{k}$ if localized at 3. \\
 (The results on $f_{k}, g_{k}$ we need are shown in the proof of Lemma\,\ref{gjlx}.)\qquad$\boxed{}$
\end{cor}

\section{$\pi_{r+k}(\s^{k}\h P^{2})$ localized at 3}

In this section,  \textit{all have been localized at 3}, whether we say after localization at 3 or not; and $\widetilde{H}_{*}(-)$  denotes  reduced homology with $\z/3$ coefficients. For an element $\mathbbm{a}\in\pi_{*}(X)$, we usually use $\overline{\,\mathbbm{a}\,}$ to denote its lifting (up to sign).\\

\subsection{  Determinations of $\pi_{r+k}(\s^{k}\h P^{2})\;\;(7\leq r\leq13$\;but\;$r\neq 11)$ }

\begin{thm} After localization at 3,\\
\[\pi _{7+k}(\Sigma ^{k}\mathbb{H}P^{2})\tg\begin{cases}
0,\;&\text{$k\geq1$}\\
\mathbb{Z}/3,\;&\text{$k=0$}
\end{cases}\
, \quad\pi _{8+k}(\Sigma ^{k}\mathbb{H}P^{2})\tg\begin{cases}
\mathbb{Z}_{(3)},\;&\text{$k\geq1$}\\
0,\;&\text{$k=0$}
\end{cases}\]\\
\[\pi _{9+k}(\Sigma ^{k}\mathbb{H}P^{2})\tg\begin{cases}
\mathbb{Z}_{(3)},\;&\text{$k=2$}\\
0,\;&\text{$\mbox{else}$}
\end{cases}, \quad
\pi _{10+k}(\Sigma ^{k}\mathbb{H}P^{2})\tg\begin{cases}
0,\;&\text{$k\geq1$}\\
\mathbb{Z}/3,\;&\text{$k=0$}
\end{cases}\]
\[\pi _{12+k}(\Sigma ^{k}\mathbb{H}P^{2})\tg\begin{cases}
0,\;&\text{$k\geq0\;\mbox{and}\;k\neq1$}\\
\mathbb{Z}_{(3)},\;&\text{$k=1$}
\end{cases}\],
\[\pi _{13+k}(\Sigma ^{k}\mathbb{H}P^{2})\tg\begin{cases}
0,\;&\text{$k\geq0\;\mbox{and} \;k\notin\n2,\,6\nn$}\\
    \mathbb{Z}_{(3)},\;&\text{$k=2\;\mbox{or}\;6$}\\
\end{cases}\]

\end{thm}

\begin{proof}
By  the corresponding exact sequences$\lon (m,k)$ given by sequence (\ref{longmk}),\,\;
 and by Lemma\,\;\ref{gjlx},\;Lemma\,\;\ref{y9jy} and the group structure of  $\pi_{*}(S^{n}) $ in \cite[Toda]{Toda}, we get the result. In these cases, the  $\z_{(3)}$-module extension problems for these  exact sequences are all trivial.
\end{proof}

\begin{lem}\label{phcp8}  $p\hc\pi_{8}^{S}(\y)=\z_{(3)}\n  3\iota\nn.$

\end{lem}

\begin{proof}  For the cofibre sequence \;$S^{8}\stackrel{\af_{1}(5)}\longrightarrow S^{5}\stackrel{i_{1}}\longrightarrow\s\h P^{2}\stackrel{p_{1}}\longrightarrow S^{9},$ \;applying $ \pi_{9}^{S}(-)$, we have the exact sequence,\\
\xymatrix@C=0.8cm{
 &&& 0 \ar[r]^{} & \pi_{8}^{S}(\y)\ar[r]^{p_{1*}} & \pi_{0}^{S}(S^{0}) \ar[r]^{\af_{1}(5)_{*}} & \pi_{3}^{S}(S^{0})   }

 \noindent  $\pi_{0}^{S}(S^{0})=\z_{(3)}\n \iota \nn$, $\af_{1}(5)_{*}(\iota)=\af_{1}$ is of order 3, then by exactness, the lemma is established.

\end{proof}

\subsection{ Determinations of $\pi_{11+k}(\s^{k}\h P^{2})$}

We recall  \cite[(13.4), p.\,176 ]{Toda}, \cite[Theorem\;(13.4), p.\,176 ]{Toda}, after localization at 3,
$\pi_{3}^{S}(S^{0})=\z/3\n \alpha_{1}  \nn$,\;$\pi_{7}^{S}(S^{0})=\z/3\n \alpha_{2}  \nn$, and\;  $ \alpha_{2}=\langle   \alpha_{1}, 3\iota, \alpha_{1}                \rangle  $.

\begin{thm}\label{3c11w} After localization at 3, $$\pi_{11}^{S}(\h P^{2})=\z/9\n\overline{\alpha_{1} }\nn,$$
where $\overline{\alpha_{1} }\in \langle j, \alpha_{1},\alpha_{1}        \rangle$,\;$j\hc \af_{2}=\pm3\overline{ \af_{1}}$.

\end{thm}

\begin{proof}
For the cofibre sequence \;$S^{8}\stackrel{\af_{1}(5)}\longrightarrow S^{5}\stackrel{j_{1}}\longrightarrow\s\h P^{2}\stackrel{p_{1}}\longrightarrow S^{9}\rightarrow\cdots$, applying  $ \pi_{12}^{S}(-)$, we have the exact sequence,\\\\\xymatrix@C=0.25cm{
 &&&& 0 \ar[r] & \pi_{7}^{S}(S^{0}) \ar[rr]^{  j_{1*} \;\;\,  } &&\pi_{11}^{S}(\h P^{2})\ar[rr]^{p_{1*}} && \pi_{3}^{S}(S^{0})  \ar[r] & 0 },\\
$$\pi_{7}^{S}(S^{0})=\z/3\n \alpha_{2}  \nn, \;\pi_{3}^{S}(S^{0})=\z/3\n \alpha_{1}  \nn,\;\text{and}\;  \alpha_{2}=\langle    \alpha_{1}, 3\iota, \alpha_{1}                \rangle.  $$ Since $\af_{1}\hc\af_{1}\in\pi_{6}^{S}(S^{0})=0$, then \;$ \langle j, \alpha_{1},\alpha_{1}        \rangle$ is well-defined. By \cite[3.9), i) and ii), p.\,33]{Toda}, we have $\x\af_{1}, \af_{1},3\iota \xx=-\af_{2}=-\langle    \alpha_{1}, 3\iota, \alpha_{1}                 \rangle  $.

\noindent Thus,\\
\indent\qquad\qquad\qquad\quad$ j_{1*}( \alpha_{2})\in  j_{1*}\langle    \alpha_{1}, 3\iota, \alpha_{1}                 \rangle
 \\ \indent\qquad\qquad\qquad \qquad\qquad =j\hc\langle    \alpha_{1}, 3\iota, \alpha_{1}                 \rangle\\\indent\qquad\qquad\qquad \qquad\qquad =-j\hc\x\af_{1}, \af_{1},3\iota \xx \\\indent\qquad\qquad\qquad \qquad\qquad =\pm\x  j,\af_{1}, \af_{1}\xx\hc 3\iota\\\indent\qquad\qquad\qquad \qquad\qquad =\pm3\x  j,\af_{1}, \af_{1}\xx.$\\\\\noindent Thus $ j_{1*} (\alpha_{2})$ can be divisible  by 3.\;we take $\overline{\alpha_{1} }\in \langle j, \alpha_{1},\alpha_{1}       \rangle$
such that  $j_{1*}( \alpha_{2})=j\hc \af_{2}=\pm3\overline{ \af_{1}}$, then we get the result.

\end{proof}

\begin{lem}\label{y9jy}
 Suppose  $\mathbb{H}P^{3}=\mathbb{H}P^{2}\cup_{h}e^{12}$  where $h:S^{11}\rightarrow\mathbb{H}P^{2}$ is the homotopy class of the attaching map,  and suppose $k\geq1$, then$$
\forall\;\mathbbm{x} \in\pi_{11+k}(\Sigma^{k}\mathbb{H}P^{2}),\;\;\Sigma^{k}h\neq3\mathbbm{x}.$$
\end{lem}
\begin{proof}
 To obtain a contradiction, assume   $\cz \mathbbm{x}\in\pi_{11+k}(\Sigma^{k}\mathbb{H}P^{2})$,\; such that $\Sigma^{k}h=3\mathbbm{x}$.
 Let $ q: \Sigma^{k}\mathbb{H}P^{2} \rightarrow\Sigma^{k}\mathbb{H}P^{2}/ \Sigma^{k}\mathbb{H}P^{1}$  be the pinch map,  we consider
 the space $$\Sigma^{k}\mathbb{H}P^{3}/ \Sigma^{k}\mathbb{H}P^{1}=S^{8+k}\cup_{q\circ\Sigma^{k }h}e^{12+k}.$$
 By the assumption  and  $\pi_{11+k}(S^{8+k})\tg\mathbb{Z}/3,\;(k\geq1)$\;\cite[Toda]{Toda},\; we have
   \begin{center}\par $q\circ\Sigma^{k }h=q\circ3\mathbbm{x}=3 (q\circ \mathbbm{x})$=0,\end{center}
then $\Sigma^{k}\mathbb{H}P^{3}/ \Sigma^{k}\mathbb{H}P^{1}$$\simeq S^{8+k}\vee S^{12+k}$.\\
  \indent On the one hand, $\widetilde{H}_{\ast}$($ \Sigma^{k}\mathbb{H}P^{3}/\Sigma^{k}\mathbb{H}P^{1}$)
   =$\mathbb{Z}/3\{a,\;b \} $, $(|a|=12+k, |b|=8+k)$, with   Steenrod operation  $P^{_1}_{\ast}(a)=b$,
   (for this  $P^{_1}_{\ast}$, see\;\cite[Example\;4.L.4, p.\,492]{Hatcher}).\\
     \indent  On the other hand, $\widetilde{H}_{\ast}(  S^{8+k}\vee S^{12+k}$)  is splitting as a Steenrod module.\\ Then, $\Sigma^{k}\mathbb{H}P^{3}/ \Sigma^{k}\mathbb{H}P^{1}$  cannot be homotopy equivalent to $ S^{8+k}\vee S^{12+k}$, a contradiction. This forces   $ \Sigma^{k}h\neq3\mathbbm{x}$,
$\forall\;\mathbbm{x} \in\pi_{11+k}(\Sigma^{k}\mathbb{H}P^{2})$.

\end{proof}

\begin{thm}\label{3c11th}
After localization at 3,\\
\[\pi _{11+k}(\Sigma ^{k}\mathbb{H}P^{2})\tg\begin{cases}
\mathbb{Z}/9,\;&\text{$k\geq1\; \mbox{but}\;k\neq4$}\\
\mathbb{Z}/9\oplus\mathbb{Z}_{(3)},\;&\text{$k=4$}\\
\mathbb{Z}/3\oplus\mathbb{Z}_{(3)},\;&\text{$k=0$}
\end{cases}\]

\end{thm}

\begin{proof}
 For $k\geq 5$,\;$\pi _{11+k}(\Sigma ^{k}\mathbb{H}P^{2})$ are in the stable range,  by Theorem\,\ref{3c11w}\;we have\;
$\pi _{11+k}(\Sigma ^{k}\mathbb{H}P^{2})\tg\z/9,\;(k\geq 5).$ \\\indent For $k=3,2$ or 1, we consider$\lon (11+k,k)$ given by sequence (\ref{longmk}),
  by Lemma\,\;\ref{gjlx},   Lemma\,\;\ref{y9jy} and the groups  $\pi_{*}(S^{n}) $  in \cite[Toda]{Toda}, we get    \,$\pi _{14}(\Sigma ^{3}\mathbb{H}P^{2})$, $\pi _{13}(\Sigma ^{2}\mathbb{H}P^{2})$ and $\pi _{12}(\Sigma \mathbb{H}P^{2})$  are all groups of order 9, the  elements $\s^{k}h$ cannot be divisible by 3, thus $\s^{k}h$   are of order 9, ($k=3,\;2$\;or\;1),  therefore,  these three groups are $\z/9\n  \s^{k}h \nn$ \;($k=3,\;2\;$or\;1).\\
\indent For $k=4$, we consider$\lon (15,4)$ given by sequence (\ref{longmk}),
by the above, $\ord(\s^{3}h)=\ord(\s^{5}h)=9$, so $\pi _{15}(\Sigma ^{4}\mathbb{H}P^{2}) $ contains the element $\s^{4}h$ of order 9. By Lemma\,\,\ref{gjlx} and the groups  $\pi_{*}(S^{n}) $  in \cite[Toda]{Toda}, we have\;$$\pi_{14}(F_{4})=\pi_{16}(S^{12})=0.$$
Since\\
 $\textcolor{white}{5555555999955555577\,7}\pi_{15}(F_{4})\\\textcolor{white}{555555555599995555}=j\hc\pi_{15}(S^{8})\\\textcolor{white}{555555555222255555}=j\hc (\z/3\n \alpha_{2}(8)\nn\jia\z_{(3)}\n [\iota_{8}, \iota_{8}]\nn)
 \\\textcolor{white}{555555514145555555}=\z/3\n j\hc\alpha_{2}(8)\nn\jia\z_{(3)}\n j\hc[\iota_{8}, \iota_{8}]\nn$,\;\\\\and\;\\
 \centerline{
 $\pi_{15}(S^{12})=\z/3\n \af_{1}(12)\nn$,}
 then,\;  $$\pi _{15}(\Sigma ^{4}\mathbb{H}P^{2}) =\z_{(3)}\n \mathbbm{x}\nn\jia\z/9\n \s^{4}h\nn\; \text{for some\;} \mathbbm{x}.$$ We know $$p_{4*}(\s^{4}h )=\pm\s (p_{3*}(\s^{3}h ))= \pm\s\af_{1}(11)=\pm\af_{1}(12),$$
 and $$i_{4*}(j\hc\alpha_{2}(8) )=\pm \s( i_{3*}(j\hc\alpha_{2}(7) )  )=\pm\s (3\s^{3}h)=\pm 3\s^{4}h.$$Then, for the relation\,\\
$$i_{4*}(  j\hc[\iota_{8}, \iota_{8}])\ty 3^{r}\ell\mathbbm{x}\m \s^{4}h,\;(\text{integer}\;r\geq0,\;\ell\;\text{is invertible in}\;\z_{(3)}),$$ 
we consider the quotient, $$(\z_{(3)
}\jia\z/9)/ \x (3^{r}, b),\,(0,3)           \xx\tg \z/3^{r}\jia\z/3$$ for any $b\in\z/9$,
by exactness, $r\geq1$ is impossible.
 Hence,\;$$\pi _{15}(\Sigma ^{4}\mathbb{H}P^{2})  = \mathbb{Z}_{(3)}\n j_{4}[\iota_{8}, \iota_{8}] \nn\oplus\mathbb{Z}/9 \n \s^{4}h  \nn.$$

\end{proof}

\subsection{ Determinations of  $\pi_{15+k}(\s^{k}\h P^{2})$ }

\begin{thm}\label{3c15w}After localization at 3, \\
$$\pi_{15}^{S}(\h P^{2})=\z/27\n \overline{\af_{2}}\nn,\;\; \overline{\af_{2}}\in\x j,\af_{1},\af_{2}                          \xx.$$
\end{thm}

\begin{proof}

For the cofibre sequence \;$S^{8}\stackrel{\af_{1}(5)}\longrightarrow S^{5}\stackrel{j_{1}}\longrightarrow\s\h P^{2}\stackrel{p_{1}}\longrightarrow S^{9}\rightarrow\cdots$, applying  $ \pi_{16}^{S}(-)$, we have the exact sequence\\\xymatrix@C=0.55cm{
 && 0 \ar[r] & \pi_{11}^{S}(S^{0}) \ar[rr]^{  j_{1*} \;\;\,  } &&\pi_{15}^{S}(\h P^{2})\ar[rr]^{p_{1*}} && \pi_{7}^{S}(S^{0})  \ar[r]^{\af_{1}(5)_{*}} &  \pi_{10}^{S}(S^{0}),}

\noindent\\ where $$\pi_{7}^{S}(S^{0})=\z/3\n \alpha_{2}  \nn,$$ $$\mathrm{Im}(\af_{1}(5)_{*})=\af_{1}\hc\pi_{7}^{S}(S^{0})=\x \af_{1}\hc\af_{2}      \xx=0,\;\text{(\cite[Theorem\,13.9, p.\,180]{Toda})}.$$
Then we have  the exact sequence
\\\\\xymatrix@C=0.25cm{
 &&&& 0 \ar[r] & \pi_{11}^{S}(S^{0}) \ar[rr]^{  j_{1*} \;\;\,  } &&\pi_{15}^{S}(\h P^{2})\ar[rr]^{p_{1*}} && \pi_{7}^{S}(S^{0})  \ar[r]&0.}\\\\ \noindent  By \cite[Theorem 14.14\;ii), Toda]{Toda1} we know  $$\pi_{11}^{S}(S^{0})=\z/9\n  \af_{3}'\nn,\; \af_{3}'\in\;-\x\af_{2},\af_{1},3\iota\xx.$$ Notice that $$\af_{1}\hc\x\af_{1},3\iota,\af_{1}\xx=\af_{1}\hc\af_{2}=0, \quad3\x\af_{1},3\iota,\af_{1}\xx=3\af_{2}=0,$$
then by the Jacobi identity (\cite[(3.7),\,p.\,33]{Toda}) for $(\af_{1},\af_{1},3\iota,\af_{1},3\iota)$, we have$$0\ty \x\x\af_{1},\af_{1},3\iota\xx,\af_{1},3\iota \xx
-\x\af_{1},\x\af_{1},3\iota,\af_{1}\xx,3\iota\xx
+\x\af_{1},\af_{1},\x3\iota,\af_{1},3\iota\xx\xx,$$
in the proof of Theorem\,\ref{3c11w}, we have already shown
$\x\af_{1}, \af_{1},3\iota \xx=-\af_{2}$,\\
hence, we have $$0\ty \af_{3}'-\x \af_{1},\af_{2},3\iota\xx \m 3\af_{3}',$$ this tells us all elements in $\x \af_{1},\af_{2},3\iota\xx$ are of order 9, equivalently speaking, any element in $\x \af_{1},\af_{2},3\iota\xx$  generates $\pi_{11}(S^{0})\tg\z/9$. we take $\overline{\af_{2}}\in\x j,\af_{1},\af_{2}\xx$, then,$$3\overline{\af_{2}}\in\x j,\af_{1},\af_{2}                          \xx\hc3\iota=\pm j\hc\x\af_{1},\af_{2}\;                          3\iota\xx.$$ Thus,
$\ord(3\overline{\af_{2}})=9$, successively \,$\ord(\overline{\af_{2}})=27.$
So\,  $\pi_{15}^{S}(\h P^{2})=\z/27\n \overline{\af_{2}}\nn.$

\end{proof}

\begin{pro}\label{a2g7} $\cz \overline{\af_{2}}(7)\in\n        j_{3}, \af_{1}(7),\af_{2}(10)\nn\bz\pi_{18}(\s^{3}\y)$, \\
$\s^{\wq}\overline{\af_{2}}(7)\ty \pm \overline{\af_{2}}\m 3\overline{\af_{2}}$,\,\;$\ord(\s^{n}\overline{\af_{2}}(7))\geq27,\;\forall n\geq0.$

\end{pro}

\begin{proof}

By \cite[Theorem\,13.8, p.\,180]{Toda}, $\pi_{17}(S^{7})\tg\z/3$ and \cite[Theorem\,13.9, p.\,180]{Toda}
we have $$ \af_{1}(7)\hc \af_{2}(10)=0.$$  Thus  the Toda bracket $\n j_{3}, \af_{1}(7),\af_{2}(10)   \nn$ is well-defined, we take $\overline{\af_{2}}(7)$ from this Toda bracket.\;By Theorem\,\ref{3c15w}, and by $\pi_{11}^{S}(S^{0})\tg\z/9=\z/3^{2}, \ord(\af_{2})=3$, we have$$\s^{\wq}\overline{\af_{2}}(7)\in\pm\x        j, \af_{1},\af_{2}           \xx\ni \pm \overline{\af_{2}} $$ $\text{mod}\; j\hc\pi_{11}^{S}(S^{0})+\pi_{8}^{S}(\y)\hc\af_{2}\bz\x3^{3-2}\;\overline{\af_{2}}\; \xx=\x3\overline{\af_{2}}\; \xx.$\\
 Thus,\;$$\s^{\wq}\overline{\af_{2}}(7)\ty \pm \overline{\af_{2}}\m 3\overline{\af_{2}},$$ so $\s^{\wq}\overline{\af_{2}}(7)$ is of order 27.\;Therefore,\;$\ord(\s^{n}\overline{\af_{2}}(7))\geq27,\;\forall n\geq0.$

\end{proof}

\begin{pro}\label{p19f4}After localization at 3,\\$\pi_{18}(F_{4})=\z/3\n j\beta_{1}(8) \nn$,\\
 $\pi_{19}(F_{4})=\z/9\n   j\af'_{1}(8)  \nn\jia\z_{(3)}\n          j\hc \overline{  3\iota_{19} }  \nn,\;\;q_{*}(\overline{  3\iota_{19} })= 3\iota_{19}$.

\end{pro}

\begin{proof}
Recall $\mathrm{PISK}(m,k)$ given by sequence (\ref{pisk1}), and $\mathrm{D}(m,k)$ given by diagram (\ref{dkm}).\\
\indent By $\mathrm{PISK}(18,4)$, we have
$\pi_{18}(F_{4})=\z/3\n j\beta_{1}(8) \nn$.\\
 \indent For $\pi_{19}(F_{4})$, we consider  $\mathrm{PISK}(19,4)$.\\
 By Corollary\,\;\ref{bqjg} we have $ \sk_{15}(J(M_{S^{8}},S^{18}))=S^{8}$.
 By  $\mathrm{D}(19,4)$,
 and $g_{4}$ is of order 3, we have  $\Ker(d_{*}:\pi_{19}(S^{19})\rightarrow \pi_{18}(F_{4}))
 =\z_{(3)}\n 3\iota_{19} \nn$. Then we have the exact sequence,\\
 \xymatrix@C=0.3cm{
   &&&0 \ar[r] & \pi_{19}(S^{8}) \ar[rr]^{j_{*}\quad} &&\pi_{19}(\sk_{29} (F_{4}))  \ar[rr]^{q_{*}} && \z_{(3)}\n 3\iota_{19} \nn \ar[r] & 0. }

\noindent It is obviously splitting, so $\pi_{19}(F_{4})=\z/9\n   j\af'_{1}(8)  \nn\jia\z_{(3)}\n      j\hc   \overline{  3\iota_{19} }  \nn$, where  $q_{*}(\overline{  3\iota_{19} })= 3\iota_{19}$.

\end{proof}

\begin{pro}\label{p1617f2}After localization at 3,
\begin{itemize}
\item[\rm(1)]  $\pi_{17}(F_{2})=\z/9\n   j\af'_{1}(6)  \nn$,\;\;$\pi_{16}(F_{2})=\z/9\n   j\beta_{1}(6)  \nn$.
\item[\rm(2)]  $\Ker(\pa_{1*_{16}}:\pi_{16}(S^{9})\rightarrow \pi_{15}(F_{1}))=0$, $\Ker(\pa_{2*_{17}}:\pi_{17}(S^{10})\rightarrow \pi_{16}(F_{2}))=0.$
\end{itemize}
\end{pro}

\begin{proof}
\begin{itemize}

\item[\rm(1)]
We consider  $\mathrm{PISK}(18,2)$ given by sequence (\ref{pisk1}). By Corollary\,\ref{bqjg},\; $$ \sk_{19}(J(M_{S^{6}},S^{14}))=S^{6}.$$By  $\mathrm{D}(18,2)$ given by diagram (\ref{dkm}), the following is commutative.\\
\indent\qquad\qquad\qquad\qquad\qquad\xymatrix@C=0.75cm{
   \pi_{17}(S^{14})\ar[d]_{\s}^{\tg} \ar[r]^{g_{2*}} & \pi_{17}(S^{6})  \ar[d]^{\tg}_{j_{*}} \\
     \pi_{18}(S^{15}) \ar[r]^{d_{*}\qquad} & \pi_{17}(J(M_{S^{6}},S^{14})).  }

We know
\begin{eqnarray}
& 
&\notag g_{2*}(   \pi_{17}(S^{14})) \\\notag& 
=& g_{2}\hc \pi_{17}(S^{14})   \\\notag&
=&  [\iota_{6},\iota_{6}]\hc\af_{1}(11)\hc \x   \af_{1}(14)   \xx  \\\notag&
=& [ \iota_{6},\iota_{6}]\hc\x\;\af_{1}(11)\hc \af_{1}(14)  \;\; \xx  \\\notag&
=& 0, 
\end{eqnarray}
then  $d_{*}(\pi_{18}(S^{15}))=0$.\,Then we have the short exact sequence induced by $\mathrm{PISK}(18,2)$\\
\xymatrix@C=0.35cm{
 &&&&&&0 \ar[r]^{} & \pi_{17}(S^{6}) \ar[r]^{j_{*}\quad\;} & \pi_{17}(\sk_{23}(F_{2})) \ar[r]^{} &0   }

\noindent Hence,\;$ \pi_{17}(\sk_{23}(F_{2}))=\z/9\n   j\af'_{1}(6) \nn$.\\
\indent For $\pi_{16}(F_{2})=\z/9\n   j\beta_{1}(6)  \nn$,  it is immediately got by $\mathrm{PISK}(16,2)$ given by sequence (\ref{pisk1}).
\item[\rm(2)] 
 We recall that  $\pa_{k*_{r}}$ is the homomorphism $\pa_{k*_{r}}: \pi_{r}( S^{8+k})\rightarrow \pi_{r-1}(F_{k})$.\\
 \indent \quad For $\Ker(\pa_{1*_{16}})$, we observe $\mathrm{BUND}(16,1)$ given by diagram (\ref{bund1}). \\ \indent \quad In $\mathrm{BUND}(16,1)$:\\ by \cite[(13.1), Serre's isomorphism, p.\,172] {Toda} and \cite[Theorem 13.9, p.\,180] {Toda}, we know,  for this $\s$  homomorphism, the restriction  $\s|_{\x\af_{2}(8)\xx}$
 is an isomorphism,  and  this  $j_{1*}$ is monomorphic (Lemma\,\,\ref{gjlx}), by \cite[Lemma\, 13.8, p.\,180]{Toda}, we know the restriction $\af_{1}(5)_{*}|_{\x\af_{2}(8)\xx}$ is monomorphic,  then  $\Ker(\pa_{1*_{16}})=0$.\\
 \indent \quad For $\Ker(\pa_{2*_{17}})$, we observe  $\mathrm{BUND}(17,2)$ given by diagram (\ref{bund1}).\\ \indent \quad  In $\mathrm{BUND}(17,2)$:\\ the   homomorphism  $\s$ is an isomorphism,  $j_{2*}$  is monomorphic, by \cite[(13.1)] {Toda}  and \cite[Lemma 13.8, p.\,180]{Toda}, we know $$\af_{1}(6)_{*}(\af_{2}(9))=-3\beta_{1}(6) ,$$thus this $\af_{1}(6)_{*}$ is monomorphic,\;therefore,\; $\Ker(\pa_{2*_{17}})=0$.

\end{itemize}

\end{proof}

\begin{thm}\label{3c15}After localization at 3,

 \[\pi _{15+k}(\Sigma ^{k}\mathbb{H}P^{2})\tg\begin{cases}
\mathbb{Z}/27,\;&\text{$k\geq9\;\mbox{or}\;k\in\n 7,6,5,3     \nn$}\\
 \mathbb{Z}/27    \oplus\mathbb{Z}_{(3)},\;&\text{$k=8\;\mbox{or}\;4$}\\
 \mathbb{Z}/9,\;&\text{$k=2$}\\
 \z/9\jia\mathbb{Z}/3,\;&\text{$k=1$}\\
\mathbb{Z}/3,\;&\text{$k=0$}
\end{cases}\]

\end{thm}
\begin{proof}

\indent For $k\geq 9$,\;$\pi _{15+k}(\Sigma ^{k}\mathbb{H}P^{2})$ are in the stable range, by Theorem\,\ref{3c15w}, we have,\;
$\pi _{15+k}(\Sigma ^{k}\mathbb{H}P^{2})\tg\z/27,\;(k\geq 9).$ \\
\indent For $k=7,6,5\;$or\;3, by Proposition\,\ref{a2g7} we know $\s^{k-3}\overline{\af_{2}}(7)\in\pi _{15+k}(\Sigma ^{k}\mathbb{H}P^{2})$ and  $\ord(\s^{k-3}\overline{\af_{2}}(7))\geq27$. Then, by$\lon (15+k,k)$ given by sequence (\ref{longmk}) and by Lemma\,\,\ref{gjlx},\;we know
$\pi _{15+k}(\Sigma ^{k}\mathbb{H}P^{2})\tg\z/27,$ ($k=7,6,5$\;or\;3), successively,\;$\ord(\s^{k-3}\overline{\af_{2}}(7))=27$,\;($k=7,6,5$\;or\;3).\\
\indent For $k=8$, by the above, $\ord(\overline{\af_{2}}(7))=\ord(\s^{6}\overline{\af_{2}}(7))=27$, so $$\ord(\s^{5}\overline{\af_{2}}(7))=27.$$ Then,\;
$\pi _{23}(\Sigma ^{8}\mathbb{H}P^{2})$ contains $\s^{5}\overline{\af_{2}}(7)$ which is of order 27. By$\lon (23,8)$ given by sequence (\ref{longmk}),\;
we have $\pi _{23}(\Sigma ^{8}\mathbb{H}P^{2})\tg\z/27\jia\z_{(3)}$.\\
\indent For $k=4$, by the proof for $k=7,6,5$\;or\;3,  we know $\ord(\s\overline{\af_{2}}(7))=27$,
$\pi _{19}(\Sigma ^{4}\mathbb{H}P^{2})$ contains $\s^{4}\overline{\af_{2}}(7)$ which is of order 27. By$\lon (19,4)$ given by sequence (\ref{longmk}),
and by Proposition\,\ref{p19f4},  we have $\pi _{19}(\Sigma ^{4}\mathbb{H}P^{2})\tg\z/27\jia\z_{(3)}$.\\
\indent For $k=2$, by Proposition\,\ref{p1617f2}, and by$\sho(17,2)$ given by sequence  (\ref{shortmk}), we have $\pi _{17}(\Sigma ^{2}\mathbb{H}P^{2})\tg\z/9$.\\
\indent For $k=1$, by Proposition\,\ref{p1617f2}, and  by$\sho(16,1)$ given by sequence  (\ref{shortmk}), we have  $\pi _{16}(\Sigma \mathbb{H}P^{2})\tg\pi_{16}(S^{5}\vee S^{13})\tg\z/9\jia\z/3$.

\end{proof}

\begin{cor}\label{ltt} If $k\geq5$ or $k=3$, then  \;$$    \pa_{k_{*15+k}}:   \pi_{15+k}(S^{8+k})\rightarrow \pi_{14+k}(F_{k})$$ is trivial.

\end{cor}

\begin{proof}

For $k\geq5$ or $k=3 $, by$\lon (15+k,k)$ given by sequence (\ref{longmk}), by Theorem\,\ref{3c15}  and by the exactness,\;our corollary  follows.

\end{proof}

\subsection{ Determinations of $\pi_{14+k}(\s^{k}\h P^{2})$ }

\begin{lem}\label{zhd}

$\pi_{18}(F_{4})=\z/3\n j \beta_{1}(8)\nn$.\\

\begin{proof}
In a similar way to the proof of  Proposition\,\ref{p1617f2}(1), we obtain this.

\end{proof}

\end{lem}

 We recall that  $\pa_{k*_{r}}$ is the homomorphism $\pa_{k*_{r}}: \pi_{r}( S^{8+k})\rightarrow \pi_{r-1}(F_{k})$.
\begin{lem}\label{ck143}
$\cok(\pa_{4*_{19}})\tg\cok(\pa_{2*_{17}})\tg\cok(\pa_{1*_{16}})\tg\z/3.$
\end{lem}

\begin{proof}
We notice that
$$\af_{1}(5)\hc\af_{2}(8)=-3\beta_{1}(5)  \text{ (\cite[Lemma 13.8, p.\,180 ]{Toda}),}\; $$
and 
$$ 3\pi_{17+n}(S^{7+n})=0,\;\text{for any}\;n\geq0.$$
Then, by observing the corresponding  $\mathrm{BUND}(m+1,k)$ given by diagram (\ref{bund1}), and by Proposition\,\ref{p1617f2}(1),  Lemma\,\ref{zhd},    the above cokernels are got.
\end{proof}

\begin{thm}After localization at 3,\\

\[\pi _{14+k}(\Sigma ^{k}\mathbb{H}P^{2})\tg\begin{cases}
\z/3,\;&\text{$k\geq4\; \mbox{or}\; k\in\n2, 1\nn$}\\
\mathbb{Z}_{(3)}\jia \z/3,\;&\text{$k=3$}\\
\z/3\jia \z/3,\;&\text{$k=0$}\\
\end{cases}\]

\end{thm}

\begin{proof}
 We recall$\sho(m,k)$ given by sequence (\ref{shortmk}).\\
\indent For $k\geq5$ or $k=3 $,  we use$\sho(14+k,k)$,
 by Corollary\,\ref{ltt}, Lemma\,\,\ref{gjlx},\;we get the result, (here,\;$\pi _{14+k}(\Sigma ^{k}\mathbb{H}P^{2})$ is in the stable range if $k\geq 8$).\\
\indent For $k=4,2\;$or\;1, by$\sho(14+k,k)$ and Lemma\,\ref{ck143}, we obtain the result.

\end{proof}

\section{$\pi_{r+k}(\s^{k}\h P^{2})$ localized at 2}

\indent \;\quad In this section,  \textit{all have been localized at 2} whether we say after localization at 2 or not; and $\widetilde{H}_{*}(-)$  denotes the reduced homology with $\z/2$ coefficients. For an element $\mathbbm{a}\in\pi_{*}(X)$, we usually use $\overline{\mathbbm{a}}$ to denote its lifting up to sign. And $\pi_{n+k}(S^{k})$ localized at 2 is also denoted by  $\pi_{n+k}^{k}$, which  naturally corresponds to  Toda's\, $\pi_{n+k}^{k}$ in  \cite{Toda}. We denote Toda's $E$ by $\s$, the  suspension functor, and we  denote Toda's    $\q$ by $P$, where $\q$ denotes  boundary homomorphisms of   $EHP$ sequence. And in this section, we often freely use the well-known  group structures of  $\pi_{n+k}^{k}$\;$ (n\leq11)$ in  \cite{Toda}. Usually, for a fixed $n$, the family  $\n\pi_{n+k}^{k}\nn_{k\geq2}$ is called the \textit{n-th} stem. For the convenience of  readers, in next paragraph we shall  give a brief introduction of Toda's naming convention  for the  generators of $\pi_{*}(S^{n};\,2\,)$ to help readers to read some of our propositions and their proofs, such a naming convention is also adopted by \cite{20STEM},\cite{2122STEM},\cite{2324STEM},\cite{Oda},\cite{32STEM} and so on. Since Toda's 2-primary component method  in \cite{Toda} is naturally corresponding to the 2-localization method, in this article  we  use the language  of the 2-localization to state the facts on  the 2-primary components.\\\indent After localization at 2, roughly speaking, 
 for $\mathbbm{x}$ 
 which represents a Greek letter,   the symbol $\mathbbm{x}_{n}$ usually  denotes one of the  generators of $\pi_{n+r}(S^{n})$ for some $r$, the  footnote $n$ indicates the codomain of $\mathbbm{x}_{n}$,  and $\mathbbm{x}_{n+k}=\s^{k}\mathbbm{x}_{n}$, $\mathbbm{x}_{n}^{\ell}$ is the abbreviation of $\mathbbm{x}_{n}\hc \mathbbm{x}_{n+r}\hc \cdots \hc \mathbbm{x}_{n+(\ell-1)r}$\,\,($\ell$ composition factors),  and $\mathbbm{x}=\s^{\wq}\mathbbm{x}_{n}$; the usages of  $\overline{\mathbbm{x}}_{n} $ and $\mathbbm{x}^{*}_{n}$ are similar to the above; in $\pi_{j+r}(S^{j})$ (not a stable homotopy group), if a  generator is just denoted by a symbol  without a footnote,  then  this  generator usually cannot survive to the stable homotopy group $\pi_{r}^{S}(S^{0})$ or its  $\s^{\wq}$-\,image is divisible by 2, for example, for  $\theta\in\pi_{24}(S^{12})$,\;$\cg'''\in\pi_{12}(S^{5})$, their  $\s^{\wq}$-\,images satisfy $\s^{\wq}\theta=0$,   $\s^{\wq}\cg'''=8\cg$  is of order 2; there is an interesting advantage of this naming convention, that is, we can observe the commutativity of unstable compositions conveniently,  $$\mathbbm{x}_{n}\hc \mathbbm{y}_{i}=\pm \mathbbm{y}_{n}\hc\mathbbm{x}_{j}  \;\;\text{ for some}\; i, j\;\;\; \text{if} \;\;n\geq a+b$$
 where $\n \mathbbm{x}_{k}\nn$ was born in $\pi_{*}(S^{a})$ and  $\n \mathbbm{y}_{\ell}\nn$ was born in $\pi_{*}(S^{b})$, (see \cite[Prop.\,3.1]{Toda}), for example, for the elements shown in \cite{Toda},  $$\cg_{n}\in\pi_{7+n}(S^{n})\,(n\geq8) \;\text{and} \;\mu_{n}\in\pi_{9+n}(S^{n})\,(n\geq3),$$ we have  $\cg_{8+3}\mu_{i}=\pm \mu_{8+3}\cg_{j}$, successively,\; $\cg_{11}\mu_{18}=\mu_{11}\cg_{20}$, ($\pm$ is not necessary, since $\mu_{3}$ is of order 2), but  $\cg_{10}\mu_{17}\neq\mu_{10}\cg_{19}$, (see \cite[p.\,156]{Toda}).
 Some common names of the generators are summarized in \cite[p.\,189]{Toda} and \cite[(1.1),\,p.\,66]{Oguchi}.

\subsection{ Some relations of generators of $\pi_{*}^{n}$ and Toda brackets    }

\begin{pro}\label{hc1}
\indent $\mathrm{(1)}$\; $P(\iota_{13})\in\n       \nu_{6},\ca_{9},2\iota_{10}  \nn.$\\
\indent $\mathrm{(2)}$\; $\s\cg'''=2\cg''$, \;$\s\cg''=\pm2\cg'$,\; $\s^{2}\cg'=2\cg_{9}.$\\
\indent $\mathrm{(3)}$\; $\cg'''=\n \nu_{5},2\nu_{8},4\iota_{10}             \nn.$\\
\indent $\mathrm{(4)}$\; $\nu_{9}\cg_{12}=2x\cg_{9}\nu_{16}$,\; $(x:\;\mathrm{odd})$,\; $\nu_{11}\cg_{14}=0.$ \\
\indent $\mathrm{(5)}$\; $\nu_{5}\s\cg'=2\nu_{5}\cg_{8}$,\quad$\nu_{5}(\s\cg')\ca_{15}=0.$ \\
\indent $\mathrm{(6)}$\; $\nu_{6}\ch_{9}=\nu_{6}\lt_{9}=2\ch_{6}\nu_{14}$,\; $\nu_{7}\ch_{10}=\nu_{7}\lt_{10}=0.$\\
\indent $\mathrm{(7)}$\; $\ck_{5}\in\n    \nu_{5}, 8\iota_{8},   \s\cg'\nn_{1}.$\\
\indent $\mathrm{(8)}$\; $\ca^{2}\mu=4\ck.$\\
\indent $\mathrm{(9)}$\; $\nu_{6}\ca_{9}=0$.\\
\indent $\mathrm{(10)}$\; $\nu_{6}\cg_{9}\ca_{16}=\nu_{6}\lt_{9}=2\ch_{6}\nu_{14}$.\\
\indent $\mathrm{(11)}$\; $P(\iota_{31})=2\cg_{15}^{2}$.
\end{pro}

\begin{proof}
$\mathrm{(1)}$\; From \cite[ Lemma 5.10, p.\,45]{Toda}.\\
\indent $\mathrm{(2)}$\; From \cite[ Lemma 5.14, p.\,48]{Toda}.\\
\indent $\mathrm{(3)}$\; By \cite[ Proposition\,2.2 (9), p.\,84]{Oda},\\
\centerline{$\cg'''\in\n \nu_{5},2\nu_{8},4\iota_{10}             \nn_{3}\bz\n \nu_{5},2\nu_{8},4\iota_{10}             \nn$},\\
\indent\qquad\qquad\quad$\ind\n \nu_{5},2\nu_{8},4\iota_{10}             \nn=\nu_{5}\hc0+4\pi_{12}^{5}=0$, hence the result.\\
\indent $\mathrm{(4)}$\; The first one is from \cite[(7.19), p.\,71]{Toda}, and the second one is from \cite[(7.20), p.\,72]{Toda}.\\
\indent $\mathrm{(5)}$\; Form \cite[(7.16), p.\,69]{Toda}.\\
\indent $\mathrm{(6)}$\; From \cite[(7.17), p.\,70]{Toda} and  \cite[(7.18), p. 70]{Toda},  $\ch_{7}\nu_{15}\in\pi_{18}^{7}$  is of order 2.\\
\indent $\mathrm{(7)}$\; From \cite[ v), p. 59 ]{Toda}.\\
\indent $\mathrm{(8)}$\; From \cite[ (7.14), p.\,69]{Toda}\\
\indent $\mathrm{(9)}$\; By $\pi_{10}^{6}=0$.\\
\indent $\mathrm{(10)}$  $\nu_{6}\cg_{9}\ca_{14}=\nu_{6}\lt_{9}$ is from \cite[  p.\,152]{Toda}, $\nu_{6}\lt_{9}=2\ch_{6}\nu_{14} $ is just a part of (6) of our proposition.\\
\indent $\mathrm{(11)}$\; From \cite[(10.10), p.\,102]{Toda}
\end{proof}

\begin{pro}\label{lm} After localization at 2,
\begin{itemize}
    \item [\rm (1)]\;$\pi_{7}^{S}(\mathbb{H}$$P^{\infty})=0,\;\pi_{8}^{S}(\mathbb{H}$$P^{\infty})\approx\mathbb{Z}_{(2)},$\;
 $\pi_{9}^{S}(\mathbb{H}$$P^{\infty}$)$\approx\pi_{10}^{S}(\mathbb{H}$$P^{\infty})\approx\mathbb{Z}/2.$
    \item [\rm (2)]$\pi_{11}^{S}(\mathbb{H}$$P^{2})\tg\mathbb{Z}/4\oplus\mathbb{Z}/16,$\\
  $\pi_{13}^{S}(\mathbb{H}P^{2})\approx(\z/2)^{2},\; \;  \pi_{14}^{S}(\mathbb{H}$$P^{2})\approx\z/2,$\\
 $\pi_{15}^{S}(\mathbb{H}$$P^{2})=\mathbb{Z}/128\n \widehat{a}\nn,\quad \widehat{a} \in \langle j,\nu,\sigma        \rangle.$
\end{itemize}
\end{pro}
\begin{proof}
    \begin{itemize}
    \item [\rm (1)] From \cite[Theorem II.9]{Liulevicius}.

    \item [\rm (2)] $\pi_{11}^{S}(\y)$ is from  \cite[p.\,198]{Mukai}, notice that $\pi_{14}^{S}(Q_{2}^{3})=\pi_{11}^{S}(\y)$.\\
     $\pi_{13}^{S}(\y)$, $\pi_{14}^{S}(\y)$ are from  \cite[Lemma 2.7]{Mukai}.\\$\pi_{15}^{S}(\y)$ is from  \cite[Lemma 3.3]{Mukai}.
    \end{itemize}
\end{proof}

\subsection{ Determinations of $\pi_{r+k}(\s^{k}\h P^{2})\;\;(7\leq r\leq10)$ }

 We recall that  $\pa_{k*_{r}}$ is the homomorphism $\pa_{k*_{r}}: \pi_{r}( S^{8+k})\rightarrow \pi_{r-1}(F_{k})$.
\begin{lem}\label{ck1}
\indent $\mathrm{(1)}$ \quad$\cok(\pa_{1*_{10}})=0$,\;$\Ker( \pa_{1*_{9}})=\z_{(2)}\n     8\iota_{8}  \nn.$\\
\indent $\mathrm{(2)}$\quad$\cok(\pa_{1*_{11}})=\Ker( \pa_{1*_{10}})=0.$\\
\indent $\mathrm{(3)}\quad\cok(\pa_{1*_{12}})=\cok(\pa_{2*_{13}})=\cok(\pa_{3*_{14}})=\Ker( \pa_{1*_{11}})=0,$\\
\indent\qquad \;\;$\Ker( \pa_{2*_{12}})\tg\Ker( \pa_{3*_{13}})\tg\z/2.$
\end{lem}

\begin{proof}
These cokernels and kernels are immediately got by  observing the corresponding diagrams $\mathrm{BUND}(m+1,k)$  given by diagram (\ref{bund1}) and the group structures of $\pi_{*}^{n}$ in \cite{Toda}.
\end{proof}

\begin{pro}\label{mzf} After localization at 2, $$\pi_{11}(\s^{2}\h P^{2})\tg\z_{(2)}.$$

\end{pro}

\begin{proof}
By Proposition\,\ref{grdl}(1)(2), we have the fibre sequence 
$$\lo\s^{2}\h P^{2}    \stackrel{d} \rightarrow  J(M_{S^{9}}, \s\y)  \rightarrow  S^{6}\stackrel{j_{2}} \longrightarrow\s^{2}\h P^{2}$$
where $\sk_{13}(J(M_{S^{9}}, \s\y))=S^{9},$ it  induces the following exact sequence  (here, the  homomorphism can be replaced by $\nu_{6_{*}}$, because for the cofiber sequence $S^{9}\stackrel{\nu_{6}}\longrightarrow S^{6}\stackrel{j_{2}}\longrightarrow \s^{2}\h P^{2}$, $\nu_{6}$ can extend to the fiber of $j_{2}$; another way, we can use the James theorem, that is, \cite[Theorem 2.1]{Ja1}),\\\\
\indent\qquad \xymatrix@C=0.5cm{&
 \pi_{11}^{9} \ar[r]^{\nu_{6_{*}}} & \pi_{11}^{6} \ar[r]^{j_{2*}\qquad}& \pi_{11}(\s^{2}\h P^{2}) \ar[r] ^{\qquad d_{*}}& \pi_{10}^{9} \ar[r]& 0.  }

\noindent\\ By  $\nu_{6_{*}}(\pi_{11}^{9})=\x \nu_{6}\ca_{9}^{2}\xx=0$ and Lemma\,\,\ref{xwjd}, we have the exact sequence,\\\\
\indent\quad \xymatrix@C=0.5cm{&&
0 \ar[r] & \pi_{11}^{6} \ar[r]^{j_{2*}\qquad}& \pi_{11}(\s^{2}\h P^{2}) \ar[r] ^{\qquad\;p_{*}}& \pi_{11}^{10} \ar[r]& 0  }

\noindent \\we know (\cite{Toda}), $$\pi_{11}^{6}= \z_{(2)}\n P(\iota_{13})\nn,\;\;\pi_{11}^{10}= \z/2\n\eta_{10}\nn.$$
 By Proposition\,\ref{hc1}(1),\; $P(\iota_{13})\in\n \nu_{6},\eta_{9},2\iota_{10}\nn$,  so \;$$ j_{2*}(P(\iota_{13}))\in  j_{2}\circ\n \nu_{6},\eta_{9},2\iota_{10} \nn=-(\n   j_{2},\nu_{6},\eta_{9}  \nn\hc 2\iota_{11})=-2\n   j_{2},\nu_{6},\eta_{9}  \nn.$$
Thus $ j_{2*}(P(\iota_{13}))$ can be divisible by 2. Hence,  $\pi_{11}(\s^{2}\h P^{2})\tg\z_{(2)}$.

\end{proof}

\begin{thm}After localization at 2,  \[\pi _{7+k}(\Sigma ^{k}\mathbb{H}P^{2})\tg\begin{cases}
0,\;&\text{$k\geq1$}\\
\z/4,\;&\text{$k=0$}\end{cases}\]
   \[\pi _{8+k}(\Sigma ^{k}\mathbb{H}P^{2})\tg\begin{cases}
\mathbb{Z}_{(2)},\;&\text{$k\geq1$}\\
\z/2,\;&\text{$k=0$}\end{cases}\]
    \[\pi _{9+k}(\Sigma ^{k}\mathbb{H}P^{2})\tg\begin{cases}
\mathbb{Z}/2,\;&\text{$k\geq3\;\mbox{or}\;k=0$}\\
\z_{(2)},\;&\text{$k=2$}\\
0,\;&\text{$k=1$}\end{cases}\]
\[\pi _{10+k}(\Sigma ^{k}\mathbb{H}P^{2})\tg\begin{cases}
\mathbb{Z}/2,\;&\text{$k\geq2$}\\
0,\;&\text{$k=1\;or\;0$}
\end{cases}\]

\end{thm}

\begin{proof}
We recall$\sho(m,k)$ given by sequence (\ref{shortmk}), and$\lon(m,k)$ given by sequence (\ref{shortmk}).\\
\indent For  $\pi _{7+k}(\Sigma ^{k}\mathbb{H}P^{2})$,\;in the case  $k\geq1$,  $\pi _{7+k}(\Sigma ^{k}\mathbb{H}P^{2})$ are in the stable range, so,$$\pi _{7+k}(\Sigma ^{k}\mathbb{H}P^{2})\tg\pi _{7}^{S}(\mathbb{H}P^{2}) \tg \pi _{7}^{S}(\mathbb{H}P^{\wq}),\;(k\geq1),$$ by Proposition\,\ref{lm}(1) we get the result.\\\indent For $\pi _{8+k}(\Sigma ^{k}\mathbb{H}P^{2})$,\;in the case\;$k\geq2$,  $\pi _{8+k}(\Sigma ^{k}\mathbb{H}P^{2})$ are in the stable range, so,\;
 $$\pi _{8+k}(\Sigma ^{k}\mathbb{H}P^{2})\tg\pi _{8}^{S}(\mathbb{H}P^{2}) \tg \pi _{8}^{S}(\mathbb{H}P^{\wq}),\;(k\geq2),$$ by Proposition\,\ref{lm}(1) we get the results for $k\geq2$.\;In the case  $k=1$, we consider$\sho(9,1)$, by Lemma\,\ref{ck1}(1) we get the result.\\
\indent For $\pi _{9+k}(\Sigma ^{k}\mathbb{H}P^{2})$, \;if $k\geq3$,  $\pi _{9+k}(\Sigma ^{k}\mathbb{H}P^{2})$ are in the stable range, so,\;
 $$\pi _{9+k}(\Sigma ^{k}\mathbb{H}P^{2})\tg\pi _{9}^{S}(\mathbb{H}P^{2}) \tg \pi _{9}^{S}(\mathbb{H}P^{\wq}),\;(k\geq3),$$ 
 by Proposition\,\ref{lm}(1) we get the result for   $k\geq3$. \;For $k=2$,  it's just Proposition\,\ref{mzf}. For $k=1$, we consider$\sho(10,1)$, and by Lemma\,\,\ref{ck1}(2),  we get the result.\\
\indent For $\pi _{10+k}(\Sigma ^{k}\mathbb{H}P^{2})$,\;if $k\geq4$, they are in the stable range, by Proposition\,\ref{lm}(1) we get the results for   $k\geq4$. For $k\in\n 3,2,1\nn$, we consider$\sho(10+k,k)$, by Lemma\,\ref{ck1}(3) we get the results.

\end{proof}

\subsection{ Determinations of $\pi_{11+k}(\s^{k}\h P^{2})$ }

 We recall that  $\pa_{k*_{r}}$ is the homomorphism $\pa_{k*_{r}}: \pi_{r}( S^{8+k})\rightarrow \pi_{r-1}(F_{k})$.
\begin{lem}\label{ck11}
\indent $\mathrm{(1)}$ \quad$\cok(\pa_{3*_{15}})=\z/8\n j\cg'\nn$,\;$\Ker( \pa_{3*_{14}})=\z/4\n     2\nu_{11}  \nn.$\\
\indent $\mathrm{(2)}$ \quad$\cok(\pa_{4*_{16}}) =\z_{(2)}\n j\cg_{8}\nn\jia\z_{8}\n j\s\cg'\nn      ,\;\;                                     \Ker( \pa_{4*_{15}})=\z/4\n     2\nu_{12}  \nn.$\\
\indent $\mathrm{(3)}$ \quad$\cok(\pa_{2*_{14}}) =\z/4\n j\cg''\nn      ,\;\;                                     \Ker( \pa_{2*_{13}})=\z/4\n     2\nu_{10}  \nn.$\\
\indent $\mathrm{(4)}$ \quad$\cok(\pa_{1*_{13}}) =\z/2\n j\cg'''\nn      ,\;\;                                     \Ker( \pa_{1*_{12}})=\z/4\n     2\nu_{9}  \nn.$\\

\end{lem}

\begin{proof}

 They are immediately got by the corresponding  diagrams  $\mathrm{BUND}(m+1,k)$  given by  diagram (\ref{bund1}), Proposition\,\ref{hc1}(2) and the corresponding   group structures of  $\pi_{*}^{n}$ shown in \cite{Toda}.

\end{proof}

\begin{thm}\label{2c11aa} After localization at 2,
\[\pi _{11+k}(\Sigma ^{k}\mathbb{H}P^{2})\tg\begin{cases}
\mathbb{Z}/16\oplus\mathbb{Z}/4,\;&\text{$k\geq5$}\\
\mathbb{Z}/8\oplus\mathbb{Z}/4\oplus\mathbb{Z}_{(2)},\;&\text{$k=4$}\\
\mathbb{Z}/8\oplus\mathbb{Z}/4,\;&\text{$k=3$}\\
\mathbb{Z}/8\oplus\mathbb{Z}/2,\,\;&\text{$k=2$}\\
\mathbb{Z}/8,\;&\text{$k=1$}\\
\mathbb{Z}_{(2)},\;&\text{$k=0$}
\end{cases}\]

\end{thm}

\begin{proof}
For $k\geq5$,  $\pi _{11+k}(\Sigma ^{k}\mathbb{H}P^{2})$ are in the stable range. Hence,
 $\pi _{11+k}(\Sigma ^{k}\mathbb{H}P^{2})\tg\pi _{11}^{S}(\mathbb{H}P^{2}) ,\;(k\geq5)$, by Proposition\,\ref{lm}(2) we get the result.\\
 \indent For  $k=3$,  we consider$ \sho(14,3)$ given by sequence (\ref{shortmk}),  by Lemma\,\,\ref{ck11}(1), we have the exact sequence,\\
 \xymatrix@C=0.3cm{
 &&&&  0 \ar[r] & \z/8\n j\cg' \nn \ar[rr]^{} && \pi_{14}(\s^{3}\y) \ar[rr]^{p_{3*}} && \z/4\n     2\nu_{11} \nn \ar[r] & 0. }

 \noindent Since\;$$p_{3}\hc \n    j, \nu_{7},2\nu_{10}      \nn=-\n   p_{3},    j, \nu_{7}     \nn\hc2\nu_{11}\ni-2\nu_{11},$$
 then $$\cz \mathbbm{x}\in\n    j, \nu_{7},2\nu_{10}      \nn, \text{\;such that}\; p_{3*}(\mathbbm{x})=-2\nu_{11},$$
 then \;$$ 4\mathbbm{x}\in\n    j, \nu_{7},2\nu_{10}      \nn\hc4\iota_{14}=-j\hc\n \nu_{7},2\nu_{10}     ,4\iota_{13}\nn.$$
 By Proposition\,\ref{hc1}(2)(3),\;$$4\cg'\in\n \nu_{7},2\nu_{10}     ,4\iota_{13}\nn\m 4\cg',$$  that is,\; $$\n \nu_{7},2\nu_{10}     ,4\iota_{13}\nn=\z/2\n 4\cg' \nn.$$
  Thus,\; $$- j\hc\n \nu_{7},2\nu_{10}     ,4\iota_{13}\nn= \z/2\n 4j\cg' \nn.$$ Hence, $\cz t\in\z,$ such that $4\mathbbm{x}=(4t)j\cg'$. Let $\mathbbm{x}'=\mathbbm{x}-tj\cg'$,
 so\;$4\mathbbm{x}'=0, $  $p_{3*}(\mathbbm{x}')=-2\nu_{11}$. Thus,  $\mathbbm{x}'$ is of order 4, then,\;$$\pi_{14}(\s^{3}\y)= \z/8\n j\cg' \nn\jia\z/4\n   \mathbbm{x}'     \nn.$$
  \indent  For $k=4$, We consider $ \sho(15,4)$ given by sequence (\ref{shortmk}),  by Lemma\,\,\ref{ck11}(2), we have the exact sequence,\\\\
 \xymatrix@C=0.3cm{
 &  0 \ar[r] & \z/8\n j\s\cg' \nn\jia\z_{(2)}\n j\cg_{8}\nn \ar[rr]^{} && \pi_{15}(\s^{4}\y) \ar[rr]^{p_{4*}} && \z/4\n     2\nu_{12} \nn \ar[r] & 0. }

  \noindent\\ Thus,\, $\pi_{15}(\s^{4}\y)\tg\z_{(2)}\jia A$,  where $$A= \z/8,\;\z/16,\;\z/8\jia\z/2,\;\z/8\jia\z/4,\;\z/16\jia\z/2\;\mbox{or}\;\z/32.$$We consider$\sho(14,3,1)$  given by diagram (\ref{shortmkt}),\;by $\mathrm{COM}(14,3,1)$ given by diagram (\ref{com1}), we know the homomorphism  $\mathcal{A}(14,3,1)$ in$ \sho(14,3,1)$ is monomorphic,  and obviously $\mathcal{C}(14,3,1)$ is monomorphic. By the snake lemma, we know $ \s :\pi_{14}(\s^{3}\y)\rightarrow\pi_{15}(\s^{4}\y)$ is monomorphic, finally, by the result of $\pi_{14}(\s^{3}\y)$ we know\; $$\pi_{15}(\s^{4}\y)\tg\z_{(2)}\jia \z/8\jia\z/4.$$
   \indent  For $k=2$, We consider$ \sho(13,2)$ given by sequence (\ref{shortmk}),  by Lemma\,\,\ref{ck11}(3), we have the exact sequence,\\\\
 \xymatrix@C=0.3cm{
 &&&&  0 \ar[r] & \z/4\n j\cg'' \nn \ar[rr]^{} && \pi_{13}(\s^{2}\y) \ar[rr]^{p_{2*}} && \z/4\n     2\nu_{10} \nn \ar[r] & 0. }

 \noindent\\ Thus,\, $\pi_{13}(\Sigma^{2}\mathbb{H}P^{2})\approx \mathbb{Z}/16 \;\mbox{or} \;\mathbb{Z}/8\oplus\mathbb{Z}/2 \;\mbox{or} \; \z/4\jia\z/4.$\;\\
 Since $$p_{2}\hc \n    j, \nu_{6},2\nu_{9}      \nn=-\n   p_{2},    j, \nu_{6}     \nn\hc2\nu_{10}\ni-2\nu_{10},$$
so there exists $ \mathbbm{y}\in\n    j, \nu_{7},2\nu_{10}      \nn$, such that\; $p_{2*}(\mathbbm{y})=-2\nu_{10}$. By Proposition\,\ref{hc1}(2)\,(3),$$ 4\mathbbm{y}\in\n    j, \nu_{6},2\nu_{9}      \nn\hc4\iota_{13}=- (j\hc\n \nu_{6},2\nu_{9}     ,4\iota_{12}\nn)\ni2j\cg''\m j\hc(\nu_{6}\hc0+4\pi_{13}^{6})=0,$$  so $4\mathbbm{y}=2j\cg''$  is of order 2.  Then  $\mathbbm{y}$ is of order 8, thus\; $\pi_{13}(\Sigma^{2}\mathbb{H}P^{2})
 \tg\mathbb{Z}/8\oplus\mathbb{Z}/2$.  We notice $$p_{2*}(2\mathbbm{y}+j\cg'')=p_{2*}(2\mathbbm{y})=-4\nu_{10}\neq0,\;2(2\mathbbm{y}+j\cg'')=4\mathbbm{y}+2j\cg''=2j\cg''+2j\cg''=0,$$
hence, $\ord(2\mathbbm{y}+j\cg'')=2$. We show
 $\x \mathbbm{y}\xx+\x2\mathbbm{y}+j\cg''\xx$ is a direct sum,
set   $$a\mathbbm{y}+b(2\mathbbm{y}+j\cg'')=0 \;\text{for some}\; a,b\in\z.$$ Applying $p_{*}$, we have $a+2b\ty0\m4$; multiplying by 2, we have $2a\ty0\m8$, then, $a\ty0\m4$, combining with the above, we know $2b\ty0\m4$, so $b\ty0\m2$, successively $a\ty0\m8$, so the sum is a direct sum.
Therefore,\; $\pi_{13}(\Sigma^{2}\mathbb{H}P^{2})=\mathbb{Z}/8\n\mathbbm{y}\nn\oplus\mathbb{Z}/2\n 2\mathbbm{y}+j\cg''\nn$.\\
   \indent  For $k=1$, we consider$ \sho(12,1)$ of sequences (\ref{shortmk}), by Lemma\,\,\ref{ck11}(4), we have the exact sequence,\\\\
 \xymatrix@C=0.3cm{
 &&&&  0 \ar[r] & \z/2\n j\cg''' \nn \ar[rr]^{} && \pi_{12}(\s\y) \ar[rr]^{p_{1*}} && \z/4\n     2\nu_{9} \nn \ar[r] & 0. }

 \noindent\\ So,\, $\pi_{12}(\Sigma\mathbb{H}P^{2})\approx \mathbb{Z}/8 \;\mbox{or} \; \z/4\jia\z/2.$\;\\
 Since $$p_{1}\hc \n    j, \nu_{5},2\nu_{8}      \nn=-\n   p_{1},    j, \nu_{5}     \nn\hc2\nu_{9}\ni-2\nu_{9},$$
 then $$\cz \mathbbm{a}\in\n    j, \nu_{5},2\nu_{8}      \nn,$$ such that\; $p_{1*}(\mathbbm{a})=-2\nu_{9}$. By Proposition\,\ref{hc1}(3),\;$$ 4\mathbbm{a}\in\n    j, \nu_{5},2\nu_{8}      \nn\hc4\iota_{12}=- j\hc\n \nu_{5},2\nu_{8}     ,4\iota_{11}\nn= j\cg'''.$$
 So \;$4\mathbbm{a}=j\cg'''$  is of order 2.   Then $\mathbbm{a}$ is of order 8,  therefore,\,$\pi_{12}(\Sigma\mathbb{H}P^{2})\approx \mathbb{Z}/8$.

\end{proof}

\subsection{  Determinations of $\pi_{12+k}(\s^{k}\h P^{2})$}

\begin{lem}\label{p14f2} $ \pi_{14}(F_{2})\tg(\z/2)^{2}.$

\end{lem}

\begin{proof}

By Lemma\,\,\ref{gjlx},  $\sk_{23}(F_{2})=S^{6}\cup _{2\ch_{6}}e^{15}$, \;by \cite{Toda}, \;$$\pi_{14}^{6}=\z/8\n \ch_{6} \nn\jia\z/2\n\lt_{6}\nn.$$
Hence,\,$ \pi_{14}(F_{2})\tg\pi_{14}(   S^{6}\cup _{2\ch_{6}}e^{15}     )\tg(\z/8\n \ch_{6} \nn\jia\z/2\n\lt_{6}\nn)/(\z/4\n 2\ch_{6}\nn)\tg(\z/2)^{2}$.

\end{proof}

\begin{thm} After localization at 2,
\[\pi _{12+k}(\Sigma ^{k}\mathbb{H}P^{2})\tg\begin{cases}
(\mathbb{Z}/2)^{2},\;&\text{$k\geq6\; \mbox{or}\;k\in\n2, \;0\nn$}\\
(\mathbb{Z}/2)^{3},\;&\text{$k=5\;\mbox{or}\;3$}\\
(\mathbb{Z}/2)^{4},\;&\text{$k=4$}\\
\mathbb{Z}_{(2)}\oplus\mathbb{Z}/2,\;&\text{$k=1$}\\
\end{cases}\]
\end{thm}
\begin{proof}

We recall$\lon (m,k)$ given by sequence (\ref{longmk}).\\
\indent For $k\geq6$,\; $\pi _{12+k}(\Sigma ^{k}\mathbb{H}P^{2})$ are in the stable range.  Hence,
 $$\pi _{12+k}(\Sigma ^{k}\mathbb{H}P^{2})\tg\pi _{12}^{S}(\mathbb{H}P^{2}) ,\;(k\geq6).$$
We consider$ \lon(12+k,k)$,  since $\pi_{4}^{S}(S^{0})=\pi_{5}^{S}(S^{0})=0$, so,
  $$\pi _{18}(\Sigma ^{6}\mathbb{H}P^{2})\tg\pi_{18}(F_{6})\tg\pi_{18}^{10}\tg\z/2\jia\z/2.$$
 \indent  For   $1\leq k\leq5$, We consider$ \lon(12+k,k)$,  since $\pi_{4}^{S}(S^{0})=\pi_{5}^{S}(S^{0})=0$, so \\  $\pi _{12+k}(\Sigma ^{k}\mathbb{H}P^{2})\tg\pi_{12+k}(F_{k}),\,(1\leq k\leq5) $,   by Lemma\,\,\ref{p14f2} and Lemma\,\,\ref{gjlx} we get the result.

\end{proof}

\subsection{ Determinations of  $\pi_{13+k}(\s^{k}\h P^{2})$ }

\begin{lem} \label{p15f22} $\pi_{15}(F_{2})=\z/2\n j\nu_{6}^{3}  \nn\jia\z/2\n j \ca_{6}\lt_{7}  \nn\jia\z/2\n  j\mu_{6}\nn\jia\z_{(2)}\n j\hc\overline{ 4\iota_{15} } \nn.$

\end{lem}

\begin{proof}
 We consider $\mathrm{PISK}(15,2)$ given by sequence (\ref{pisk1}),\\\\
\xymatrix@C=0.49cm{
 &\pi_{16}^{15} \ar[r]^{d_{*}\qquad\;\quad} &  \pi_{15}(J(M_{S^{6}},S^{14})) \ar[r]^{} & \pi_{15}(\sk_{23}(F_{2})) \ar[r]^{} &\pi_{15}^{15}\ar[r]^{d_{*}\qquad\;\quad} &   \pi_{14}(J(M_{S^{6}},S^{14})),  }

 \noindent\\ by $\mathrm{D}(16,2)$ given by diagram (\ref{dkm}),  we have $$ \cok(d_{*} :\pi_{16}^{15}  \rightarrow\pi_{15}(J(M_{S^{6}},S^{14})))= \z/2\n j\nu_{6}^{3}  \nn\jia\z/2\n j \ca_{6}\lt_{7}  \nn\jia\z/2\n  j\mu_{6}\nn.$$  By $\mathrm{D}(15,2)$,  we have $$ \Ker(d_{*}:\pi_{15}^{15} \rightarrow \pi_{14}(J(M_{S^{6}},S^{14})))=\z_{(2)}\n     4\iota_{15}  \nn.$$  Then, $$ \pi_{15}(\sk_{23}(F_{2})) =\z/2\n j\nu_{6}^{3}  \nn\jia\z/2\n j \ca_{6}\lt_{7}  \nn\jia\z/2\n  j\mu_{6}\nn\jia\z_{(2)}\n \overline{ 4\iota_{15} } \nn,$$ hence the result follows.

\end{proof}

 We recall that  $\pa_{k*_{r}}$ is the homomorphism $\pa_{k*_{r}}: \pi_{r}( S^{8+k})\rightarrow \pi_{r-1}(F_{k})$.

\begin{lem} \label{ck13}

$\cok(\pa_{6*_{20}})\tg (\z/2)^{2}\jia\z_{(2)}$,\quad
$\cok(\pa_{5*_{19}})\tg (\z/2)^{3}$, \\
$\cok(\pa_{4*_{18}})\tg (\z/2)^{4}$,\quad
$\cok(\pa_{3*_{17}})\tg (\z/2)^{3}$, \\
$\cok(\pa_{2*_{16}})\tg (\z/2)^{2}\jia\z_{(2)}$,\quad
$\cok(\pa_{1*_{15}})\tg (\z/2)^{3}$.
\end{lem}

\begin{proof}
For  $\cok(\pa_{r*_{r+14}}) ,\;(r\neq2)$ in this lemma, We consider $\mathrm{BUND}(r+14,r)$  given by    diagram (\ref{bund1}), we obtain the result.\\
\indent For  $\cok(\pa_{2*_{16}})$, by Lemma\,\,\ref{p15f22} and  the diagram $\mathrm{BUND}(16,2)$, we get the result.

\end{proof}
\begin{thm} After localization at 2,\[\pi _{13+k}(\Sigma ^{k}\mathbb{H}P^{2})\tg\begin{cases}
(\z/2)^{2},\;&\text{$k\geq7$}\\
 (\z/2)^{2}   \oplus\mathbb{Z}_{(2)},\;&\text{$k=6\;\mbox{or}\;2$}\\
(\z/2)^{3},\;&\text{$k=5,\;3,\;1\;\mbox{or}\;0$}\\
(\z/2)^{4},\,\;&\text{$k=4$}\\
\end{cases}\]
\end{thm}
\begin{proof}

\indent For $k\geq7$,  $\pi _{13+k}(\Sigma ^{k}\mathbb{H}P^{2})$ are in the stable range.  Thus,\,
 $\pi _{13+k}(\Sigma ^{k}\mathbb{H}P^{2})\tg\pi _{13}^{S}(\mathbb{H}P^{2}) ,\;(k\geq7)$,  by  Proposition\,\ref{lm}(2) we get the result.\\
 \indent For  $1\leq k\leq6$, by use of $ \lon(13+k,k)$ given by sequence (\ref{longmk}),  since $\pi_{5}^{S}(S^{0})$ is trivial, then$$\pi _{13+k}(\Sigma ^{k}\mathbb{H}P^{2})\tg \cok(\pa_{k*_{14+k}}),\,(1\leq k\leq6),$$   by Lemma\,\,\ref{ck13}  we get the results.

\end{proof}

\subsection{ Determinations of  $\pi_{14+k}(\s^{k}\h P^{2})$ }

\begin{lem} \label{p16f22}$\mathrm{(1)}$\; $\pi_{16}(F_{2})=\z/8\n  j\nu_{6}\cg_{9}    \nn\jia\z/2\n j \ca_{10}\mu_{11}    \nn\jia\z/2\n     j\hc \overline{\ca_{15}}  \nn$,\;\\
$\mathrm{(2)}$\; $\pi_{18}(F_{4})=\z/8  \n     j\cg_{8}\nu_{15}  \nn\jia\z/2\n j\ca_{8}\mu_{9}\nn$,\; $j\nu_{8}\cg_{11}=2t'\cg_{8}\nu_{15}$, $(t': \mathrm{odd})$.
\end{lem}

\begin{proof}
\begin{itemize}
    \item [\rm(1)]
By Lemma\,\ref{gjlx}, \;$\sk_{23}(F_{2})=S^{6} \cup_{2\ch_{6}} e^{15}$. We consider $\mathrm{PISK}(17,2)$ given by sequence (\ref{pisk1}), and we observe \;$\mathrm{D}(18,2)$ and $\mathrm{D}(17,2)$ given by diagram (\ref{dkm}),
 we know $d_{*}( \pi_{17}^{15})=d_{*}( \pi_{16}^{15})=0$.  Since $$ \s (S^{6} \cup_{2\ch_{6}} e^{15})=\s (S^{6} \cup_{\pm[\iota_{6},\nu_{6}]} e^{15})=S^{7}\vee S^{16},$$ we have the homotopy-commutative diagram with rows fibre sequences, \\\\
\xymatrix@C=1.2cm{
   & J(M_{S^{6}},S^{14})  \ar[d]^{} \ar[r]^{}& \sk_{23}(F_{2})     \ar[d]_{\lo\s} \ar[r]^{} & S^{15} \ar[d]^{\lo\s} \\
&  \lo J(M_{S^{7}},S^{15})  \ar[r]^{} &\lo(S^{7}\vee S^{16})  \ar[r]^{ }& \lo S^{16}  }

\noindent  then we have the commutative diagram with exact rows, (by the Hilton-Milnor formula,  the second row splits),\\\\
\xymatrix@C=0.9cm{
&0\ar[r] & \pi_{16}^{6} \ar[d]_{\s} \ar[r]^{} &\pi_{16}(\sk_{23}(F_{2}))   \ar[d]_{\s} \ar[r]^{} &  \pi_{16}^{15} \ar[d]_{\s}^{\tg} \ar[r] & 0  \\
  & 0 \ar[r] & \pi_{17}^{7}  \ar[r]^{}  & \pi_{17}^{7}  \jia\pi_{17}^{16}   \ar[r]^{} & \pi_{17}^{16} \ar[r] & 0   }

\noindent  By \cite[Theorem\,7.3, p.\,66]{Toda} we know\; $\pi_{16}^{6}\stackrel{\s}\longrightarrow\pi_{17}^{7}$ is isomorphic. By the Five Lemma, \\ $$\pi_{16}(\sk_{23}(F_{2}))\stackrel{\s}\longrightarrow \pi_{17}^{7}  \jia\pi_{17}^{16}$$ is  isomorphic,  and the first row  also splits, then our lemma is established.
    \item [\rm(2)] Similarly to the proof of Lemma\,\,\ref{p14f2}, we obtain the result.

\end{itemize}

\end{proof}

 We recall that  $\pa_{k*_{r}}$ is the homomorphism $\pa_{k*_{r}}: \pi_{r}( S^{8+k})\rightarrow \pi_{r-1}(F_{k})$.
\begin{lem} \label{ck14}

\indent$\mathrm{(1)}$\quad$\cok(\pa_{7*_{22}})\tg (\z/2)^{2}$. \quad\\
\indent$\mathrm{(2)}$\quad$\cok(\pa_{6*_{21}})\tg (\z/2)^{2}$, \quad
$\cok(\pa_{5*_{20}})\tg (\z/2)^{2}.$\\
\indent$\mathrm{(3)}$\quad$\cok(\pa_{4*_{19}})\tg (\z/2)^{2}$, \quad
$\cok(\pa_{3*_{18}})\tg \z/2\jia\z_{(2)}$, \\
\indent\qquad\;$\cok(\pa_{2*_{17}})\tg (\z/2)^{2}.$\\
\indent$\mathrm{(4)}$\quad$\cok(\pa_{1*_{16}})\tg (\z/2)^{2}.$\\
\indent$\mathrm{(5)}$\quad$\Ker (\pa_{r*_{r+14}})=0,\; (1\leq r\leq 7).$

\end{lem}

\begin{proof}

We recall that $\mathrm{BUND}(m+1,k)$ is given by diagram (\ref{bund1}).\\
\indent$\mathrm{(1)}$\quad             We consider  $\mathrm{BUND}(22,7)$, notice  $\nu_{11}\cg_{14}=0$  in Proposition\,\ref{hc1}(4), then we get the result..\\
\indent$\mathrm{(2)}$\quad  We consider the corresponding
$\mathrm{BUND}(r+15,r)$, by   $\nu_{9}\cg_{12}=2x\cg_{9}\nu_{16}\;(x:\,\mathrm{odd})$  (Proposition\,\ref{hc1}(4)), we get the result.\\
\indent$\mathrm{(3)}$\quad We consider the corresponding $\mathrm{BUND}(r+15,r)$  and use Lemma\,\,\ref{p16f22}, we get the result..\\
\indent$\mathrm{(4)}$\quad  We consider $\mathrm{BUND}(16,1)$  and use Proposition\,\ref{hc1}(5), we obtain the result.\\
\indent$\mathrm{(5)}$\quad We consider the corresponding $\mathrm{BUND}(r+14,r)$ and use Lemma\,\ref{p15f22},  we obtain the result.

\end{proof}

\begin{thm} After localization at 2,
\[\pi _{14+k}(\Sigma ^{k}\mathbb{H}P^{2})\tg\begin{cases}
\mathbb{Z}/2,\;&\text{$k\geq8$}\\
(\mathbb{Z}/2)^{2},\;&\text{$k=7,\;6\;,5,\;4,\;2\;\mbox{or}\;1$}\\
\z/2\jia\z_{(2)},\;&\text{$k=3$}\\
\z/8\jia\mathbb{Z}/4\oplus\mathbb{Z}/2,\;&\text{$k=0$}
\end{cases}\]
\end{thm}
\begin{proof}

\indent For $k\geq8$,  $\pi _{14+k}(\Sigma ^{k}\mathbb{H}P^{2})$ are in the stable range. Hence,
 $\pi _{14+k}(\Sigma ^{k}\mathbb{H}P^{2})\tg\pi _{14}^{S}(\mathbb{H}P^{2}) ,\;(k\geq8)$, by Proposition\,\ref{lm}(2) we get the result.\\
 \indent For $1\leq k\leq7$, by $ \sho(14+k,k)$ given by sequence (\ref{shortmk}), and by Lemma\,\ref{ck14}, we get the results, (in this case, all the $\z_{(2)}$-module extension problems are trivial).

\end{proof}

\subsection{ Determinations of $\pi_{15+k}(\s^{k}\h P^{2})$ }

\begin{lem}\label{zc15ys}
\begin{itemize}
\item[\rm(1)]  \; $\cz \overline{\cg_{15}}\in\n   j_{7},\nu_{11},\cg_{14}                           \nn\bz\pi_{22}(\s^{7}\y)$,\;\\and\, $\ord(\s^{n}\,\overline{\cg_{15}}\;)\geq128,\;\,\forall n\geq 0.$
\item[\rm(2)]\;$\cz \overline{2\cg_{14}}\in\n   j_{6},\nu_{10},2\cg_{13}                           \nn\bz\pi_{21}(\s^{6}\y)$,\;and\,  $ \ord(\s^{n}\,\overline{2\cg_{14}}\;)\geq64,\;\forall n\geq 0.$
\item[\rm(3)]\;$\cz \overline{4\cg_{13}}\in\n   j_{5},\nu_{9},4\cg_{12}                           \nn\bz\pi_{20}(\s^{5}\y)$,\;and\,  $ \ord(\overline{4\cg_{13}})\geq32.$
\item[\rm(4)]\;$\cz \overline{8\cg_{9}}\in\n   j_{1},\nu_{5},8\cg_{8}                           \nn\bz\pi_{16}(\s\y)$,\;and\, $ \ord(\s^{n}\,\overline{8\cg_{9}}\;)\geq16,\;\forall n\geq 0.$
\end{itemize}
\end{lem}
\begin{proof}

\begin{itemize}
\item[\rm(1)]    By Proposition\,\ref{hc1}(4), we know the Toda bracket  $ \n   j_{7},\nu_{11},\cg_{14}                           \nn$ is well-defined, we take  $\overline{\cg_{15}}$ \,from it.
  By Proposition\,\ref{lm}(2), we have$$\s^{\wq}\,\overline{\cg_{15}}\in \pm \x  j,\nu,\cg                         \xx\ni \pm\hat{a} $$
  $\text{mod} \;j\hc \pi_{11}^{S}(S^{0})+\pi_{8}^{S}(\y)\hc \cg\bz 8\pi_{15}^{S}(\y)=\x         8\hat{a} \xx.$ Thus\;  $$ \s^{\wq}\,\overline{\cg_{15}}\;\ty \pm\hat{a}\m8\hat{a}.$$ Since $ \ord(\hat{a})=128$, then,\,$\s^{\wq}\,\overline{\cg_{15}}$ is also of order 128, hence the result holds.

\item[\rm(2)] By Proposition\,\ref{hc1}(4)  and  $\cg_{10}\nu_{17}\in\pi_{20}^{10}$ is of order 4,  we know the Toda bracket  $ \n   j_{6},\nu_{10},2\cg_{13}                           \nn$ is well-defined, we take   $\overline{2\cg_{14}}$\, from it.
   By Proposition\,\ref{lm}(2), we have $$\s^{\wq}\,\overline{2\cg_{14}}\in \pm \x   j,\nu,2\cg                         \xx\ni \pm2\hat{a}$$
   $\text{mod} \; j\hc \pi_{11}^{S}(S^{0})+\pi_{8}^{S}(\y)\hc2 \cg\bz 16\pi_{15}^{S}(\y)=\x         16\hat{a} \xx$.\;Thus\;  $$ \s^{\wq}\,\overline{2\cg_{14}}\ty \pm2\hat{a}\m16\hat{a}.$$ Since $ \ord(\hat{a})=128$, then $\s^{\wq}\,\overline{2\cg_{14}}$ is of order 64, hence the result follows.
\item[\rm(3)] By Proposition\,\ref{hc1}(4)  and $\cg_{9}\nu_{16}\in\pi_{19}^{9}$  is of order 8,  we know the Toda bracket  $ \n   j_{5},\nu_{9},4\cg_{12}                           \nn$ is well-defined, we take $\overline{4\cg_{13}}$\, from it.
   By Lemma\,\ref{lm}(2), we have $$\s^{\wq}\,\overline{4\cg_{13}}\in \pm \x   j,\nu,4\cg                         \xx\ni \pm4\hat{a}$$
       $\text{mod} \;j\hc \pi_{11}^{S}(S^{0})+\pi_{8}^{S}(\y)\hc4 \cg\bz 16\pi_{15}^{S}(\y)=\x         16\hat{a} \xx.$ Thus\;  $$ \s^{\wq}\,\overline{4\cg_{13}}\ty \pm4\hat{a}\m16\hat{a}.$$Since $ \ord(\hat{a})=128$,  then,\, $\s^{\wq}\,\overline{4\cg_{13}}$ is of order 32, hence the result holds.

  \item[\rm(4)] Since $\nu_{5}\cg_{8}  \in\pi_{15}^{5}$  is of order 8, we know the Toda bracket  $ \n   j_{5},\nu_{5},8\cg_{8}                           \nn$  is well-defined, we take \,$\overline{8\cg_{9}}$ \,from it.
   By Proposition\,\ref{lm}(2),  we have$$\s^{\wq}\,\overline{8\cg_{9}}\in \pm \x   j,\nu,8\cg                         \xx\ni \pm8\hat{a}$$
   $\text{mod} \; j\hc \pi_{11}^{S}(S^{0})+\pi_{8}^{S}(\y)\hc8\cg\bz \x         16\hat{a} \xx$.\; Thus\;  $$ \s^{\wq}\,\overline{8\cg_{9}}\ty \pm8\hat{a}\m16\hat{a}.$$ Since $ \ord(\hat{a})=128$, then  $\s^{\wq}\,\overline{8\cg_{9}}\;$ is of order 16, hence the result follows.

\end{itemize}

\end{proof}

\begin{lem}\label{p19f42}
\begin{itemize}

\item[\rm(1)]  \;$ \pi_{20}(F_{5})=\z/8\n   j\ck_{9}  \nn.$

\item[\rm(2)]  \;$ \pi_{19}(F_{4})=\z/8\n   j\ck_{8}  \nn\jia\z/2\n   j\ch_{8}\nu_{16}   \nn\jia\z_{(2)}\n           j\hc  \overline{8\iota_{19} }\nn.$
\item[\rm(3)] \;$ \pi_{17}(F_{2})=\z/8\n   j\ck_{6}  \nn\jia\z/2\n   j\ch_{6}\nu_{14}   \nn\jia\z/2\n           j\hc  \overline{ \ca_{15}^{2} }\nn.$

\end{itemize}
\end{lem}
\begin{proof}
\begin{itemize}

\item[\rm(1)] 

Since
\begin{eqnarray}
& 
&\notag\pi_{20}(\sk_{32}(F_{5}))  \\\notag& 
=&  \pi_{20}(S^{9}\cup_{\ch_{9}\nu_{17}}e^{21})   \\\notag&
=& j\hc\pi_{21}^{9} \\\notag& 
\tg&  \frac{\pi_{21}^{9}}{\x  \ch_{9}\nu_{17}\xx} \\\notag&
=& \frac{ \z/8\n\ck_{9}\nn\jia\z/2\n\ch_{9}\nu_{17}\nn}{\x  \ch_{9}\nu_{17}\xx}  \\\notag&
\tg& \z/8,  
\end{eqnarray}

hence the result follows.

\item[\rm(2)] We examine the exact sequence $\mathrm{PISK}(19,4)$ given by sequence (\ref{pisk1}).\\
By Corollary\,\ref{bqjg} we have $\sk_{25}(J(M_{S^{8}},S^{18}))=S^{8}$, we observe  $\mathrm{D}(20,4)$ given by diagram (\ref{dkm}),  we have $$
\cok (d_{*}:   \pi_{20}^{19}\rightarrow \pi_{19}(  J(M_{S^{8}},S^{18})                        )=\z/8\n   j\ck_{8}  \nn\jia\z/2\n   j\ch_{8}\nu_{16}   \nn,$$  here, we notice $$f_{4_{*}}(\ca_{18})=\nu_{8}\cg_{11}\ca_{18}-2t'\cg_{8}\nu_{15}\ca_{18}=(\nu_{8}\ca_{11})\cg_{12}-0=0.$$
By $\mathrm{D}(19,4)$ given by diagram (\ref{dkm}), $$\Ker(d_{*}:   \pi_{19}^{19}\rightarrow \pi_{18}(  J(M_{S^{8}},S^{18})                        )=\z_{(2)}\n  8\iota_{19}\nn.$$Hence,   $\mathrm{PISK}(19,4)$ induces the short exact sequence,\\
\xymatrix@C=0.25cm{
&0 \ar[r] & \z/8\n   j\ck_{8}  \nn\jia\z/2\n   j\ch_{8}\nu_{16}   \nn \ar[rr]^{\qquad\,\bz} && \pi_{29}(\sk _{19}(F_{4}))  \ar[rr]^{} && \z_{(2)}\n  8\iota_{19}\nn \ar[r] & 0 }

\noindent It  obviously splits, hence the result holds.
\item[\rm(3)]We recall
$\mathrm{PISK}(m,k)$ given by sequence (\ref{pisk1}) and
  $\mathrm{D}(m+1,k)$ given by diagram (\ref{dkm}). We consider    $\mathrm{PISK}(17,2)$,
by Corollary\,\ref{bqjg}, we know $$\sk_{19}(J(M_{S^{6}},S^{14}))=S^{6},$$ by $\mathrm{D}(18,2)$ we know
$$\cok (d_{*}:   \pi_{18}^{15}\rightarrow \pi_{17}(  J(M_{S^{6}},S^{14})                        )=\z/8\n   j\ck_{6}  \nn\jia\z/2\n   j\ch_{6}\nu_{14}   \nn.$$
 By\,$\mathrm{D}(17,2)$ we have, $$\Ker(d_{*}:   \pi_{17}^{15}\rightarrow \pi_{16}(  J(M_{S^{6}},S^{14})                        )=\z/2\n  \ca_{15}^{2}\nn.$$
Hence   $\mathrm{PISK}(17,2)$ induces the short exact sequence,\\
\xymatrix@C=0.25cm{
&  0 \ar[r] &\z/8\n   j\ck_{6}  \nn\jia\z/2\n   j\ch_{6}\nu_{14}   \nn\ar[rr]^{\qquad\bz} && \pi_{17}(\sk _{23}(F_{2}))  \ar[rr]^{} && \z/2\n  \ca_{15}^{2}\nn\ar[r] & 0. }

\noindent  Similar to  the proof of  Lemma\,\ref{p16f22}(1),  this sequence splits, hence the result holds.

\end{itemize}

\end{proof}

 We recall that  $\pa_{k*_{r}}$ is the homomorphism $\pa_{k*_{r}}: \pi_{r}( S^{8+k})\rightarrow \pi_{r-1}(F_{k})$,\;
 and recall $\mathrm{BUND}(m+1,k)$  given by diagram (\ref{bund1}).
\begin{lem}\label{ck15}
\begin{itemize}

\item[\rm(1)]  \;
$\cok(\pa_{8*_{24}})\tg\z/8\jia\z_{(2)}$,\quad$\Ker( \pa_{8*_{23}})\tg\z/16,$\\
$\cok(\pa_{6*_{22}})\tg\z/8$,\quad$\Ker( \pa_{6*_{21}})\tg\z/8,$\\
 $\cok(\pa_{5*_{21}})\tg\z/8$,\quad$\Ker( \pa_{5*_{20}})\tg\z/4,$\\
 $\cok(\pa_{3*_{19}})\tg\z/8\jia(\z/2)^{2}$,\quad$\Ker( \pa_{3*_{18}})\tg\z/2.$
\item[\rm(2)] \;$\cok(\pa_{4*_{20}})= \pi_{19}(F_{4})=\z/8\n   j\ck_{8}  \nn\jia\z/2\n   j\ch_{8}\cg_{16}   \nn\jia\z_{(2)}\n           j\hc  \overline{8\iota_{19} }\nn$,\\\indent$\Ker( \pa_{4*_{19}})\tg\z/4.$
\item[\rm(3)] \;$\cok(\pa_{2*_{18}})\tg\z/8\jia(\z/2)^{2}$,\quad$\Ker( \pa_{2*_{17}})\tg\z/2.$
\item[\rm(4)] \;$\cok(\pa_{1*_{17}})\tg(\z/8)^{2}$,\quad$\Ker( \pa_{1*_{16}})\tg\z/2.$

\end{itemize}
\end{lem}

\begin{proof}

\begin{itemize}

\item[\rm(1)]  We consider the corresponding diagrams $\mathrm{BUND}(m+1,k)$  and  Proposition\,\ref{hc1}(4),\,(6), we get the result.
\item[\rm(2)] We consider \;$\mathrm{BUND}(20,4)$ and $\mathrm{BUND}(19,4)$,  by Lemma\,\ref{p19f42}(2), Proposition\,\ref{hc1}(4) and Lemma\,\ref{p16f22}(2),\, we get the result.
\item[\rm(3)]We consider \;$\mathrm{BUND}(18,2)$ and $\mathrm{BUND}(17,2)$,  by Lemma\,\ref{p19f42}(3), Lemma\,\ref{p16f22}(1)  and Proposition\,\ref{hc1}(6),(10),\, we get the result.
\item[\rm(4)]  By Proposition\,\ref{hc1}(6),(10) and the group structures of $\pi_{16}^{5},\pi_{17}^{6}$, we obtain $$\nu_{5}\cg_{8}\ca_{15}\ty0\m\nu_{5}\ch_{8},\nu_{5}\lt_{8},$$
(consider its $\s$-image).
Combining with Proposition\,\ref{hc1}(5), and  by use of \;$\mathrm{BUND}(17,1)$ and $\mathrm{BUND}(16,1)$, we get the result.

\end{itemize}

\end{proof}

We recall$\sho(m,k)$ given by sequence (\ref{shortmk}) and$\lon (m,k)$ given by sequence (\ref{longmk}).
\begin{thm}  After localization at 2,
\[\pi _{15+k}(\Sigma ^{k}\mathbb{H}P^{2})\tg\begin{cases}
\mathbb{Z}/128,\;&\text{$k\geq7\; \mbox{and}\; k\neq8$}\\
\mathbb{Z}/128\oplus\mathbb{Z}_{(2)},\;&\text{$k=8$}\\
\z/64,\;&\text{$k=6$}\\
\z/32,\;&\text{$k=5$}\\
\z_{(2)}\jia\z/16\jia\z/2,\;&\text{$k=4$}\\
\mathbb{Z}/16\oplus(\mathbb{Z}/2)^{2},\;&\text{$k=3\;\mbox{or}\;2$}\\
\z/16\jia\z/8,\;&\text{$k=1$}\\
(\mathbb{Z}/2)^{2}\oplus\mathbb{Z}/4,\;&\text{$k=0$}\\
\end{cases}\]\\

\end{thm}
\begin{proof}  For $k\geq9$,  $\pi _{15+k}(\Sigma ^{k}\mathbb{H}P^{2})$ are in the stable range. Hence,
 $$\pi _{15+k}(\Sigma ^{k}\mathbb{H}P^{2})\tg\pi _{15}^{S}(\mathbb{H}P^{2}),\;(k\geq9),$$ by Proposition\,\ref{lm}(2) we get the result.\\
 \indent   For $k=7$,  we consider$ \lon(22,7)$, by Lemma\,\ref{zc15ys}(1),   we know  $\pi_{22}(\s^{7}\y)$ contains the element $\overline{\cg_{15}}$ of order at least 128,
  Since $\pi_{22}(F_{7})\tg\pi_{22}^{11}\tg\z/8$,\; $\pi_{22}^{15}\tg\z/16$, so, in$ \lon(22,7)$,  $i_{7*}$ is monomorphic,  $p_{7*}$  is epimorphic and   $\pi_{22}(\s^{7}\y)=\z/128\n     \overline{\cg_{15}}   \nn$.\\
 \indent   For $k=8$,  we consider$ \sho(23,8)$,  by $\pi_{22}(\s^{7}\y)=\z/128\n     \overline{\cg_{15}}   \nn$ and Lemma\,\ref{zc15ys}(1),   we know  $\pi_{23}(\s^{8}\y)$
 contains the element $\s\,\overline{\cg_{15}}$ of order  128.\;By Lemma\,\ref{ck15}(1),
$$\cok(\pa_{8_{*24}})\tg\z/8\jia\z_{(2)},\;\quad\Ker(\pa_{8_{*23}})\tg\z/16.$$  Since $\pi_{22}(\s^{7}\y)$  contains an element  of order  128,  then\;$$\pi_{23}(\s^{8}\y)\tg\z/128\jia\z_{(2)}.$$
\indent   For $k=6$,  we consider$ \sho(21,6)$, by Lemma\,\ref{zc15ys}(2),    $\pi_{21}(\s^{6}\y)$ contains the element  $\overline{2\cg_{14}}$ of order at least 64,
by Lemma\,\ref{ck15}(1),  we have $\pi_{21}(\s^{6}\y)\tg\z/64$.\\
\indent   For $k=5$,  we consider$ \sho(20,5)$, by Lemma\,\ref{zc15ys}(3),    $\pi_{20}(\s^{5}\y)$ contains the element $\overline{4\cg_{13}}$  of order at least 32, by Lemma\,\ref{ck15}(1),  we have  $\pi_{20}(\s^{5}\y)\tg\z/32$.\\
\indent   For $k=1$,  we consider$ \sho(16,1)$, by Lemma\,\ref{zc15ys}(4),    $\pi_{16}(\s\y)$  contains the element $\overline{8\cg_{9}}$  of order at least 16,  by Lemma\,\ref{ck15}(4),  $\pi_{16}(\s\y)\tg\z/16\jia\z/8$.   Successively, by Lemma\,\ref{zc15ys}(4) we know $$\s^{n}\overline{8\cg_{9}},\;(n\geq0)$$ are all of order 16.\\
\indent   For $k=3$ or 2, we consider$ \sho(15+k,k)$, by the above proof for $k=1$,   $\pi_{18}(\s^{3}\y)$  contains the element   $\s^{2}\,\overline{8\cg_{9}}$ of order 16, $\pi_{17}(\s^{2}\y)$  contains the element  $\s\,\overline{8\cg_{9}}$ of order 16,\;by Lemma\,\ref{ck15}(1),(3) respectively, we get the results.\\
\indent For $k=4$, we consider$ \sho(19,4)$, by Lemma\,\ref{ck15}(2)  and Lemma\,\ref{ext}(2),  we have $\pi_{19}(\s^{4}\mathbb{H}P^{2}) \tg \z_{(2)}\jia K$, where\\  
$K\in\n\z/8\jia\z/2, \;\z/8\jia\z/4,\;\z/16\jia\z/2,\;\z/8\jia(\z/2)^{2},\;\z/32\jia\z/2,\;(\z/8)^{2},\\\;\z/8\jia\z/2\jia\z/4\nn.$\\
By the proof for $k=1$, we know $\s^{3}\,\overline{8\cg_{9}}\in\pi_{19}(\s^{4}\y)$  is of order 16, then$$\pi_{19}(\s^{4}\mathbb{H}P^{2})\tg\z_{(2)}\jia \z/32\jia\z/2\;\, \text{or}\; \z_{(2)}\jia \z/16\jia\z/2.$$ 
 To exclude to the extra solution, we consider the following fibre sequence, where $ \mathcal{U}_{4}=J(M_{S^{11}}, \s^{3}\h P^{2}) $ denotes the fibre of  $  j_{4}: S^{8}\hookrightarrow\s^{4}\h P^{2}$, \\ \xymatrix@C=1cm{
 &&&  \mathcal{U}_{4} \ar[r] & S^{8} \ar[r]^{j_{4}\quad} &\s^{4}\h P^{2}.}
\noindent \\\noindent By checking homology, and by  $\s\mathcal{U}_{4}$ is splitting and by $ \s :\pi_{17}(S^{11})\rightarrow\pi_{18}(S^{12})$ is an isomorphism, we get $\sk_{21}(\mathcal{U}_{4})=S^{11}\vee S^{18}$. So there's a fibre sequence up to dimension 20 \;(here, $ \mathbf{a}=\;\,\nu_{8}\vee f$ for some homotopy class $f$, since $\nu_{8}$ can  extend to \;$\mathcal{U}_{4}$),  \\
   \xymatrix@C=1.7cm{
 &  S^{11}\vee S^{18} \ar[r]^{\quad \mathbf{a}=\;\,\nu_{8}\vee f} & S^{8} \ar[r]^{j_{4}\quad} &\s^{4}\h P^{2}, }

\noindent\\ it  induces the exact sequence,\\
   \xymatrix@C=0.85cm{
&\pi_{19} (S^{8} )\ar[r]^{j_{4_{*}}\quad} &\pi_{19} (\s^{4}\h P^{2})\ar[r]^{} &  \pi_{18}(S^{11})\jia \pi_{18}(S^{18})\ar[r]^{\qquad\quad\mathbf{a}_{*}} _{\qquad\quad (\nu_{8_{*}}, f_{*} )}& \pi_{18}(S^{8}).}

\noindent\\\noindent We denote $(x\cg_{11},y\iota_{18})\in \pi_{18}(S^{11})\jia \pi_{18}(S^{18})$  by $(x,y)$ for simplicity, ($x,y\in\z$).\\
By the result of  $\cok(\partial _{4_{*20}}  )$ shown in Lemma\;\ref{ck15}(2), we know  this $j_{4_{*}}$\,is monic, so $$j_{4_{*}} (\pi_{19} (S^{8} ) )\tg\pi_{19} (S^{8} ) \tg\z/8\jia\z/2.$$ We know $\nu_{8_{*}}: \z/16\n \sigma_{11}\nn\rightarrow \z/8\n\nu_{8}\sigma_{11} \nn\jia else$ has image $ \z/8\n\nu_{8}\sigma_{11} \nn$,
that is, $\mathbf{a}_{*}(1,0)=\nu_{8}\sigma_{11}$. Suppose $\mathbf{a}_{*}(0,1)\ty m\nu_{8}\sigma_{11}\m \cg_{8}\nu_{15},\ca_{8}\mu_{9}$, ($0\leq m\leq7 $ is an integer), then, $$\mathbf{a}_{*}(-m,1)= \ell\cg_{8}\nu_{15}+\ell'\ca_{8}\mu_{9}\text{\;for some\;} \ell,\ell'\in\z.$$ We know  $$\pi_{18}(S^{11})\jia \pi_{18}(S^{18})=\x (1,0)\xx\jia\x (-m,1)\xx\tg \z/16\jia\z_{(2)},$$ 
recall(\cite{Toda}), $$\pi_{18} (S^{8})=\z/8\n \nu_{8}\sigma_{11} \nn\jia\z/8\n \cg_{8}\nu_{15}\nn\jia\z/2\n\ca_{8}\mu_{9}\nn,$$
 therefore,
$$\Ker(\mathbf{a}_{*})=\x(8,0)\xx\jia\x b(-m,1) \xx\tg\z/2\jia\z_{(2)},$$ where $b=\ord(\ell\cg_{8}\nu_{15}+\ell'\ca_{8}\mu_{9})\in\z_{+}$. So we have the exact sequence induced by  above,\\\\\indent
\xymatrix@C=0.4cm{
  &&0 \ar[r] & \z/8\jia\z/2 \ar[r] & \pi_{19} (\s^{4}\h P^{2}) \ar[r] &\z/2\jia\z_{(2)} \ar[r] & 0. }

 \noindent \\\noindent This tells us   $\pi_{19} (\s^{4}\h P^{2})$ cannot have elements of order 32. We have already shown  $$\pi_{19}(\s^{4}\mathbb{H}P^{2})\tg\z_{(2)}\jia \z/32\jia\z/2 \;\,\text{ or\; }\z_{(2)}\jia \z/16\jia\z/2,$$hence, $\pi_{19}(\s^{4}\mathbb{H}P^{2})\tg\z_{(2)}\jia \z/16\jia\z/2.$

\end{proof}

\subsection{ Determinations of some unstable  $\pi_{36+k}(\s^{k}\h P^{2})$ }

$\pi_{36+k}(\s^{k}\h P^{2})$ is in the stable range   $\iff k\geq 36-8+2=30$.\\We will determine $\pi_{36+29}(\s^{29}\h P^{2})$, $\pi_{36+28}(\s^{28}\h P^{2})$ and $\pi_{36+27}(\s^{27}\h P^{2})$, and study $\pi_{36+11}(\s^{11}\h P^{2})$.\\

\indent In this subsection, we will use the following freely,
$$\pi_{64}^{32}\tg(\z/2)^{6},\;\quad\pi_{63}^{31}\tg\pi_{65}^{33}\tg(\z/2)^{5},\;\quad \pi_{66}^{34}\tg(\z/2)^{4},$$
see\, \cite[Theorem 1.1]{32STEM},
and \;$$\pi_{28}^{S}(S^{0})=\z/2\n \lt\ct\nn,\;\,\pi_{29}^{S}(S^{0})=0,$$ \,see\,  \cite[Theorem 2(b), p.\,81]{Oda} and  \cite[Theorem 3(a), p.\,105]{Oda}.
\begin{lem}\label{nlkg} $\n    \nu_{31+n},\lt_{34+n}\ct_{42+n},2\iota_{62+n}          \nn=0,$  for any $ n\geq0$.

\end{lem}

\begin{proof} Since $\pi_{12}^{S}(S^{0})=0$,\;then,\; for any $ n\geq0,$ we have
\begin{eqnarray}& 
&\notag\n        \nu_{31+n},\lt_{34+n}\ct_{42+n},2\iota_{62+n}            \nn  \\\notag& 
\dyd& \n    \nu_{31+n},\lt_{34+n},2\ct_{42+n}          \nn \\\notag&
\dyd& \n    \nu_{31+n},\lt_{34+n},2\iota_{42+n}          \nn\hc\ct_{43+n}   \\\notag&
\xyd& \pi_{43+n}^{31+n}\hc\ct_{43+n}   \\\notag&
=&  0\hc\ct_{43+n}   \\\notag&
=& 0
\end{eqnarray}

$\ind\n        \nu_{31+n},\lt_{34+n}\ct_{42+n},2\iota_{62+n}            \nn=\nu_{31+n}\hc\pi^{34+n}_{63+n}+2\pi_{63+n}^{31+n}=\nu_{31+n}\hc0+0=0$. 

Hence our lemma is established.

\end{proof}

\begin{thm}

 After localization at 2,
\[\pi _{36+k}(\Sigma ^{k}\mathbb{H}P^{2})\tg\begin{cases}
(\z/2)^{6},\;&\text{$k=29$}\\
(\z/2)^{7},\;&\text{$k=28$}\\
(\z/2)^{6},\;&\text{$k=27$}\\
\end{cases}\]\\

\end{thm}
\begin{proof} We only show the case  $k=27$,   similar  proofs work for the other two cases.\\
 \indent  For $\pi _{63}(\Sigma ^{27}\mathbb{H}P^{2})$,  we use $\lon(63,27)$ given by sequence (\ref{longmk}). By Lemma\,\ref{sk0}(1),\;$$\sk_{64}(F_{11})=S^{31}.$$ Since $$\pi_{64}^{35}=\pi_{29}^{S}(S^{0})=0,$$so,\; in$\lon(63,27)$, $i_{27*}$ is monomorphic. Since $$\pi_{28}^{S}(S^{0})=\z/2\n \lt\ct\nn
,\;\nu\lt=0,$$  we observe $\mathrm{BUND}(63,27)$ given by diagram (\ref{bund1}), we have $\pa_{27_{*63}}=0$,\; so in$\lon(63,27)$, $p_{27*}$ is epimorphic. Then$\lon(63,27)$ gives a short exact sequence,\;\\\\\xymatrix@C=0.3cm{
 &&& 0 \ar[r] &(\z/2)^{5}  \ar[rr]^{  } && \pi _{63}(\Sigma ^{27}\mathbb{H}P^{2}) \ar[rr]^{p_{_{_{27*}}}} && \z/2\n \lt_{35}\ct_{43}\nn \ar[r] & 0 }

\noindent\\ Since $$ p_{27*}\n  j_{27},    \nu_{31},\lt_{34}\ct_{42}   \nn= p_{27}\hc\n  j_{27},    \nu_{31},\lt_{34}\ct_{42}  \nn=\n p_{27},  j_{27},    \nu_{31}\nn\hc \lt_{35}\ct_{43}\ni \lt_{35}\ct_{43},$$then, $$\cz\,\overline{\lt_{35}\ct_{43}}\in\n  j_{27},    \nu_{31},\lt_{34}\ct_{42}  \nn,$$such that $$ p_{27*}(  \,\overline{\lt_{35}\ct_{43}}   )= \lt_{35}\ct_{43}.$$
By Lemma\,\ref{nlkg},\,we have\,$$2\,\overline{\lt_{35}\ct_{43}}\in\n  j_{27},    \nu_{31},\lt_{34}\ct_{42}  \nn\hc2 \iota_{63}=-( j_{27}\hc\n       \nu_{31},\lt_{34}\ct_{42} ,2 \iota_{62}\nn)=0,$$
thus \;$2\,\overline{\lt_{35}\ct_{43}}=0$,\; successively\;  $\overline{\lt_{35}\ct_{43}}$ is of order 2. Therefore,\, $$ \pi _{63}(\Sigma ^{27}\mathbb{H}P^{2})\tg(\z/2)^{6}.$$

\end{proof}

\begin{rem} In similar fashion,  $ \pi _{36}^{S}(\mathbb{H}P^{2})$ can also be got, while, this paper pays more attention to  unstable cases.

\end{rem}

For the following lemmas, we point out that, $N.\,Oda$  denotes  the generator\;$\bar{\beta}\in\pi_{20+19}^{19}$ defined in   \cite{20STEM} by $C_{2} $,\;that is,  $C_{2}=\bar{\beta}$,\; see\, \cite[(2.1),\;p.\,52]{Oda}.\\
 \indent  We recall  \cite[Theorem 4.3 (a),\; p.\,104]{Oda},$$\pi_{48}^{19}=\z/2\n   C_{2}\mu_{39}\nn\jia\z/2\n   \cg_{19}^{*}\cg_{41}\nn\jia\z/2\n  \s^{3}P_{1}\nn.\;$$
 To get \;$\nu_{16}\hc\pi_{48}^{19}\; $  and $\cok ( \pa_{11_{*48}})$ we care about,\;
 we give the following lemmas.

\begin{lem}\label{clast}

\begin{itemize}

\item[\rm(1)]
$\nu_{16}C_{2}\in\z/2\n   P(\lt_{33}) \nn\jia\z/2\n  P(\ch_{33}) \nn$,\\
\[\nu_{16}C_{2}\mu_{39}=\begin{cases}
0,\;&\text{ if \;$\nu_{16}C_{2}\ty0\m P(\ch_{33})$}\\
\s\cs^{*},\;&\text{if  $\nu_{16}C_{2}\ty P(\lt_{33})\m P(\ch_{33})\; $}\\
\end{cases}\]
\item[\rm(2)] $         \nu_{16}\cg_{19}^{*}\cg_{41}=\xi_{16}\cg_{24}^{2}\in\pi_{48}^{16}.$
\item[\rm(3)] 
$\nu_{13}\hc P_{1}\ty \lambda\cn_{31}\m \s\pi_{44}^{12}$,\\
$\nu_{16}\hc \s^{3}P_{1}\ty (\s^{3}\lambda)\cn_{34}\m \s^{4}\pi_{44}^{12}.$
\item[\rm(4)]  $[\iota_{15},\nu_{15}]=0$,\;  $\sk_{50}(F_{11})=S^{15}\vee S^{33}$.
    \item[\rm(5)] $\pi_{47}(F_{11}) = j(\hat{i}_{15}\hc \pi_{47}^{15}\jia \hat{i}_{33}\hc \pi_{47}^{33}\jia [\hat{i}_{15},\hat{i}_{33}]\hc\pi_{47}^{47})\tg (\z/2)^{8}\jia(\z/2)^{2}\jia\z_{(2)}$, \,generated by\\\\ $ (j\hc \hat{i}_{15})_{*}\n   \overline{  \nu\ct }_{15},\;   \bar{\ccs}_{15},\;\;\ca^{*\prime}\ca_{31}^{*},\;\;                        \mu^{*\prime\prime }     ,\;\;\s\cs^{*},\;  (\s^{2} \lambda)\cn_{33}, \;\xi_{15}\cg_{33}^{2},\;\cg_{15}\mu_{3,22}\nn,\\\;(j\hc \hat{i}_{33})_{*}\n\cn_{33},\;\cg_{33}^{2}  \nn \;\text{and}\; \;(j\hc[\hat{i}_{15},\hat{i}_{33}])_{*} \iota_{47}.$
    $$\pi_{48}(\s F_{11})= (\s  (j\hc \hat{i}_{15}))\hc \pi_{48}^{16}\jia (\s  (j\hc \hat{i}_{33}))\hc \pi_{48}^{34}\tg \pi_{48}^{16}\jia  \pi_{48}^{34},$$
Here, $\,\hat{i}_{m}:\,S^{m}\hookrightarrow S^{15}\vee S^{33}$,\,(m=15,33) are the inclusions.
\item[\rm(6)]
 If  $\nu_{16}C_{2}\ty0\m P(\ch_{33})$, then, 
 $$\cok ( \pa_{11_{*}}:  \pi_{48}^{19}\rightarrow \pi_{47}(F_{11}))\tg(\z/2)^{6}\jia(\z/2)^{2}\jia\z_{(2)},$$ generated by  $$ (j\hc \hat{i}_{15})_{*}\n  \overline{  \nu\ct }_{15},\;   \bar{\ccs}_{15},\;\ca^{*\prime}\ca_{31}^{*},\;\;                        \mu^{*\prime\prime }     ,\;\s\cs^{*},\;  \cg_{15}\mu_{3,22}\nn,$$ $$ (j\hc \hat{i}_{33})_{*}\n\cn_{33},\;\cg_{33}^{2}\nn\;\; \text{and}\;\;\; (j\hc[\hat{i}_{15},\hat{i}_{33}])_{*} \iota_{47}.$$
 If  $\nu_{16}C_{2}\ty P(\lt_{33})\m P(\ch_{33}) $, then,
$$\cok ( \pa_{11_{*}}:  \pi_{48}^{19}\rightarrow \pi_{47}(F_{11}))\tg(\z/2)^{5}\jia(\z/2)^{2}\jia\z_{(2)},$$ generated by \;$$ (j\hc \hat{i}_{15})_{*}\n  \;  \overline{  \nu\ct }_{15},\;  \bar{\ccs}_{15},\;\ca^{*\prime}\ca_{31}^{*},\;                        \mu^{*\prime\prime }     ,\;  \cg_{15}\mu_{3,22}\nn$$ $$ (j\hc \hat{i}_{33})_{*}\n\cn_{33},\;\cg_{33}^{2}\nn\;\;\text{and}\;\;\; (j\hc[\hat{i}_{15},\hat{i}_{33}])_{*} \iota_{47}.$$
    
\end{itemize}
\end{lem}

\begin{proof}

\begin{itemize}
\item[\rm(1)]  By \cite[Lemma\, 16.4, p.\,53]{20STEM}, we know $\s C_{2}=P(\ca_{41}) $. Then, $$\nu_{17}\s C_{2}= \nu_{17}\hc P(\ca_{41})
=\nu_{17}\hc [\iota_{20},\iota_{20}]\ca_{39}=[\iota_{17},\iota_{17}]\nu_{33}^{2}\ca_{39}=0.$$ Then, \;
 $$\nu_{16}C_{2}\in \Ker(\s:    \pi_{39}^{16}\rightarrow \pi_{40}^{17})= P(\pi_{41}^{33}).$$   By \cite[(3.4),\,p.\,12 ]{2324STEM} we know $P(\pi_{41}^{33})=\z/2\n  \; P(\lt_{33}) \nn\jia\z/2\n \; P(\ch_{33}) \nn$.  By \cite[Theorem 14.1]{Toda}, we have $\lt\mu=\ca^{2}\rho$ and $\ch\mu=0$.
By \cite[(5.8)]{32STEM},  $P(\ca_{33}^{2}\rho_{35})=\s\cs^{*}$. Hence,   $P(\lt_{33})\hc\mu_{39}=\s\cs^{*}$ and  $P(\ch_{33})\hc\mu_{39}=0$. Hence the result.

\item[\rm(2)]   By  \cite[ Proposition\, 3.5 (3), p. 60 ]{Oda} we know the two relations,\;\\the first one,  $\nu_{13}\cg_{16}^{*}\ty t\xi_{13}\cg_{31}\m \lambda \cg_{31},\;(t:\; \mbox{odd})$, \\ the second one, $(\s^{3} \lambda) \cg_{34}=0$.\; \\Then,\; $\nu_{16}\cg_{19}^{*}\ty t\xi_{16}\cg_{34}\m 0,\;(t:\; \mbox{odd})$.\;  Hence, $$\nu_{16}\cg_{19}^{*}\cg_{41}=t\xi_{16}\cg_{34}^{2}=\xi_{16}\cg_{34}^{2} ,\; (\ord(\cg_{34}^{2})=2).$$

\item[\rm(3)]  We noticed these two facts,

 the first one,  $H(P_{1})=\cn_{31}$,\; \cite[Proposition\, 3.3 (4), p. 90]{Oda};\; \\ the second one,  $H(\lambda\cn_{31})=\nu_{25}^{2}\cn_{31} $,\;
 \cite[the proof of $\pi_{45}^{13}$]{32STEM}, \;\\
   then   $H(\nu_{13}P_{1})=\s (\nu_{12}\wedge \nu_{12})\hc \cn_{31}=\nu_{25}^{2}\cn_{31}$.\;\\
   Hence,   $\nu_{13}P_{1}\ty \lambda\cn_{31}\m \s\pi_{44}^{12}$,\; then,  $\nu_{16}\s^{3}P_{1}\ty (\s^{3} \lambda)\cn_{34}\m \s^{4}\pi_{44}^{12}$.

\item[\rm(4)] By Proposition\,\ref{hc1}(11), (4), we know  $[\iota_{15}, \nu_{15}]=\pm P(\iota_{31})\hc \nu_{29}=2\cg_{15}^{2}\nu_{29}=0$. \;By Lemma\,\ref{sk0}, we have $\sk_{50}(F_{11})=S^{15}\vee S^{33}$.
\item[\rm(5)] By the Hilton-Milnor formula, \;       $$ \lo( S^{15}\vee S^{33})\simeq\lo S^{15}\times\lo S^{33}\times\lo S^{47}\times\lo S^{61}\times\cdots,$$
 then, up to isomorphism,\,$$\pi_{47}(\sk_{50}(F_{11}))=\pi_{46}(\lo S^{15}\times\lo S^{33}\times\lo S^{47})=\pi^{15}_{47}\jia\pi_{47}^{33}\jia\pi_{47}^{47},$$
and the isomorphism $\pi_{47}(\sk_{50}(F_{11}))\tg \pi^{15}_{47}\jia\pi_{47}^{33}\jia\pi_{47}^{47}$ are induced by the inclusions and iterated Whitehead products,\;(see \cite[Theorem 8.1,\,p.\,533]{GW}).\\
 By \;\cite[Theorem 1.1]{32STEM}:\\ \\\centerline{$\pi_{47}^{15}\tg(\z/2)^{8}$,}\\\\ generated by $ \overline{  \nu\ct }_{15},\;   \bar{\ccs}_{15},\;\ca^{*\prime}\ca_{31}^{*},\;                        \mu^{*\prime\prime }     ,\;\s\cs^{*},\;(\s^{2} \lambda)\cn_{33},\; \xi_{15}\cg_{33}^{2},   \;\cg_{15}\mu_{3,22} $.\\
By \cite{Toda},\; $\pi_{47}^{33}=\z/2\n \cn_{33}\nn\jia\z/2\n \cg_{33}^{2}\nn$,\\ hence,\;
$\pi_{47}(\sk_{50}(F_{11})) \tg (\z/2)^{8}\jia(\z/2)^{2}\jia\z_{(2)}$\;\,generated by\\\\$ (j\hc \hat{i}_{15})_{*}\n    \overline{  \nu\ct }_{15},  \; \bar{\ccs}_{15},\;\ca^{*\prime}\ca_{31}^{*},\;                        \mu^{*\prime\prime }     ,\;\s\cs^{*},\; (\s^{2} \lambda)\cn_{33},\; \xi_{15}\cg_{33}^{2},\;\cg_{15}\mu_{3,22}\nn,\\ (j\hc \hat{i}_{33})_{*}\n\cn_{33},\;\cg_{33}^{2}  \nn \;\text{and}\; \;(j\hc[\,\hat{i}_{15},\hat{i}_{33}])_{*} \iota_{47}.$\\ \\
Since $\sk_{51}(\s F_{11})=S^{16}\vee S^{34}$,\;  then, $$\pi_{48}(\s F_{11})= (\s  (j\hc \hat{i}_{15}))\hc \pi_{48}^{16}\jia (\s  (j\hc \hat{i}_{33}))\hc \pi_{48}^{34}\tg \pi_{48}^{16}\jia  \pi_{48}^{34}.$$

\item[\rm(6)] By Lemma\,\ref{xh2}, there exists the commutative diagram,\\

\indent\qquad\qquad\qquad\qquad \xymatrix@C=1.5cm{
  \pi_{48}^{19} \ar[d]_{\nu_{16*}} \ar[r]^{\pa_{11*}\;\;} & \pi_{47}(F_{11}))\ar[d]^{\s} \\
    \pi_{48}^{16}\ar[r]^{(\s j)_{*}\quad} & \pi_{48}(\s F_{11}))   }

    \noindent that is,\\

\indent\qquad\xymatrix@C=2.9cm{
  \pi_{48}^{19} \ar[d]_{\nu_{16*}} \ar[r]^{\pa_{11*}\;\;\qquad\qquad\qquad\qquad} &  j\hc \hat{i}_{15}\hc \pi_{47}^{15}\jia j\hc \hat{i}_{33}\hc \pi_{47}^{33}\jia j\hc [\,\hat{i}_{15},\hat{i}_{33}]\hc\pi_{47}^{47}\ar[d]^{\s} \\
    \pi_{48}^{16}\ar[r]^{(\s j)_{*}\qquad} &\s  (j\hc \hat{i}_{15})\hc \pi_{48}^{16}\jia \s  (j\hc \hat{i}_{33}) \hc \pi_{48}^{34}  }

\noindent By \;\cite[Theorem 1.1]{32STEM},  $\pi_{32+15}^{15}\stackrel{\s}\longrightarrow \pi_{32+16}^{16}$  is monomorphic, we notice $$(\mathrm{Im}\pa_{11*})\cap (j\hc [\,\hat{i}_{15},\hat{i}_{33}]\hc\pi_{47}^{47})=0,$$
then, by (1),(2),(3) and (5)  of this Lemma, we get the result.

\end{itemize}
\end{proof}

 \indent We recall \cite[Theorem 2 (b), p.\,81]{Oda}, $$\pi_{47}^{19}=\z/2\n   C_{2}^{\,(1)}\nn\jia\z/2\n  \s\mathbf{F}_{2}\nn\jia\z/2\n  \lt_{19}\ct_{27}\nn,$$ and \cite{2324STEM},$$\pi_{38}^{15}\tg\z/16\jia\z/8\jia(\z/2)^{4}$$
 generated by\;$\bar{\rho}_{15},\;\nu_{15}\ct_{18},\;\phi_{15},\;\psi_{15},\;\bar{\lt}^{*\prime},\;\ch^{*\prime}$.\;\;($\mathbf{F}_{2}\in\pi_{46}^{18}$ is the generator, not the fibre $F_{2}$).\;\; \\\indent To study $\nu_{16}\hc\pi_{47}^{19}$  and $\Ker ( \pa_{11_{*47}})$,\;
 we give the following lemma, we will freely use the exactness of \textit{EHP} sequence (or say $\s$\textit{HP} sequence ) of $\pi_{*}(S^{n})$. 
 \begin{lem}\label{c21}

\begin{itemize}

\item[\rm(1)]  $  C_{2}^{\,(1)}\in \n  C_{2},\;2\iota_{39},\;8\cg_{39}                                     \nn_{1}.$

\item[\rm(2)]  $D_{2}^{\,(1)}\cg_{39}\ty  \n  \bar{\lt}^{*\prime},\;2\iota_{38},\;8\cg_{38}                                    \nn_{1}\m \s\pi_{46}^{15}$.
\end{itemize}
\end{lem}
\begin{proof}
\begin{itemize}
\item[\rm(1)]   
As a matter of fact, it's just the definition of $C_{2}^{\,(1)}$ given by \textit{N.\,Oda}, see \cite[Definition 3.6]{m11}.
\item[\rm(2)] We recall the following facts,
$$D_{1}^{\,(1)}\ty \mu^{*\prime}\m \s\pi_{38}^{14},\,\text{(\cite[(2.3),\,p.\,52]{Oda})};
H(\mu^{*\prime})=\ca_{29}\mu_{30},(\text{\cite[(3.11)]{2324STEM}),}$$
it's easy to get $\ca\mu\cg=\lt\mu$ by \cite[Theorem 14.1, p.\,190]{Toda}, then,\;$$H(D_{2}^{\,(1)}\cg_{39})=\ca_{29}\mu_{30}\cg_{39}=\lt_{29}\mu_{37}.$$  
Since $H(\bar{\lt}^{*\prime} )=\ca_{29}\lt_{30}=\lt_{29}\ca_{37}$,\;(\cite[(3.4)]{2324STEM}), $\mu\in\x \ca,2\iota,8\cg\xx$,\;(\cite[p.\,189]{Toda})\\ then,
\begin{eqnarray}
& 
&\notag H\n \bar{\lt}^{*\prime},\;2\iota_{38},\;8\cg_{38}                                    \nn_{1}  \\\notag& 
\xyd& \n \lt_{29}\ca_{37},\;2\iota_{38},\;8\cg_{38}                                     \nn_{1}  \\\notag&
\dyd& \lt_{29}\hc \n\ca_{37},\;2\iota_{38},\;8\cg_{38}                                     \nn_{1}   \\\notag&
\ni& \lt_{29}\mu_{37}. 
\end{eqnarray}
$\lt^{2}=\lt\ch=0$,\;(\cite[Theorem 14.1]{Toda}), $8\pi_{29+10}^{29}=0$, thus,
$$\ind \n\lt_{29}\ca_{37},\;2\iota_{38},\;8\cg_{38}\nn_{1}= \lt_{29}\ca_{37}\hc\pi_{46}^{38}=\lt_{29}\ca_{37}\hc \x \lt_{38},\ch_{38}\xx=0.$$
Hence,\;$H\n  \bar{\lt}^{*\prime},\;2\iota_{38},\;8\cg_{38}\nn_{1}=\lt_{29}\mu_{37}$. 
Then the result holds.
\end{itemize}
\end{proof}

\begin{lem}\label{klast}

\begin{itemize}

\item[\rm(1)]  $\nu_{15} \mathbf{F}_{2}=\nu_{16}\lt_{19}\ct_{27}=0.$
\item[\rm(2)]
 \[\nu_{16} C_{2}^{\,(1)}\begin{cases}
=0,\;&\text{ if \;$\nu_{16}C_{2}\ty0\m P(\ch_{33})$}\\
\ty (\s D_{1}^{\,(1)})\hc\cg_{40}\m \w^{*}_{16},\,\cn_{16}^{*},\co_{3,16},\;&\text{if  $\nu_{16}C_{2}\ty P(\lt_{33})\m P(\ch_{33})\; $}\\
\end{cases}\]
\item[\rm(3)] $\Ker(  \pa_{11_{*}}:  \pi_{47}^{19}\rightarrow \pi_{48}(F_{11}))=\z/2\n   \s \mathbf{F}_{2}      \nn\jia\z/2\n      \lt_{19}\ct_{27}\nn\jia G$,\;

 \[G=\begin{cases}
\z/2\n C_{2}^{\,(1)} \nn,\;&\text{ if \;$\nu_{16}C_{2}\ty0\m P(\ch_{33})$}\\
0\;&\text{if  $\nu_{16}C_{2}\ty P(\lt_{33})\m P(\ch_{33})\; $}\\
\end{cases}\]
\end{itemize}
\end{lem}
\begin{proof}

\begin{itemize}

\item[\rm(1)]   $ \nu_{16}\lt_{19}\ct_{27}=0$ is obvious, see Proposition\,\ref{hc1}(6).\;Next, we consider $\nu_{15} \mathbf{F}_{2}$.\\
 By \cite[Proposition\, 3.2(2), p.\,89]{Oda}, we know $H(\mathbf{F}_{2})=a\ck_{35},\;(a : \mbox{odd})$.\;By \cite[Proposition 2.2,\,18]{Toda} and \cite[Proposition 3.1,\,25]{Toda}, we have
$$H(\nu_{15} \mathbf{F}_{2})= \s(\nu_{14}\wedge \nu_{14})  \hc a\ck_{35}=a\nu_{29}^{2}\ck_{35}=a\nu_{29}(\nu_{32}\ck_{35}),$$
By \cite[Theorem 14.1  ii)]{Toda},\;$ \nu\ck=0$,\;so  $H(\nu_{15} \mathbf{F}_{2})=0$. Then\;$\nu_{15} \mathbf{F}_{2}\in \s\pi_{45}^{14}$.\\
    By \cite[Theorem 3 (c),\; p.\,106]{Oda},  we know\; $$\Ker(\s^{\wq}: \pi_{45}^{14}\rightarrow \pi_{31}^{S}(S^{0}) )=\z/4\n 2\w^{*}_{14} \nn,$$
    and  $\w^{*}_{15}$ is of order 2.\\
Since \;$\s^{\wq}(\nu_{15}\mathbf{F}_{2})=\nu\hc\s^{\wq} \mathbf{F}_{2}\in \nu\hc \pi_{28}^{S}(S^{0}) =\x  \nu\lt\ct\xx=0$,\\
then \;$\nu_{15}\mathbf{F}_{2}\in \x2\w^{*}_{15}   \xx=0$.\;
\item[\rm(2)]
By Lemma\,\ref{c21}(1), we know,
$\nu_{16}C_{2}^{\,(1)} \in\nu_{16}\hc\n C_{2},\;2\iota_{39},\;8\cg_{39}\nn_{1}$.   \\
\indent\qquad If  $\nu_{16}C_{2}\ty0\m P\ch_{33}$, \\  we have $\nu_{16}C_{2}= tP\ch_{33}$ for some $t\in\z$; by \cite[Proposition\, 3.2 (1), p. 57]{Oda}, we have $\x   \ch,2\iota,8\cg                 \xx=0$, then,\\\\
$\textcolor{white}{7777777777888777767}\;\;\,\nu_{16}C_{2}^{\,(1)} \\\textcolor{white}{777777777788877777}\in\nu_{16}\hc\n C_{2},\;2\iota_{39},\;8\cg_{39}  \nn_{1}\\\textcolor{white}{777777777788877777}\xyd
  \n    tP(\ch_{33})           ,\;2\iota_{39},\;8\cg_{39}                                  \nn_{1}\\\textcolor{white}{777777777788877777}\dyd P\n            t(\ch_{33})                 ,\;2\iota_{41},\;8\cg_{41} \nn_{3}\\\textcolor{white}{777777777788877777}=0,  $\;\\\\ 
  we know $8\pi_{40}^{16}=0$, so
$$\ind  \n      tP(\ch_{33})           ,\;2\iota_{39},\;8\cg_{39} \;                                    \nn_{1}= tP(\ch_{33})   \hc \x \lt_{39}, \ch_{39}\xx+\pi_{40}^{16}\hc8\cg_{40}=0,$$then,$$\nu_{16}C_{2}^{\,(1)}=0,\;  \text{(if}\; \nu_{16}C_{2}\ty0\m P\ch_{33}).$$
Essentially,  we also shown $ \n  P(\ch_{33}),\;2\iota_{39},\;8\cg_{39}\nn_{1}=0.$\\
\indent\qquad If $\nu_{16}C_{2}\ty P\lt_{33}\m P\ch_{33}\; $,\\
by   Lemma\,\,\ref{c21}\,(2) and  $P(\lt_{33})=\s \bar{\lt}^{*\prime}$,(\cite[(3.4)]{2324STEM}), we have, for some $z\in\z$,  
\begin{eqnarray}
& 
&\notag\nu_{16}C_{2}^{\,(1)}  \\\notag& 
\in& \nu_{16}\hc\n C_{2},\;2\iota_{39},\;8\cg_{39}                                    \nn_{1}  \\\notag& 
\xyd& \n  P(\lt_{33})+zP(\ch_{33}),\;2\iota_{39},\;8\cg_{39}                                     \nn_{1}  \\\notag&
\xyd&  \n  P(\lt_{33}),\;2\iota_{39},\;8\cg_{39}                                    \nn_{1} + \n  zP(\ch_{33}),\;2\iota_{39},\;8\cg_{39}                                     \nn_{1}   \\\notag&
=& \n  P(\lt_{33}),\;2\iota_{39},\;8\cg_{39}                                     \nn_{1} +0  \\\notag&
=& \n  P(\lt_{33}),\;2\iota_{39},\;8\cg_{39}                                     \nn_{1}     \\\notag&
=& \n\s \bar{\lt}^{*\prime},\;2\iota_{39},\;8\cg_{39}                                    \nn_{1}    \\\notag&
\dyd& \s \n  \bar{\lt}^{*\prime},\;2\iota_{38},\;8\cg_{38}                                     \nn,   \\\notag&
\ni& \s ( D_{1}^{\,(1)}\cg_{39}+\s\mathbbm{x})   \\\notag&
 =& (\s D_{1}^{\,(1)})\cg_{40}+\s^{2}\mathbbm{x}\;\;\text{for some}\; \mathbbm{x}\in\pi_{45}^{15}, 
\end{eqnarray}
$\ind \n  P(\lt_{33}), 2\iota_{39}, 8\cg_{39}\nn_{1}=P(\lt_{33})   \hc \x \lt_{39}, \ch_{39}\xx+\pi_{40}^{16}\hc8\cg_{40}=0$; then,\,
$\nu_{16}C_{2}^{\,(1)} \ty   (\s D_{1}^{\,(1)})\cg_{40}\m \s^{2} \pi_{45}^{14}$,  (if $\nu_{16}C_{2}\ty P\lt_{33}\m P\ch_{33}).$

\item[\rm(3)]  We have shown that $\sk_{50}(F_{11})=S^{15}\vee S^{33}$  in Lemma\,\ref{clast}(4). By use of Lemma\;\ref{xh2}, there exists the commutative diagram, \,(here, by abuse of notation,  $j$ denotes the inclusions although they are different inclusions, in this case, these two $j$ and these two  $\s j$ all induce monomorphisms.)\\
\indent\qquad\qquad\qquad \xymatrix@C=2.9cm{
  \pi_{47}^{19} \ar[d]_{\nu_{16*}} \ar[r]^{\pa_{11*}\;\;\qquad} &  j\hc \pi_{46}^{15}\jia j\hc \pi_{46}^{33} \ar[d]^{\s} \\
    \pi_{47}^{16}\ar[r]^{(\s j)_{*}\qquad} &(\s  j)\hc \pi_{47}^{16}\jia (\s j) \hc \pi_{47}^{34}  }

\noindent By  \cite[Theorem\,3 (c),\; p.\,106]{Oda},\; we know   $\pi_{46}^{15}\stackrel{\s}\longrightarrow\pi_{47}^{16}$ is a monomorphism,\;
then \;$ j\hc \pi_{46}^{15}\jia j\hc \pi_{46}^{33}   \stackrel{\s}\longrightarrow         (\s  j)\hc \pi_{47}^{16}\jia (\s j) \hc \pi_{47}^{34} $  is a monomorphism.
So  $ \Ker(\nu_{16*})=\Ker(\pa_{11*})$. Then, by (1),(2) of this lemma, we get the result.
\end{itemize}

\end{proof}

 We recall that $\nu_{16}C_{2}\in\z/2\n   P(\lt_{33}) \nn\jia\z/2\n  P(\ch_{33}) \nn$,  equivalently speaking, either   $\nu_{16}C_{2}\ty0\m P(\ch_{33})$ or  $\nu_{16}C_{2}\ty P(\lt_{33})\m P(\ch_{33})$.\\
 \indent By Lemma\,\ref{clast}\,(6) and Lemma\,\ref{klast}\,(3), we have
\begin{thm}

 After localization at 2, \\

 if \;$\nu_{16}C_{2}\ty0\m P(\ch_{33})$,  then\;$ \pi_{47}(\s^{11}\y)$ is decided by the following short exact sequence, \\
\resizebox{12.3cm}{0.40cm}{
\xymatrix@C=0.3cm{
 0 \ar[r] & (\z/2)^{8}\jia\z_{(2)} \ar[rr]^{\xyd} && \pi_{47}(\s^{11}\y) \ar[rr]^{ p_{_{11*}}\qquad\qquad\quad} && \z/2\n   \s \mathbf{F}_{2}      \nn\jia\z/2\n      \lt_{19}\ct_{27}\nn\jia \z/2\n C_{2}^{\,(1)}\nn\ar[r] & 0; }}{\textcolor{white} {.}}\\\\
 \indent  if  $\nu_{16}C_{2}\ty P(\lt_{33})\m P(\ch_{33})$, then\;$ \pi_{47}(\s^{11}\y)$ is decided by the following short exact sequence, \\
\xymatrix@C=0.3cm{
 0 \ar[r] & (\z/2)^{7}\jia\z_{(2)} \ar[rr]^{\xyd} && \pi_{47}(\s^{11}\y) \ar[rr]^{\qquad p_{_{11*}}\qquad\qquad\quad} && \z/2\n   \s \mathbf{F}_{2}      \nn\jia\z/2\n      \lt_{19}\ct_{27}\nn\ar[r] & 0, }{\textcolor{white} {.}}\\\\
 \noindent where $(\z/2)^{8}\jia\z_{(2)}$ and  $(\z/2)^{7}\jia\z_{(2)}$ is  the abbreviation of  $ i_{11*}( \pi_{47}(F_{11}))=\cok ( \pa_{11_{*_{48}}})$  shown in Lemma\,\,\ref{clast}\,(6). \\\indent\qquad\qquad\qquad\qquad\qquad\qquad\qquad\qquad\qquad\qquad\qquad\qquad\qquad\qquad\qquad\,\quad\boxed{}
\end{thm}

\section{Applications}\subsection{The Hopf fibration $S^{11}\rightarrow\y$  localized at 2}

We notice on page 38 of \cite{James76}, its  $\mathrm{P}_{n,k}:=\h P^{n+1}/ \h P^{n+1-k}$,  see \cite[p.\,21]{James76}. Then we have the following proposition.
\begin{pro}\label{2bgx} \cite[p.\,38]{James76} Let    $$h_{n}: S^{4n+3}\rightarrow\h P^{n}$$ be the homotopy class of the Hopf fibration, $$ p^{\h}_{n}:  \h P^{n}\rightarrow S^{4n}$$ be the pinch map. Then,\;
$p^{\h}_{n}\hc h_{n}=n(\s^{4n-4}\mathbbm{x}) $, \;where  $\mathbbm{x}\in\pi_{7}(S^{4})$  is   the homotopy class of the Hopf fibration.\;
Of course, after localization at 2, we can take $p^{\h}_{n}\hc h_{n}=n\nu_{4n}.$
\end{pro}

\begin{thm}
After localization at 2, for the  homotopy class of  the Hopf fibration  $h: S^{11}\rightarrow\y,$  \;then the folowing relation holds $$\s h=x\mathbbm{a} \, \text{\;for some odd}\;x, $$ \,where  $  \mathbbm{a}\in  \n   j_{1},\nu_{5},2\nu_{8}                           \nn  $ generates  $\pi_{12}(\s\y)\tg\z/8$.

\end{thm}
\begin{proof}
By \cite[Remark 3\,(3)]{Mukai1}, we know, $\ord(\s^{\wq}h)=\frac{6!}{2}=8\times9\times5$,\;(has not been localized).
Then, after localization at 2,  $\ord(\s^{\wq}h)=8$.
By Theorem\,\ref{2c11aa} we know  $\pi_{12}(\s\y)\tg\z/8\n    \mathbbm{a} \nn$,\,$  \mathbbm{a}\in  \n   j_{1},\nu_{5},2\nu_{8}                           \nn $, so   $\s h$ is of order 8 and    $\s h=x\mathbbm{a} \,$ for some odd $x$. Hence the result holds.

\end{proof}

By the above theorem, we have
\begin{cor} After localization at 2, $$\s\h P^{3}\simeq \s\y\cup_{\mathbbm{a}} e^{13},$$\; \; where   $\mathbbm{a}\in  \n   j_{1},\nu_{5},2\nu_{8}                           \nn\xyd  \pi_{12}(\s\y)=\z/8\n  \mathbbm{a} \nn.\bzd$
\end{cor}

\subsection{Two classification theorems  }

 \begin{definition}\;For a \textit{CW} complex $ X$,  if $\widetilde{ H}_{*}(X;\z)\tg \widetilde{ H}_{*}(\h P^{n};\z)$ as graded groups, we call $X$  a homology $n$-dimensional quaternionic projective space.
  \end{definition}

\indent We recall that, after localization at 3,\; $$\pi_{12}(S^{5})=\z/3\n \alpha_{2}(5)\nn,\;\;\pi_{12}(S^{9})=\z/3\n \alpha_{1}(9)\nn,$$ $h: S^{11}\rightarrow \y$ is   the attaching class of  $\h P^{3}=\y \cup_{h} e^{12}$.\;We also recall Theorem\,\ref{3c11th}, after localization at 3,\;$$\pi _{15}(\Sigma ^{4}\mathbb{H}P^{2})  = \mathbb{Z}_{(3)}\n j_{4}[\iota_{8}, \iota_{8}] \nn\oplus\mathbb{Z}/9 \n \s^{4}h  \nn.$$After localization at 3, let\\ \centerline{$A=S^{5}\cup_{\alpha_{2}(5)}\, e^{13}$,}\; \\ let\\
\centerline{$i_{5}': S^{5}\rightarrow S^{5}\vee S^{9}$\;\,and\;\, $ i_{9}': S^{5}\rightarrow S^{5}\vee S^{9}$\quad be the inclusions, }
\\ let \\  \centerline{$c_{k}= \s^{k-1}(  i_{5}'\hc\alpha_{2}(5)    +i_{9}'   \hc \alpha_{1}(9)    )$,}\\\\
by use of these notations, we give the following classification theorem.\\
\begin{thm}\label{fldl} After localization at 3,  up to homotopy,\\  the $k$-fold ($k\geq1$ but $k\neq4$) suspension of the simply-connected homology $3$-dimensional  quaternionic projective spaces can be classified as the following, \\
\indent$ \Sigma ^{k}\mathbb{H}P^{3}, \;\quad \Sigma ^{k}\mathbb{H}P^{2}\cup_{3\Sigma^{k}h}e^{12+k},\;\quad\Sigma ^{k}\mathbb{H}P^{2}\vee S^{12+k},$\\
\indent$\Sigma^{k-1}A\vee S^{8+k},\; S^{4+k}\vee\Sigma^{4+k}\mathbb{H}P^{2},$ \;$(S^{4+k}\vee S^{8+k})\cup _{c_{k}} e^{12+k}\;$  \\and\;$S^{4+k}\vee S^{8+k} \vee S^{12+k}.$

\end{thm}
\begin{proof}
 After localization at 3, for $k\geq1$ but $k\neq4$.\\\indent Firstly, $\pi_{11+k}(\s^{k}\y)=\z/9\n \s^{k}h \nn$,\;(Theorem\,\ref{3c11th}), by Corollary\,\ref{ygd}
and Corollary\,\ref{lgd}, any $k$-fold  suspension of a simply-connected homology $3$-dimensional  quaternionic projective space,  that is, a \textit{CW}  complex of type $$ S^{4+k}\cup e^{8+k}\cup e^{12+k},$$ must be homotopy equivalent a  \textit{CW} complex listed in our theorem.\\\indent
Next,\;we observe the  Steenrod module structures of  the \textit{CW} complexes shown in our theorem.
 For convenience, we use $P^{n}_{*}[\,i\,]$ to denote $$\;\;P^{n}_{*} :\;\widetilde{H}_{i}(-)\rightarrow \widetilde{H}_{i-4n}(-).$$
 For $\Sigma ^{k}\mathbb{H}P^{2}\cup_{3\Sigma^{k}h}e^{12+k}$,\;\,$P^{1}_{*}[12+k]=0$,  $P^{1}_{*}[8+k]\neq0$.
On  $P^{1}_{*}[12+k]=0$,  because  $$( \Sigma ^{k}\mathbb{H}P^{2}\cup_{3\Sigma^{k}h}e^{12+k})/ S^{4+k}\simeq  S^{8+k}\vee S^{12+k}.$$
For $(S^{4+k}\vee S^{8+k})\cup _{c_{k}} e^{12+k}$ ,  \;$P^{1}_{*}[12+k]\neq0$,  $P^{1}_{*}[8+k]=0$,  because 
       $$(S^{4+k}\vee S^{8+k})\cup _{c_{k}} e^{12+k})/S^{4+k}\simeq  \Sigma ^{k+4}\mathbb{H}P^{2}.$$
 For $    \Sigma^{k-1}A\vee S^{8+k}$,\;$P^{1}_{*}[12+k]=P^{1}_{*}[8+k]=0.$\\ \indent Then, we observe the homotopy groups,
 $$\pi_{11+k}(\Sigma^{k-1}A\vee S^{8+k})=\pi_{11+k}(          (    S^{4+k}\cup_{\alpha_{2}(4+k)}\, e^{12+k}   )\vee S^{8+k}) \tg \z/3,$$  
    $$\pi_{11+k}(S^{4+k}\vee S^{8+k} \vee S^{12+k})\tg(\z/3)^{2},$$$$
     \pi_{11+k}(\Sigma ^{k}\mathbb{H}P^{2}\cup_{3\Sigma^{k}h}e^{12+k})\tg \frac{\z/9\n \Sigma^{k}h  \nn}{\x 3\Sigma^{k}h\xx}\tg\z/3,$$
     $$\pi_{11+k}(\Sigma ^{k}\mathbb{H}P^{2}\vee S^{12+k})\tg \z/9.$$
Hence, we only need to show    $$S^{4+k}\vee\Sigma^{4+k}\mathbb{H}P^{2} \;\cancel{\simeq}\;(S^{4+k}\vee S^{8+k})\cup _{c_{k}} e^{12+k}.$$
Since  $$(S^{4+k}\vee\Sigma^{4+k}\mathbb{H}P^{2})/S^{8+k}\simeq  S^{4+k}\vee S^{12+k},$$
      $$((S^{4+k}\vee S^{8+k})\cup _{c_{k}} e^{12+k})/ S^{8+k}\simeq  S^{4+k}\cup_{\alpha _{2}(4+k)}e^{12+k},$$
     $$ \pi_{11+k}(S^{4+k}\vee S^{12+k})\tg\z/3,\; \;  \pi_{11+k}(   S^{4+k}\cup_{\alpha _{2}(4+k)}e^{12+k}    )=0,$$
    then, $S^{4+k}\vee\Sigma^{4+k}\mathbb{H}P^{2}$    cannot be  homotopy equivalent to $(S^{4+k}\vee S^{8+k})\cup _{c_{k}} e^{12+k}$.\\
The above facts imply  the \textit{CW} complexes listed in our theorem    cannot be  homotopy equivalent to each other.

\end{proof}

\begin{thm}\label{wjy} After localization at 3, suppose $Z$ is a $k$-fold ($k\geq1, \;but\,k\neq4$) suspension of a simply connected homology $3$-dimension  quaternionic projective space,  
\begin{itemize}
\item[\rm(1)] if  $ \widetilde{H}_{*}(Z)$ has nontrivial has nontrivial Steenrod operations $P^{1}_{\ast}$ of dimension $8+k$ and $12+k$,\;
then,$$Z\simeq \s^{k}\h P^{3}.$$

\item[\rm(2)]  if  $ \widetilde{H}_{*}(Z)\tg   \widetilde{H}_{*}( \s^{k}\h P^{3})$  as  Steenrod modules,  then
$$Z\simeq \s^{k}\h P^{3}.$$

\end{itemize}
\end{thm}
\begin{proof}
\begin{itemize}
\item[\rm(1)]
After localization at 3, for  $k\geq1,k\neq4$,\; for convenience, we use  $P^{n}_{*}[\,i\,]$ to denote $$\;\;P^{n}_{*} :\;\widetilde{H}_{i}(-)\rightarrow \widetilde{H}_{i-4n}(-).$$ Since $P^{1}_{*}[8+k]\neq0$ and $\pi_{7+k}(S^{4+k})\tg\z/3$,  we can take
$$\sk_{11+k}(Z)=\s^{k}\y,\;Z=\s^{k}\y\cup_{f} e^{12+k},$$ where\;$f\in\pi_{11+k}(\s^{k}\y)=\z/9\n \s^{k}h\nn$. Now, to obtain a contraction, suppose $f$ can be divisible by 3, that is, $\cz\mathbbm{x}\in\pi_{11+k}(\s^{k}\y)$,\;such  that $f=3\mathbbm{x}$.
Then  $$Z/S^{4+k}\simeq S^{8+k}\cup _{q\hc3\mathbbm{x} } e^{12+k}\simeq S^{8+k}\cup _{0 } e^{12+k}\simeq S^{8+k}\vee S^{12+k},$$
here $ q:  \s^{k}\y\rightarrow ( \s^{k}\y)/S^{4+k}$  is the pinch, and we notice that  $\pi_{11+k}(S^{8+k})\tg\z/3$. By the given,\;for $\widetilde{H}_{*}(Z),\,P^{1}_{*}[12+k]\neq0$,   then, for $\widetilde{H}_{*}(Z/S^{4+k})\tg \widetilde{H}_{*}( S^{8+k}\vee S^{12+k})$ \;satisfying $ P^{1}_{*}[12+k]\neq0$, it's impossible.\,This forces \;$f\in\pi_{11+k}(\s^{k}\y)=\z/9\n \s^{k}h\nn$ is of order 9. By Corollary\,\ref{ygd}, we get the result.

\item[\rm(2)]  Immediately got by (1) of this theorem.
\end{itemize}
\end{proof}
\begin{rem} The theorem above would be failed if   $k=4$  was allowed. A counterexample is, after localization at 3, let $u=j_{4}[\iota_{8}, \iota_{8}]\in\pi_{15}(\s^{4}\y)$,  then $C_{u+\Sigma^{4}h}=S^{8}\cup e^{12}\cup e^{16}$ has nontrivial Steenord operations $P^{1}_{\ast}$ of dimension $12$ and $16$,\; by checking $\pi_{15}(-)$, we know
$C_{u+\Sigma^{4}h}$ and $ \Sigma ^{4}\mathbb{H}P^{3}$ are not  homotopy equivalent.
\end{rem}

\subsection{Homotopy decompositions of suspended self smashes }
 
Suppose $\n A_{n}\nn_{n=1}^{\wq} $ is a family of spaces,\;$ \n\mathbbm{f}_{n}: A_{n}\rightarrow  A_{n+1}\nn_{n=1}^{\wq} $ is a family of maps,\,
 we denote the homotopy colimit of the sequence\\  \indent \qquad \xymatrix@C=0.6cm{
 &&& A_{1} \ar[r]^{\mathbbm{f}_{1}} & A_{2} \ar[r]^{\mathbbm{f}_{2}} & A_{3}\ar[r]^{\mathbbm{f}_{3}}  & \cdots   } \textcolor{white}{.}\\
\noindent by\;$\mathop{\mathrm{hocolim}}\limits_{\mathbbm{f}_{n}}A_{n}$. If  $A_{n}=A_{1}$  and\; $\mathbbm{f}_{n}=\mathbbm{f}_{1}$ for all $n\geq1$, then $\mathop{\mathrm{hocolim}}\limits_{\mathbbm{f}_{n}}A_{n}$ is usually denoted by $\mathop{\mathrm{hocolim}}\limits_{\mathbbm{f}_{1}}A_{1}.$

\indent Let  $X$ be a $p$-local path-connected \textit{CW} complex, and   let the symmetric group $\mathcal{S}_{n}$ act on $\emph{X}$$^{\wedge n }$ by permuting positions. Thus, for each $\tau\in\mathcal{S}_{n}$, we have a map $\tau: \emph{X}^{\wedge n }\rightarrow\emph{X}^{\wedge n }$.\;Successively,  by We consider the group structure of $[\Sigma X^{\wedge n },\Sigma X^{\wedge n }$] we obtain a map  $$\hat{k}   : \Sigma\emph{X}^{\wedge n }\rightarrow\Sigma\emph{X}^{\wedge n }$$for any $k$  in the group ring $\mathbb{Z}_{(p)}[\mathcal{S}_{n}]$, by abuse of notation,  $\hat{k}$ and $k$ are both denoted by $k$.
 Then,\, the $\z/p$ coefficient reduced homology $ \widetilde{H}_{*}(\s X^{\wedge n })$ becomes a module over $\mathbb{Z}_{(p)}[\mathcal{S}_{n}]$, the module structure is decided  by permuting factors in graded sense for each $\tau\in\mathcal{S}_{n} $. \,Let  $1=\Sigma_{\alpha}\, e_{\alpha}$ be
an orthogonal decomposition of the identity in $\mathbb{Z}_{(p)}[\mathcal{S}_{n}]$ in terms of primitive idempotents.
 The composition\\
\centerline{ $\s X^{ \wedge n}\stackrel{\mu'}   \longrightarrow            \mathop{\bigvee}\limits_{\af}\s X^{ \wedge n}
\rightarrow            \mathop{\bigvee}\limits_{\af}         \mathop{\mathrm{hocolim}}\limits_{e_{\af}}\s X^{ \wedge n}$}\\
\noindent is a homotopy equivalence, ($\mu'$ is the classical  comultiplication of the suspension as a $co$-$H$ space),  because its induced map on the singular chains over
$\mathbb{Z}_{(p)}$\, is a natural  homotopy equivalence with respect to $X$, (see \cite[section 3.1, p.\,32-34]{Wu} and \cite[p.\,11-12]{Wu1} for  more details).\\
\indent The following is just a special and the most simple  case  of the above by taking $n=2$, $p$  an odd prime and $ 1=\frac{1+(12)}{2}+\frac{1-(12)}{2}$.
\begin{lem}\label{lsd}  Let $p$ be an odd prime, suppose $X$ is a path connected \textit{CW} complex, after localization at $p$,

\begin{itemize}
\item[\rm(1)]\label{fjgs}(\cite{Wu})
there exists the natural decomposition with respect to $X$,\\

\indent\qquad\qquad\qquad\qquad$\s X^{\wedge2}\simeq \mathop{\mathrm{hocolim}}\limits_{ \frac{1+(12)}{2} }\s X^{\wedge2 }\bigvee \mathop{\mathrm{hocolim}}\limits_{\frac{1-(12)}{2}}\s X^{\wedge2}$.
\item[\rm(2)]\label{fjgs}(\cite{Wu})\quad
$  \widetilde{H}_{\ast}( \s X^{\wedge2}  )$ becomes  a module over the group ring $\z_{(p)}[\mathcal{S}_{2}]$,
the module structures are given by:\\ for $(12)\in\mathcal{S}_{2}\xyd\z_{(p)}[\mathcal{S}_{2}]$ which decides a map $ \s X^{\wedge2}\stackrel{(12)\;}\longrightarrow\s X^{\wedge2}$,
 $$(12 )_{*} :   \widetilde{H}_{\ast}( \s X^{\wedge2}  )\rightarrow  \widetilde{H}_{\ast}( \s X^{\wedge2}  ),\; \cg (a\otimes b)\mapsto (-1)^{| a|\cdot|b|}\cg( b\otimes a);$$for $k\in\z_{(p)}\xyd\z_{(p)}[\mathcal{S}_{2}]$ which decides a map $ \s X^{\wedge2}\stackrel{k}\longrightarrow\s X^{\wedge2}$,
 $$k_{*}:   \widetilde{H}_{\ast}( \s X^{\wedge2})\rightarrow  \widetilde{H}_{\ast}( \s X^{\wedge2}  ),\; \cg (a\otimes b)\mapsto k(\cg (a\otimes b)).$$
Here $\cg: \widetilde{H}_{*}(-)\rightarrow \widetilde{H}_{*+1}(\s\,-)$ denotes the suspension isomorphism.
\item[\rm(3)]\label{fjgs}(\cite{Wu1})$$\widetilde{H}_{\ast}(\mathop{\mathrm{hocolim}}\limits_{\frac{1\pm(12)}{2}}\s X^{\wedge2})=
 \mathrm{Im}((    \frac{1\pm(12)}{2}    )_{*} :   \widetilde{H}_{\ast}( \s X^{\wedge2}  )\rightarrow  \widetilde{H}_{\ast}( \s X^{\wedge2}  )  ).  $$

\end{itemize}
\end{lem}
\begin{thm}
After localization at $3$,
\begin{center}
\par
 $\Sigma\mathbb{H}P^{2}\wedge\mathbb{H}P^{2}\,\simeq S^{13}\vee\Sigma^{5}\mathbb{H}P^{3}$.
\end{center}
\end{thm}

\begin{proof}
We denote $a\otimes b$ by  $ab$ for short.\\
After localization at $3$,\;\\
 We know   $\widetilde{H}_{\ast}$($\mathbb{H}P^{2}$)=$\z/3$$\{x,y\}:=V$,\;where  $|x|=4,|y|=8$,\;\;$ P_{*}^{1}(y)=x$.\;\\
  then $\widetilde{H}_{\ast}(\mathbb{H}P^{2}\wedge\mathbb{H}P^{2})= V\otimes V=\z/3\{xx,xy,yx,yy\}$.
By Lemma\,\ref{lsd}, we take $X=\y$, we have $$\Sigma\mathbb{H}P^{2}\wedge\mathbb{H}P^{2}\simeq \mathop{\mathrm{hocolim}}\limits_{\frac{1+(12)}{2}}(\Sigma\mathbb{H}P^{2}\wedge\mathbb{H}P^{2})\vee \mathop{\mathrm{hocolim}}\limits_{\frac{1-(12)}{2}}(\Sigma\mathbb{H}P^{2}\wedge\mathbb{H}P^{2}).$$
 \noindent
 $\widetilde{H}_{\ast}(\mathop{\mathrm{hocolim}}\limits_{\frac{1+(12)}{2}}(\Sigma\mathbb{H}P^{2}\wedge\mathbb{H}P^{2}))
 =\mathrm{Im}(\frac{1+(12)}{2})_{\ast}=\z/3\{            \sigma(xx), \;  \sigma(xy+yx),\;\sigma(yy)\}$,\\
 where\;$$| \sigma(xx)|=9,|\sigma (xy+yx)|=13,|\sigma (yy)|=17,$$
 $$P^{_1}_{\ast}(\sigma(yy))=\sigma(xy+yx),\;P^{_1}_{\ast}(\sigma(xy+yx))=-\sigma(xx);$$
\noindent
$\widetilde{H}_{\ast}(\mathop{\mathrm{hocolim}}\limits_{\frac{1-(12)}{2}}(\Sigma\mathbb{H}P^{2}\wedge\mathbb{H}P^{2}))
 =\mathrm{Im}(\frac{1-(12)}{2})_{\ast}=\z/3\{\sigma(xy-yx)\}$,\; where $$|\sigma(xy-yx)|=13.$$
We know $\pi_{0}(-)=\widetilde{H}_{0}(-;\z)$, and  by Proposition\,\ref{clhp1},   homotopy colimit commutes with $ \pi_{1}(-)$,,\;so\;
$\mathop{\mathrm{hocolim}}\limits_{\frac{1\pm(12)}{2}}(\Sigma\mathbb{H}P^{2}\wedge\mathbb{H}P^{2})$
are simply connected.\;Then,\;$$\mathop{\mathrm{hocolim}}\limits_{\frac{1-(12)}{2}}(\Sigma\mathbb{H}P^{2}\wedge\mathbb{H}P^{2})\simeq S^{13}.$$
By Theorem\,\ref{wjy}(1)\; we have,$$\mathop{\mathrm{hocolim}}\limits_{\frac{1+(12)}{2}}(\Sigma\mathbb{H}P^{2}\wedge\mathbb{H}P^{2})\simeq \s^{5}\h P^{3},$$ hence the result holds.
\end{proof}

\begin{thm} After localization at 3,
\begin{center}
\par
$\Sigma\mathbb{H}P^{3}\wedge\mathbb{H}P^{3}\,\simeq\, \Sigma^{9}\mathbb{H}P^{3}\vee Y$,
\end{center}
where $Y$ is a 6-cell  \textit{CW} complex and $\sk_{13}(Y)=\Sigma^{5}\mathbb{H}P^{2}$.
\end{thm}
\begin{proof}
We denote $a\otimes b$  by  $ab$.\\
After localization at 3,
 $\widetilde{H}_{\ast}(\mathbb{H}P^{3})=\z/3\{x,y,z\}:=V$, $P^{1}_{*}(z)=-y,\;P^{2}_{*}(z)=x,\;P^{1}_{*}(y)=x$,
where $|x|$=4, $|y|$=8 and $|z|$=12.\;\\And $\widetilde{H}_{\ast}(\mathbb{H}P^{3}\wedge\mathbb{H}P^{3})=V\otimes V=\z/3\{xx,xy,xz,yx,yy,yz,zx,zy,zz\}$.\\
By Lemma\,\ref{lsd}, we take $X=\h P^{3}$, we have \\\\
\centerline{$\Sigma\mathbb{H}P^{3}\wedge\mathbb{H}P^{3}\simeq \mathop{\mathrm{hocolim}}\limits_{\frac{1+(12)}{2}}(\Sigma\mathbb{H}P^{3}\wedge\mathbb{H}P^{3})\vee \mathop{\mathrm{hocolim}}\limits_{\frac{1-(12)}{2}}(\Sigma\mathbb{H}P^{3}\wedge\mathbb{H}P^{3})$.}\\\\
By Proposition\,\ref{clhp1},   $\mathrm{hocolim}(-)$ commutes with $ \pi_{1}(-)$ and  $\;\pi_{0}(-)=\widetilde{H}_{0}(-;\z)$\;,\,so\,
$\mathop{\mathrm{hocolim}}\limits_{\frac{1\pm(12)}{2}}(\Sigma\mathbb{H}P^{3}\wedge\mathbb{H}P^{3})$
are simply connected.

$$\widetilde{H}_{\ast}(\mathop{\mathrm{hocolim}}\limits_{\frac{1-(12)}{2}}(\Sigma\mathbb{H}P^{3}\wedge\mathbb{H}P^{3}))=\mathrm{Im}(\frac{1-(12)}{2})_{\ast}=\z/3\{\sigma(xy-yx),\sigma (xz-zx),\sigma (zy-yz)\}$$
where $|\sigma(xy-yx)|=13,|\sigma (xz-zx)|=17,|\sigma (zy-yz)|=21$.\;There are nontrivial   Steenrod operations  $$P^{_1}_{\ast}(\sigma(zy-yz))=-\sigma(xz-zx),\;P^{_1}_{\ast}(\sigma(xz-zx))=\sigma(yx-xy).$$
By  Theorem\,\ref{wjy}(1),\;we get
 $$ \mathop{\mathrm{hocolim}}\limits_{\frac{1-(12)}{2}}(\Sigma\mathbb{H}P^{3}\wedge\mathbb{H}P^{3}) \simeq  \s^{9}\h P^{3}.$$
Let $Y=\mathop{\mathrm{hocolim}}\limits_{\frac{1+(12)}{2}}(\Sigma\mathbb{H}P^{3}\wedge\mathbb{H}P^{3})$,\;then
\, $$\widetilde{H}_{\ast}(Y)=\mathrm{Im}(\frac{1+(12)}{2})_{\ast}=      \z/3\{\sigma(xy+yx),\sigma(xz+zx),\sigma(zy+yz),\sigma(xx),\sigma(yy),\sigma(zz)\}.$$
 Here,\;the degrees of the basis $\sigma(xy+yx),\sigma(xz+zx),\sigma(zy+yz),\sigma(xx),\sigma(yy),\sigma(zz)  $ are 13,\;17,\;21,\;9,\;17,\;25 respectively. We know $P^{_1}_{\ast}(\sigma(xy+yx))=-\sigma(xx)$,\; by Theorem\,\ref{wjy}(1), we can take $\sk_{13}(Y)=\Sigma^{5}\mathbb{H}P^{2}$ up to homotopy. Then the result follows.
\end{proof}

\textcolor{white}{.}\\\\\\\\\

\noindent
Juxin Yang\\School of Mathematical Sciences, Hebei Normal University,\\
Shijiazhuang, 050024, P.R. China,\\
and Yanqi Lake Beijing Institute of Mathematical Sciences and Applications (BIMSA),\;
Beijing, 101408, P.R. China\\
yangjuxin@bimsa.com\\\\
\noindent
Juno Mukai\\
Shinshu University\\
3-1-1 Asahi, Matsumoto, Nagano 390-8621, Japan\\
jmukai@shinshu-u.ac.jp\\\\
\noindent
Jie Wu\\
BIMSA,\\
Beijing, 101408, P.R. China\\
wujie@bimsa.com

\end{document}